\documentclass[]{amsart}
\usepackage[utf8]{inputenc}
\usepackage[T1]{fontenc}
\usepackage{amsmath,amsthm,amsfonts,amssymb,amscd}
\usepackage{tikz-cd}
\usepackage[margin=1.5in]{geometry}
\usepackage{graphicx}
\usepackage{adjustbox}
\usepackage{fancyhdr}
\usepackage{mathrsfs}
\usepackage{xcolor}
\usepackage{enumerate}
\usepackage{listings}
\usepackage{float}
\usepackage{hyperref}
\usepackage{pgffor}
\usepackage{cite}
\usepackage{tikz,xifthen}
\usepackage{pgf}
\usepackage{subcaption}
\usepackage{boldline}
\usepackage{mathtools}
\usepackage{parskip}
\usepackage{tabularray}
\usepackage{lscape}
\usepackage{xfrac}
\UseTblrLibrary{diagbox}
\usepackage{nicematrix}
\usepackage{longtable}

\usetikzlibrary{arrows,positioning,decorations.markings,calc,decorations.pathmorphing}

\tikzset{parallel arrow/.style={latex-,
		shorten >=2mm, shorten <=2mm, 
		decoration={sl,raise=1cm},decorate},shaded/.style={fill=red!10!blue!20!gray!30!white},
	Tbox/.style = {circle, draw, very thick, fill=white, opaque},
	PAdefn/.style = {scale=.7,baseline},
    diagonal fill/.style 2 args={fill=#2, path picture={
\fill[#1, sharp corners] (path picture bounding box.south west) -|
                         (path picture bounding box.north east) -- cycle;}}}

\title{New hyperfinite subfactors with infinite depth}
\author{Dietmar Bisch}
\author{Julio C{\'a}ceres}

\address{Department of Mathematics, Vanderbilt University, 1326 Stevenson Center, Nashville, TN 37240, USA}
\email{dietmar.bisch@vanderbilt.edu}
\email{julio.e.caceres.gonzales@vanderbilt.edu}

\tikzset{
	shaded/.style = {fill=red!10!blue!20!gray!30!white},
	Tbox/.style = {circle, draw, very thick, fill=white, opaque},
	PAdefn/.style = {scale=.7,baseline},
}

\newcommand{\ra}{\rightarrow}
\DeclareMathOperator{\tr}{tr}

\DeclareMathOperator{\ID}{id}

\DeclareMathOperator{\End}{End}

\DeclareMathOperator{\spn}{span}

\DeclareMathOperator{\opp}{opp}

\newcommand{\R}{\mathbb{R}}

\newcommand{\C}{\mathbb{C}}
\newcommand{\Z}{\mathbb{Z}}
\newcommand{\ZZ}{\mathcal{Z}}

\newcommand{\CC}{\mathcal{C}}
\newcommand{\MM}{\mathcal{M}}

\newcommand{\V}{\mathcal{V}}

\newcommand{\sst}{\subset}

\newcommand{\pint}[2]{\langle #1,#2\rangle}

\usepackage{expl3}
\ExplSyntaxOn
\cs_new_eq:NN \Repeat \prg_replicate:nn
\ExplSyntaxOff

\newcommand{\bulmat}[2]{
\begin{smallmatrix}
	\Repeat{#1}{\Repeat{#2}{\bullet}\\}
\end{smallmatrix}
}
\newcommand{\rot}[1]{\rotatebox[origin=c]{90}{$#1$}}

\allowdisplaybreaks

\theoremstyle{definition}
\newtheorem{defi}{Definition}[section]
\newtheorem{rem}[defi]{Remark}

\theoremstyle{plain}
\newtheorem{theo}[defi]{Theorem}
\newtheorem*{theo*}{Theorem}
\newtheorem{prop}[defi]{Proposition}
\newtheorem{lem}[defi]{Lemma}
\newtheorem{coro}[defi]{Corollary}
\newtheorem{algo}[defi]{Algorithm}
\numberwithin{equation}{section}

\begin{document}

\begin{abstract} 
We construct new hyperfinite subfactors with Temperley-Lieb-Jones
(TLJ) standard invariant and Jones indices between $4$ and $3 + \sqrt{5}$.
Our subfactors occur at all indices between $4$ and $5$ at which
finite depth, hyperfinite subfactors exist. The presence of
hyperfinite subfactors with TLJ standard invariant, and hence of
infinite depth, at these special indices comes as a surprise. 
In particular, we obtain the first 
example of a hyperfinite subfactor with ``trivial'' (i.e. TLJ)
standard invariant and integer index $>4$, namely with index $5$.
This is achieved by constructing explicit families of nondegenerate 
commuting squares of multi-matrix algebras based on exotic inclusion 
graphs with appropriate norms. We employ our graph planar algebra
embedding theorem \cite{caceres2023}, a result of Kawahigashi
\cite{kawahigashi2023characterization} and other techniques to 
determine that the resulting commuting square subfactors have indeed 
TLJ standard invariant.

We also show the existence of a hyperfinite subfactor with
trivial standard invariant and index $3 + \sqrt{5}$. It is interesting
that this index features at least two distinct infinite depth, irreducible
hyperfinite subfactors and several finite depth ones as well.
Our work leads to the conjecture that every Jones index $>4$ of an 
irreducible, hyperfinite (finite depth) subfactor can be realized 
by one with TLJ standard invariant.
\end{abstract}

    \maketitle

    \section{Introduction}
Vaughan Jones introduced in his breakthrough paper
\cite{Jones1983} the notion of \textit{index} for an 
inclusion of II$_1$ factors $N \subset M$. He showed
that the index is quantized when $\le 4$ and constructed
irreducible, hyperfinite subfactors at 
each allowed index $4 \cos^2 \pi /n$, $n \ge 3$.
If the index is $<4$, the subfactor is automatically
irreducible, that is $N' \cap M = \C$. 
Constructing {\it irreducible}, hyperfinite subfactors
with Jones index $>4$ turns out to be a very subtle and
hard problem. In fact, Jones asked in \cite{Jones1983} 
to determine the set of indices of \emph{irreducible} subfactors of the
hyperfinite II$_1$ factor, a question that remains open 
to this day. Small, partial progress has been made over the
years. For instance, it was shown in \cite{bisch1994} and later
in \cite{stojanovic2021} that rational, non-integer numbers can be
Jones indices of irreducible, hyperfinite subfactors. Moreover,
it is known that some of these indices are accumulation points of
indices of irreducibe, hyperfinite subfactors (e.g. $4.5$
by \cite{bisch1994} and \cite{schou2013commuting}). However,
it is unknown if the set of Jones indices of irreducible,
hyperfinite subfactors contains an interval, or transcendental
numbers etc.

Note that the situation is quite different if one drops the condition
that the subfactors must be hyperfinite. In \cite{popa1993markov}, Popa
showed that for any $\lambda>4$, one can construct irreducible, 
non-hyperfinite subfactors with index $\lambda$ and Temperley-Lieb-Jones 
standard invariant. Popa's axiomatization of the standard
invariant of a subfactor and the reconstruction theorem proved in
\cite{popa1995axiomatization} led him and Shlyakhtenko to show
that the free group factor $L(F_\infty)$ is universal in the
theory of subfactors \cite{popashlyakhtenko2003}. More precisely,
they showed that every $\lambda>4$ is the index
of a subfactor of $L(F_\infty)$ which itself is isomorphic
to $L(F_\infty)$ and has trivial (i.e. TLJ) standard invariant.
Even better, {\it every} abstract standard invariant can be realized as 
the standard invariant of a subfactor of $L(F_\infty)$.


Indices of finite depth subfactors of the hyperfinite II$_1$ factor
are square norms of their principal graphs, hence must be algebraic 
integers. Thus, in order to construct irreducible, hyperfinite subfactors 
with more interesting Jones indices, one has to turn to infinite depth 
subfactors. 

A classical method to construct hyperfinite subfactors is via
(nondegenerate) commuting squares of multi-matrix algebras
(see for instance \cite{goodman2012coxeter}). A particularly
interesting class of
examples is obtained from so-called {\it spin model commuting
squares} that are built from complex Hadamard matrices. 
These are the easiest commuting squares, and they yield irreducible,
hyperfinite subfactors with integer index. There is at least one
such spin model commuting square subfactor for every non-negative integer, 
and probably many more. It has been a long-standing
problem to determine the standard invariant of these subfactors, and
Jones developed the theory of {\it planar algebras} as a tool to compute 
the ``quantum symmetries'' (or standard invariant) encoded by these 
subfactors \cite{jones1999planar2}, \cite{jones1999planar}. 
Despite partial progress (see e.g. \cite{burstein2015}, 
\cite{montgomery2023}), these subfactors remain mysterious. More
generally, computing the standard invariant of a subfactor that
is obtained from iterating the Jones basic construction of a nondegenerate
commuting square of multi-matrix algebras, called a {\it commuting
square subfactor} in this paper, remains to be a formidable task.
While, in theory, the problem is solved by Ocneanu compactness 
\cite{evanskawahigashi1998} (see also \cite{ocneanu1988quantized}), 
in practice, the computation
surmounts the computational power we have. The problem is hard even in 
the case of low index subfactors, such as the
Haagerup subfactor, which is the finite depth subfactor with 
smallest index $>4$. It has index $\frac{5+\sqrt{13}}{2}$ 
\cite{asaedahaagerup1999}.

In \cite{caceres2023} we use Jones' planar algebra technology to show 
that the standard invariant of a commuting square subfactor has to 
embed, as a planar algebra, into the graph
planar algebra (\cite{jones2000palgbipartite}, \cite{palgnotes})
of one of the inclusion graphs in the commuting square.
Thus, we established an obstruction for the structure of the standard
invariant of such a subfactor
that in particular cases is sufficient to determine the entire
structure of the standard invariant without having to do any
computations whatsoever! More precisely, we proved the following 
embedding theorem \cite{caceres2023} (see section 2 for the notation
used here):

\begin{theo}\label{embeddingT}
		Consider a nondegenerate commuting square of finite 
dimensional C*-algerbas $\{A_{ij}\}_{i,j}$. Let $P_\bullet$ be the 
subfactor planar algebra associated to the commuting square subfactor
$A_{0,\infty}\sst A_{1,\infty}$ and let $G_\bullet$ the graph planar 
algebra associated to the Bratelli diagram of $A_{0,0}\sst A_{1,0}$. 
Then the subfactor planar algebra $P_\bullet$ embeds into $G_\bullet$ 
as a planar algebra.
\end{theo}

In this paper, we use our theorem to construct new irreducible, 
hyperfinite subfactors with TLJ standard invariant at Jones indices
that were only known to be Jones indices of hyperfinite, {\it finite depth}
subfactors. These subfactors have of course $A_\infty$ principal 
graphs, and we will therefore call them $A_\infty$-subfactors. 
The classification
of subfactor planar algebras with small loop parameters 
(see e.g. \cite{popa1994classification}, 
\cite{jones2014classification}, \cite{izumi2015subfactors}, 
\cite{afzaly2015classification}), or, 
equivalently, standard invariants of subfactor with small Jones
indices, is used in a crucial way in our work.


This classification shows that there are only few Jones indices between 
$4$ and $3 + \sqrt{5}$ where finite depth subfactors occur. Moreover, all 
irreducible, infinite depth (hyperfinite) subfactors with indices 
in $\left[4, 5.25\right]$ except these few and $3 +\sqrt{5}$, 
must be $A_\infty$-subfactors.
We display the Jones indices where finite depth subfactors exist
in Table \ref{indices}. Recall that the index of the Extended Haagerup 
subfactor is the largest root of the polynomial $x^3 - 8x^2 + 17x -5$,
which is $\approx 4.37720$. Note that any standard invariant (or,
equivalently, subfactor planar algebra) of a finite depth subfactor,
hyperfinite or not, can be realized as that of a {\it hyperfinite},
finite depth subfactor which arises as a commuting square subfactor
from its canonical commuting square of higher relative commutants
\cite{popafinitedepth}.

\begin{table}[h]
    \centering
    \begin{NiceTabular}{|c|c|c|}
            \CodeBefore
            \rowcolors{1}{blue!15}{}
            \Body
            \hline
            Index & \# of subfactors & Name\\
            $\frac{1}{2}(5+\sqrt{13})$& 2& Haagerup\\
            $\approx 4.37720$& 2&Extended Haagerup\\
            $\frac{1}{2}(5+\sqrt{17})$& 2& Asaeda-Haagerup\\
            $3+\sqrt{3}$& 2& 3311\\
            $\frac{1}{2}(5+\sqrt{21})$& 2& 2221\\
            $5$& 7 & -\\
            $\approx 5.04892$& 2& $\mathfrak{su}(2)_5$ and $\mathfrak{su}(3)_4$\\
            $3+\sqrt{5}$& 11 & -\\
            \hline
        \end{NiceTabular}
    \caption{Indices of (hyperfinite) finite depth subfactors between 
$4$ and $5.25$}
    \label{indices}
\end{table}

An irreducible, hyperfinite subfactor $N\sst M$ with index 
$\frac{5+\sqrt{13}}{2}$ was constructed from a commuting square by 
Haagerup and Schou in \cite[Chapter 7]{schou2013commuting}. It could 
in principle have finite depth and be one of the two Haagerup subfactors. The 
Bratteli diagram for the first vertical inclusion in their commuting square is
\[ \Gamma=\begin{tikzpicture}[scale=0.5]
	\node[circle,fill=black,inner sep=0pt,minimum size=1.5mm] at (1,0) (A2) {};
	\node[circle,fill=black,inner sep=0pt,minimum size=1.5mm] at (2,0) (A3) {};
	\node[circle,fill=black,inner sep=0pt,minimum size=1.5mm] at (3,0) (A4) {};
	\node[circle,fill=black,inner sep=0pt,minimum size=1.5mm] at (4,0) (A5) {};
	\node[circle,fill=black,inner sep=0pt,minimum size=1.5mm] at (5,0) (A6) {};
	\node[circle,fill=black,inner sep=0pt,minimum size=1.5mm] at (5,1) (B1) {};
	\node[circle,fill=black,inner sep=0pt,minimum size=1.5mm] at (5,2) (B2) {};
	\node[circle,fill=black,inner sep=0pt,minimum size=1.5mm] at (6,0) (C1) {};
	\node[circle,fill=black,inner sep=0pt,minimum size=1.5mm] at (7,0) (C2) {};
	\node[circle,fill=black,inner sep=0pt,minimum size=1.5mm] at (8,0) (C3) {};
	\node[circle,fill=black,inner sep=0pt,minimum size=1.5mm] at (9,0) (C4) {};
	\draw (A2) -- (A3) -- (A4) -- (A5) -- (A6) -- (C1) -- (C2) -- (C3) -- (C4);
	\draw (A6) -- (B1) -- (B2);
\end{tikzpicture} \] 
Our embedding theorem \ref{embeddingT} shows that the subfactor planar 
algebra of this commuting square subfactor must embed into the graph planar 
algebra of $\Gamma$. However, Peters proved that the Haagerup subfactor 
planar algebra does {\it not} embed into the graph planar algebra of $\Gamma$ 
in \cite[Theorem 6.8]{peters2010planar}. In fact, corollary 1.4 
in \cite{grossman2018extended} shows that the Haagerup subfactor planar 
algebra only embeds in the graph planar algebra associated to the 
following three graphs:
\begin{align*}
& \begin{tikzpicture}[scale=0.6]
\node[circle,fill=black,inner sep=0pt,minimum size=1.5mm] at (1,0) (A1) {};
\node[circle,fill=black,inner sep=0pt,minimum size=1.5mm] at (2,0) (A2) {};
\node[circle,fill=black,inner sep=0pt,minimum size=1.5mm] at (3,0) (A3) {};
\node[circle,fill=black,inner sep=0pt,minimum size=1.5mm] at (5,0.5) (B1) {};
\node[circle,fill=black,inner sep=0pt,minimum size=1.5mm] at (6,0.5) (B2) {};
\node[circle,fill=black,inner sep=0pt,minimum size=1.5mm] at (7,0.5) (B3) {};
\node[circle,fill=black,inner sep=0pt,minimum size=1.5mm] at (4,0) (P) {};
\node[circle,fill=black,inner sep=0pt,minimum size=1.5mm] at (5,-0.5) (C1) {};
\node[circle,fill=black,inner sep=0pt,minimum size=1.5mm] at (6,-0.5) (C2) {};
\node[circle,fill=black,inner sep=0pt,minimum size=1.5mm] at (7,-0.5) (C3) {};
\draw (A1) -- (A2) -- (A3) -- (P) -- (B1) -- (B2) -- (B3);
\draw (P) -- (C1) -- (C2) -- (C3);
\end{tikzpicture}\\
& \begin{tikzpicture}[scale=0.6]
\node[circle,fill=black,inner sep=0pt,minimum size=1.5mm] at (1,0) (A1) {};
\node[circle,fill=black,inner sep=0pt,minimum size=1.5mm] at (2,0) (A2) {};
\node[circle,fill=black,inner sep=0pt,minimum size=1.5mm] at (3,0) (A3) {};
\node[circle,fill=black,inner sep=0pt,minimum size=1.5mm] at (5,0.5) (B1) {};
\node[circle,fill=black,inner sep=0pt,minimum size=1.5mm] at (4,0) (P) {};
\node[circle,fill=black,inner sep=0pt,minimum size=1.5mm] at (5,-0.5) (C1) {};
\node[circle,fill=black,inner sep=0pt,minimum size=1.5mm] at (6,0) (C2) {};
\node[circle,fill=black,inner sep=0pt,minimum size=1.5mm] at (6,-1) (C3) {};
\draw (A1) -- (A2) -- (A3) -- (P) -- (B1);
\draw (P) -- (C1) -- (C2);
\draw (C1) -- (C3);
\end{tikzpicture}\\
& \begin{tikzpicture}[scale=0.6]
\node[circle,fill=black,inner sep=0pt,minimum size=1.5mm] at (2,0) (A2) {};
\node[circle,fill=black,inner sep=0pt,minimum size=1.5mm] at (3,0) (A3) {};
\node[circle,fill=black,inner sep=0pt,minimum size=1.5mm] at (5,0.5) (B1) {};
\node[circle,fill=black,inner sep=0pt,minimum size=1.5mm] at (4,0) (P) {};
\node[circle,fill=black,inner sep=0pt,minimum size=1.5mm] at (5,-0.5) (C1) {};
\node[circle,fill=black,inner sep=0pt,minimum size=1.5mm] at (5,0) (C2) {};
\draw (A2) -- (A3) -- (P) -- (B1);
\draw (P) -- (C1);
\draw (P) -- (C2);
\end{tikzpicture}
\end{align*}

The first two graphs are the principal graphs of the Haagerup 
subfactor. Thus, we can conclude that the commuting square subfactor 
$N\sst M$ constructed by Haagerup and Schou cannot be the Haagerup 
subfactor and therefore, by the classification of small index subfactor 
planar algebras, must be an $A_\infty$-subfactor, i.e. has trivial
standard invariant.

This is not the only index where this technique can be applied. It follows
from theorem 1.3 in \cite{grossman2018extended} that the Extended 
Haagerup subfactor planar algebra embeds in the graph planar algebras 
of the following four graphs:

\begin{table}[h!]
\centering
\begin{tblr}[caption={Extended Haagerup module graphs}]{
colspec={lrc},
vline{3},
hline{2},
row{1,2}={abovesep=5pt,belowsep=5pt}
}
Principal graph&\begin{tikzpicture}[scale=0.5,baseline=0]
\node[circle,fill=black,inner sep=0pt,minimum size=1.5mm] at (-1,0) (D1) {};
\node[circle,fill=black,inner sep=0pt,minimum size=1.5mm] at (-2,0) (D2) {};
\node[circle,fill=black,inner sep=0pt,minimum size=1.5mm] at (-3,0) (D3) {};
\node[circle,fill=black,inner sep=0pt,minimum size=1.5mm] at (0,0) (A0) {};
\node[circle,fill=black,inner sep=0pt,minimum size=1.5mm] at (1,0) (A1) {};
\node[circle,fill=black,inner sep=0pt,minimum size=1.5mm] at (2,0) (A2) {};
\node[circle,fill=black,inner sep=0pt,minimum size=1.5mm] at (3,0) (A3) {};
\node[circle,fill=black,inner sep=0pt,minimum size=1.5mm] at (5,0.5) (B1) {};
\node[circle,fill=black,inner sep=0pt,minimum size=1.5mm] at (6,0.5) (B2) {};
\node[circle,fill=black,inner sep=0pt,minimum size=1.5mm] at (7,0.5) (B3) {};
\node[circle,fill=black,inner sep=0pt,minimum size=1.5mm] at (4,0) (P) {};
\node[circle,fill=black,inner sep=0pt,minimum size=1.5mm] at (5,-0.5) (C1) {};
\node[circle,fill=black,inner sep=0pt,minimum size=1.5mm] at (6,-0.5) (C2) {};
\node[circle,fill=black,inner sep=0pt,minimum size=1.5mm] at (7,-0.5) (C3) {};
\draw (D3) -- (D2) -- (D1) -- (A0) -- (A1) -- (A2) -- (A3) -- (P) -- (B1) -- (B2) -- (B3);
\draw (P) -- (C1) -- (C2) -- (C3);
\end{tikzpicture} &\begin{tikzpicture}[scale=0.5,baseline=0]
\node[circle,fill=black,inner sep=0pt,minimum size=1.5mm] at (1,1) (A1) {};
\node[circle,fill=black,inner sep=0pt,minimum size=1.5mm] at (2,0.5) (A2) {};
\node[circle,fill=black,inner sep=0pt,minimum size=1.5mm] at (1,-1) (B1) {};
\node[circle,fill=black,inner sep=0pt,minimum size=1.5mm] at (2,-0.5) (B2) {};
\node[circle,fill=black,inner sep=0pt,minimum size=1.5mm] at (3,0) (C1) {};
\node[circle,fill=black,inner sep=0pt,minimum size=1.5mm] at (4,0) (C2) {};
\node[circle,fill=black,inner sep=0pt,minimum size=1.5mm] at (5,0) (C3) {};
\node[circle,fill=black,inner sep=0pt,minimum size=1.5mm] at (6,0) (C4) {};
\node[circle,fill=black,inner sep=0pt,minimum size=1.5mm] at (7,0.5) (D1) {};
\node[circle,fill=black,inner sep=0pt,minimum size=1.5mm] at (7,-0.5) (D2) {};
\node[circle,fill=black,inner sep=0pt,minimum size=1.5mm] at (8,-1) (D3) {};
\draw (A1) -- (A2) -- (C1) -- (C2) -- (C3) -- (C4) -- (D1);
\draw (B1) -- (B2) -- (C1);
\draw (C4) -- (D2) -- (D3);
\end{tikzpicture}\\
Dual principal graph&\begin{tikzpicture}[scale=0.5,baseline=0]
\node[circle,fill=black,inner sep=0pt,minimum size=1.5mm] at (-1,0) (D1) {};
\node[circle,fill=black,inner sep=0pt,minimum size=1.5mm] at (-2,0) (D2) {};
\node[circle,fill=black,inner sep=0pt,minimum size=1.5mm] at (-3,0) (D3) {};
\node[circle,fill=black,inner sep=0pt,minimum size=1.5mm] at (0,0) (A0) {};
\node[circle,fill=black,inner sep=0pt,minimum size=1.5mm] at (1,0) (A1) {};
\node[circle,fill=black,inner sep=0pt,minimum size=1.5mm] at (2,0) (A2) {};
\node[circle,fill=black,inner sep=0pt,minimum size=1.5mm] at (3,0) (A3) {};
\node[circle,fill=black,inner sep=0pt,minimum size=1.5mm] at (5,0.5) (B1) {};
\node[circle,fill=black,inner sep=0pt,minimum size=1.5mm] at (4,0) (P) {};
\node[circle,fill=black,inner sep=0pt,minimum size=1.5mm] at (5,-0.5) (C1) {};
\node[circle,fill=black,inner sep=0pt,minimum size=1.5mm] at (6,0) (C2) {};
\node[circle,fill=black,inner sep=0pt,minimum size=1.5mm] at (6,-1) (C3) {};
\draw (D3) -- (D2) -- (D1) -- (A0) -- (A1) -- (A2) -- (A3) -- (P) -- (B1);
\draw (P) -- (C1) -- (C2);
\draw (C1) -- (C3);
\end{tikzpicture}&\begin{tikzpicture}[scale=0.5,baseline=0]
\node[circle,fill=black,inner sep=0pt,minimum size=1.5mm] at (2,0.5) (A2) {};
\node[circle,fill=black,inner sep=0pt,minimum size=1.5mm] at (2,-0.5) (B2) {};
\node[circle,fill=black,inner sep=0pt,minimum size=1.5mm] at (3,0) (C1) {};
\node[circle,fill=black,inner sep=0pt,minimum size=1.5mm] at (4,0) (C2) {};
\node[circle,fill=black,inner sep=0pt,minimum size=1.5mm] at (5,0) (C3) {};
\node[circle,fill=black,inner sep=0pt,minimum size=1.5mm] at (6,0) (C4) {};
\node[circle,fill=black,inner sep=0pt,minimum size=1.5mm] at (7,0) (C5) {};
\node[circle,fill=black,inner sep=0pt,minimum size=1.5mm] at (8,0) (C6) {};
\node[circle,fill=black,inner sep=0pt,minimum size=1.5mm] at (9,0.5) (D1) {};
\node[circle,fill=black,inner sep=0pt,minimum size=1.5mm] at (9,-0.5) (D2) {};
\node[circle,fill=black,inner sep=0pt,minimum size=1.5mm] at (6,-1) (E1) {};
\draw (A2) -- (C1) -- (C2) -- (C3) -- (C4) -- (C5) -- (C6) -- (D1);
\draw (B2) -- (C1);
\draw (C6) -- (D2);
\draw (C4) -- (E1);
\end{tikzpicture}
\end{tblr}
\caption{Extended Haagerup module graphs}
\label{exthaagMGs}
\end{table}


More generally, theorem 1.2 in \cite{grossman2018extended} 
characterizes those finite connected bipartite graphs into
whose graph planar algebra a given finite depth subfactor planar
algebra embeds as certain fusion graphs for the identity
object in the box space $P_{1,+}$ of the planar algebra.
We explain this in more detail in section $3$ and refer
to these graphs as {\it module graphs}.

This theorem implies that if we can find a subfactor whose subfactor 
planar algebra embeds into the graph planar algebra of a graph that 
is {\it not} a module graph for any finite depth subfactor with the 
same index, then this subfactor must have infinite depth. We use this 
idea to construct irreducible, hyperfinite subfactors with trivial 
standard invariant and indices $\sim 4.37720$ and $\frac{5+\sqrt{17}}{2}$. 
To this end, in section \ref{newcs}, we construct commuting squares whose 
Bratteli diagram of the first vertical inclusion has norm $\sim 4.37720$ 
(respectively $\frac{5+\sqrt{17}}{2}$) and is not one of the module graphs 
for the Extended Haagerup (respectively Asaeda-Haagerup) subfactor. In 
the case of the Asaeda-Haagerup subfactor, in section \ref{modgraphsembed}, 
we use the results from \cite{grossman2018asaeda} and 
\cite{grossman2016brauer} together with the techniques outlined 
in \cite[section 4]{grossman2018extended} to compute a finite list of 
all potential module graphs.

To handle the remaining indices from Table \ref{indices} we use a different approach. In section \ref{1param} we construct a one-parameter family of non-equivalent commuting squares associated to 4-stars of the form $G=S(i,i,j,j)$. Kawahigashi in \cite{kawahigashi2023characterization} characterized all commuting squares that produce a fixed finite depth subfactor. This characterization shows that a finite depth subfactor can only be realized as a commuting square subfactor in countably many ways. Thus, by a cardinality argument, we can conclude that 
our one-parameter families 
of commuting squares must produce infinite depth subfactors with 
indices $\frac{5+\sqrt{17}}{2},3+\sqrt{3},\frac{5+\sqrt{21}}{2}$, $5$ 
and $3+\sqrt{5}$. The classification of small index subfactors then shows 
that all but the one with index $3+\sqrt{5}$ must have trivial standard 
invariant. 

It was shown in \cite{bisch1997intermediate} that there is a
standard invariant with infinite principal graphs and index
$3+\sqrt{5}$. It is a Fuss-Catalan algebra of Bisch and Jones
and arises as the free composition of an $A_3$ and an $A_4$ 
standard invariant. By a result of Popa \cite{popa1995freeindependent},
this Fuss-Catalan standard invariant can be realized as the
standard invariant of a hyperfinite (necessarily irreducible and
infinite depth) subfactor with index $3+\sqrt{5}$. This subfactor 
has an intermediate subfactor. In Section 6, we show that a commuting 
square subfactor that has an intermediate subfactor is always obtained 
from an intermediate commuting square, which might have non-connected
inclusion graphs.
We then show that the commuting square subfactors with index $3+\sqrt{5}$ 
coming from our explicit one-parameter family of commuting squares do not 
admit intermediate commuting squares. This implies that the infinite
depth subfactor obtained at index $3+\sqrt{5}$ does not have an intermediate 
subfactor and hence, by classification, must be an $A_\infty$-subfactor.

The results obtained in this paper lead to the conjecture that for every
Jones index $>4$ that is the index of an irreducible, finite depth subfactor, 
there exists also a hyperfinite A$_\infty$-subfactor with that index. 
It may even be true that every number $>4$ that is in Jones' set of indices 
of irreducible,
hyperfinite subfactors can be realized as the index of a hyperfinite, 
infinite depth subfactor with trivial standard invariant, i.e. a hyperfinite
A$_\infty$-subfactor. The classification of small index standard invariants,
combined with the results of this paper, establish these conjectures for
all indices in the interval $[4,5]$ and for index $3 + \sqrt{5}$. If true,
it would mean that progress on Jones' 40-year old index problem could 
come from focussing efforts on hyperfinite $A_\infty$-subfactors.

{\bf Acknowledgements.} D.B and J.C. were partially supported by US ARO 
grant W911NF-23-1-0026 during this project. The material of section 6 is 
based on work that was also supported by 
the National Science Foundation under Grant No. DMS-1928930, while D.B. 
was in residence at the Simons Laufer Mathematical Sciences Institute in 
Berkeley, California, during July 2024. Section 5 is based on discussions
with Uffe Haagerup in the early 1990's during D.B.'s frequent visits to 
Odense University (now the University of Southern Denmark). D.B. would like 
to credit Uffe Haagerup for teaching him how to compute bi-unitary 
connections and for sharing his insights pertaining to commuting squares 
associated to $4$-stars. This led D.B. to compute a one-parameter family of
connections for $S(1,1,3,3)$ back in 1992.

    \section{Subfactors and commuting squares}    
For the convenience of the reader, we briefly recall in this section 
the construction of hyperfinite, finite index subfactors from 
nondegenerate commuting squares of multi-matrix
algebras (for details, see \cite{goodman2012coxeter}, 
\cite{evanskawahigashi1998}, \cite{jones1997introduction}). 
Commuting squares are a useful tool to compute invariants of
a subfactor $N\sst M$ when it is approximated by finite-dimensional 
C*-algebras in the following way
	\[ \begin{array}{ccc}
	N&\sst &M\\
	\rotatebox[origin=c]{90}{$\sst$}& &\rotatebox[origin=c]{90}{$\sst$}\\
	A_{n+1}&\sst &B_{n+1}\\
	\rotatebox[origin=c]{90}{$\sst$}& &\rotatebox[origin=c]{90}{$\sst$}\\
	A_{n}&\sst &B_{n}
	\end{array} \]
	where $N=\left(\bigcup_{n\geq 0} A_n\right)''$ and $M=\left(\bigcup_{n\geq 0} B_n\right)''$ \cite{pimsner1986entropy}. 

	
	\begin{defi}
		Let $A_0\sst B_0$, $A_1\sst B_1$ be finite von Neumann algebras such that
		\begin{equation}\label{csquare}
		\begin{array}{ccc}
		A_{1}&\sst &B_{1}\\
		\rotatebox[origin=c]{90}{$\sst$}& &\rotatebox[origin=c]{90}{$\sst$}\\
		A_{0}&\sst &B_{0}
		\end{array}
		\end{equation}
		and let $\tr$ be a faithful normal trace in $B_1$, which
induces traces on the subalgebras by restriction. 
Then (\ref{csquare}) is called a \emph{commuting square} if 
$E_{B_0}E_{A_1}=E_{A_0}$ where $E_{B_0},E_{A_1},E_{A_0}$ are the 
unique tr-preserving conditional expectations.
	\end{defi}

	\begin{prop}
		Consider a commuting square as in (\ref{csquare}) of finite-dimensional C*-algebras  and let $\tr$ be a Markov trace for $B_0\sst B_1$. Let $B_2=\langle B_1,e_{B_0}\rangle$ be the basic construction and let $A_2=\{A_1,e_{B_0}\}''$. Then
		\begin{equation*}
		\begin{array}{ccc}
		A_{2}&\sst &B_{2}\\
		\rotatebox[origin=c]{90}{$\sst$}& &\rotatebox[origin=c]{90}{$\sst$}\\
		A_{1}&\sst &B_{1}
		\end{array}
		\end{equation*}
		is also a commuting square.
	\end{prop}
	
	\begin{defi}
		Consider a commuting square
		\begin{equation}\label{csfindim}
		\begin{array}{ccc}
		A_{1,0}&\stackrel{\text{K}}{\sst} &A_{1,1}\\
		\rotatebox[origin=c]{90}{$\sst$}_G& &\rotatebox[origin=c]{90}{$\sst$}_L\\
		A_{0,0}&\stackrel{\text{H}}{\sst} &A_{0,1}
		\end{array}
		\end{equation}
		of finite-dimensional C*-algebras with a faithful normal trace on $A_{1,1}$ and inclusion graphs $G$, $H$, $K$ and $L$. We say that (\ref{csfindim}) is a \emph{nondegenerate} or \emph{symmetric} commuting square if $GH=KL$ and $HL^t=G^tK$.
	\end{defi}
	
It is easy to see that whenever we have a {\it nondegenerate} commuting square, 
the inclusion $A_0\sst A_1\sst A_2$ from the previous proposition is 
isomorphic to the basic construction 
$A_0\sst A_1\sst \langle A_1, e_{A_0}\rangle$. 
	 
	\begin{prop}
		Suppose (\ref{csfindim}) is a nondegenerate commuting square 
with connected graphs $G$, $H$, $K$ and $L$. 
Then $\|H\|=\|K\|$ and $\|G\|=\|L\|$. Moreover, $\tr$ is a Markov trace 
for $A_{1,0}\sst A_{1,1}$, $\tr|_{A_{0,1}}$ is a Markov trace for 
$A_{0,0}\sst A_{0,1}$, $\tr$ is a Markov trace for 
$A_{0,1}\sst A_{1,1}$ and $\tr|_{A_{1,0}}$ is a Markov trace 
for $A_{0,0}\sst A_{1,0}$.
	\end{prop} 
	
	
Starting from a nondegenerate commuting square and a Markov trace, 
we can iterate Jones' basic construction in the horizontal and vertical 
directions. Thus, we obtain a lattice of inclusions of finite-dimensional
C$^*$-algebras of the form

	\begin{equation}\label{infgrid}
	\begin{array}{ccccccccc}
	A_{\infty,0} & \sst & A_{\infty,1} & \sst & A_{\infty,2} & \sst & \cdots & \sst & A_{\infty,\infty} \\
	\cup &  & \cup &  & \cup &  &  &  & \cup \\
	\vdots &  & \vdots &  & \vdots &  &  &  & \vdots \\
	\cup &  & \cup &  & \cup &  &  &  & \cup \\
	A_{2,0} & \sst & A_{2,1} & \sst & A_{2,2} & \sst & \cdots & \sst & A_{2,\infty}  \\
	\cup &  & \cup &  & \cup &  &  &  &  \\
	A_{1,0}& \sst & A_{1,1} & \sst & A_{1,2} & \sst & \cdots & \sst & A_{1,\infty} \\
	\cup &  & \cup &  & \cup &  &  &  & \cup \\
	A_{0,0} & \sst & A_{0,1} & \sst & A_{0,2} & \sst & \cdots & \sst & A_{0,\infty} 
	\end{array}
	\end{equation}

where $A_{n,\infty}=\left(\bigcup_{k} A_{n,k}\right)''$ and 
$A_{\infty,k}=\left(\bigcup_{n} A_{n,k}\right)''$ are hyperfinite 
von Neumann algebras obtained from the GNS construction applied to the
extension of Markov trace in the well-known way. If all inclusion 
matrices are connected, $A_{n,\infty}$ and $A_{\infty,k}$ are factors. 
If the initial commuting square is as in (2.2) with connected
inclusion matrices, the indices are
$[A_{1,\infty}:A_{0,\infty}]=\|G\|^2$ and
$[A_{\infty, 1}:A_{\infty, 0}]=\|H\|^2$. We call the subfactors
$A_{0,\infty} \subset A_{1,\infty}$ and
$A_{\infty, 0} \subset A_{\infty, 1}$ {\it commuting square subfactors}.
We will focus on the horizontal subfactor
$A_{0,\infty} \subset A_{1,\infty}$ as it is obvious how to 
restate the results for the vertical subfactor as well.

The next result shows how to compute, in principle, the standard invariants
of commuting square subfactors.
	
	\begin{theo}[Ocneanu compacntess]\label{ocneanu}
		Given a lattice of inclusions arising from a nondegenerate commuting square as in (\ref{infgrid}), we have
		\[ A'_{0,\infty}\cap A_{n,\infty}= A'_{0,1}\cap A_{n,0},\quad A'_{\infty,0}\cap A_{\infty,n}= A'_{1,0}\cap A_{0,n},\text{ for all }n\geq 0.\]
	\end{theo}


Schou proved in \cite{schou2013commuting} an easy criterion that will 
ensure irreducibilty of commuting square subfactors. It is a 
generalization of Wenzl's criterion \cite[Theorem 1.6]{wenzl1988hecke}. 

    \begin{theo}[Wenzl's criterion]\label{wenzcrit} 
        $\dim(A_{0,\infty}'\cap A_{1,\infty})\leq (\min\{\text{1-norm of rows of $G$ and $L$}\})^2$.
    \end{theo}

    Hence, if either $G$ or $L$ have at least one vertex of degree 1 then $\dim(A_{0,\infty}'\cap A_{1,\infty})=1$, that is, $A_{0,\infty}\sst A_{1,\infty}$ is irreducible.

    \section{Module graphs and embeddings}\label{modgraphsembed}

Consider $\CC$ a fusion category and $\{\CC_i\}_i$ its Morita equivalence class. 
We want to determine all possible left $\CC$-module $C^*$-categories $\mathcal{M}$ for a fixed unitary multifusion category $\CC$ coming from a subfactor. In our case, $\CC$ will be such that $1=1_+\oplus 1_-$ and consequently we can write
$$\CC=\begin{pNiceMatrix}
\CC_{11}&\CC_{12}\\
\CC_{21}&\CC_{22}
\end{pNiceMatrix}$$

where $\CC_{ij}=1_\pm\otimes \CC\otimes 1_\pm$. In fact, for $\CC$ coming from an irreducible finite depth subfactor $\CC_{12}$ will be a Morita equivalence between the fusion categories $\CC_{11}$ and $\CC_{22}$. Now, if $\MM$ is an indecomposable left $\CC$-module category then $\MM_1:=1_+\otimes \MM$ and $\MM_2:=1_-\otimes \MM$ will be indecomposable left $\CC_{11}$-module and $\CC_{22}$-module categories respectively. Following theorem 3.3 in 
\cite{grossman2018extended}, we have:
\begin{theo}[\cite{etingof2010fusion,etingof2016tensor}]
If $\mathcal{A}$ is a fusion category and $\mathcal{N}$ is a semisimple $\mathcal{A}$-module category, then $\mathcal{A}-\mathcal{B}$ bimodule category structures on $\mathcal{N}$ which extend the $\mathcal{A}$-module structure correspond exactly to functors $\mathcal{F}:\mathcal{A}\ra \End_\mathcal{A}(\mathcal{N})$, and such a bimodule is a Morita equivalence if and only if $\mathcal{F}$ is an equivalence of multitensor categories. Two such bimodule categories are equivalent if and only if the functors differ by an inner autoequivalence. Furthermore, $\End_\mathcal{A}(\mathcal{N})$ is a tensor category (with simple unit object) if and only if $\mathcal{N}$ is indecomposable.
\end{theo}
This implies there is a bijection between indecomposable left $\mathcal{A}$-modules and Morita equivalences between $\mathcal{A}$ and fusion categories $\mathcal{B}$. Therefore, to determine all possible indecomposable left $\CC$-module categories $\MM$, it is enough to understand the Morita equivalence class of $\CC_{11}$ (or $\CC_{22}$). This information is captured by the Brauer-Picard 
groupoid \cite{etingof2010fusion}.

For the Extended Haagerup fusion categories, the Brauer-Picard groupoid is 
computed in \cite{grossman2018extended}. Moreover, they show the following
theorem that is crucial for the next section:
\begin{theo}\label{fusionembed}
        Suppose $P_\bullet$ is a finite depth subfactor planar algebra. Let $\mathcal{C}$ denote the unitary multifusion category of projections in $P_\bullet$, with distinguished object $X=\ID_{1,+}\in P_{1,+}$, and the standard unitary pivotal structure with respect to $X$. There is an equivalence between:
        \begin{enumerate}
            \item Planar algebra embeddings $P_\bullet \ra G_\bullet$, where $G_\bullet$ is the graph planar algebra associated to a finite connected bipartite graph $\Gamma$, and
            \item indecomposable finitely semisimple pivotal left $\mathcal{C}$-module $C^*$ categories $\mathcal{M}$ whose fusion graph with respect to $X$ is $\Gamma$.
        \end{enumerate}
    \end{theo}
We will refer to the fusion graphs of $\mathcal{M}$ with respect to $X$ 
as \emph{module graphs}. The combinatorial information for all these 
Morita equivalences is then used to compute the four module graphs, 
depicted in Table \ref{exthaagMGs}, associated to the Extended Haagerup 
subfactor. We need to determine the module graphs for the Asaeda-Haagerup
subfactor. In fact, it suffices to compute a finite list that contains 
all of them. 

\subsection{Asaeda-Haagerup module graphs}

Let $\mathcal{AH}_1$ and $\mathcal{AH}_2$ be fusion categories corresponding to the even and odd part of the Asaeda-Haagerup subfactor. In \cite{grossman2016brauer}, the authors identify a fusion category $\mathcal{AH}_3$ which is also in the Morita equivalence class of $\mathcal{AH}_1$ and $\mathcal{AH}_3$. Moreover, they determine all the Morita equivalences between them. Denote the fusion rings of $\mathcal{AH}_i$ by $AH_i$.
\begin{theo}[\cite{grossman2016brauer}]
    There are exactly $4$ invertible bimodule categories over each not-necessarily-distinct pair $\mathcal{AH}_i$-$\mathcal{AH}_j$, up to equivalence. These realize the following fusion bimodules, which are each realized uniquely:
    \begin{itemize}
        \item $10_{11},12_{11},13_{11},14_{11},2_{12},5_{12},8_{12},9_{12},2_{13},3_{13},6_{13},7_{13},$
        \item $2_{21},5_{21},8_{21},9_{21},8_{22},11_{22},12_{22},13_{22},1_{23},3_{23},4_{23},6_{23},$
        \item $2_{31},3_{31},6_{31},7_{31},1_{32},2_{32},4_{32},6_{32},8_{33},11_{33},12_{33},13_{33}.$
    \end{itemize}
\end{theo}
Here $a_{ij}$ denotes the $a^{th}$ fusion bimodule on a list of potential $AH_i$-$AH_j$ fusion bimodules. Using this notation, $9_{12}$ is the fusion bimodule realized by the Asaeda-Haagerup subfactor. It is generated by a single element $X$ which corresponds to $\ID_{1,+}\in P_{1,+}$ in the subfactor planar algebra. 

\begin{defi}[\cite{grossman2016brauer}]
A \emph{multiplication map} on a triple of fusion bimodules 
\newline
$(_{A}K_{B},_{B}L_C,_{A}M_{C})$ is a homomorphism from $_{A}K\otimes_{B}L_C$ to $_{A}M_C$ which takes tensor products of basis elements in $K$ and $L$ to non-negative combinations of basis elements of $M$ and preserves dimension and multiplication by duals. The triple $(_{A}K_{B},_{B}L_C,_{A}M_{C})$ is said to be \emph{multiplicatively compatible} if there exists such a multiplication map.
\end{defi}

The authors, in a supplementary file, provide a list of multiplicative compatible bimodules, in particular:
\begin{align*}
9_{12}\cdot 2_{21}&=10_{11},& 9_{12}\cdot 5_{21}&=13_{11},& 9_{12}\cdot 8_{21}&=12_{11},& 9_{12}\cdot 9_{21}&=14_{11}\\
9_{12}\cdot 8_{22}&=2_{12},& 9_{12}\cdot 11_{22}&=5_{12},& 9_{12}\cdot 12_{22}&=8_{12},& 9_{12}\cdot 13_{22}&=9_{12}\\
9_{12}\cdot 1_{23}&=2_{13},& 9_{12}\cdot 3_{23}&=6_{13},& 9_{12}\cdot 4_{23}&=3_{13},& 9_{12}\cdot 6_{23}&=7_{13}.
\end{align*}
To determine all the module graphs with respect to $X$ we need to determine the multiplication maps for the triples above. 

Let $(_{A}K_{B},_{B}L_C,_{A}M_{C})$ be a triple of fusion bimodules, with bases $\xi_i$ for $1\leq i\leq l$, $\eta_j$ for $1\leq j\leq m$, and $\mu_k$ for $1\leq k\leq n$, respectively. We will consider the lexicographic order for ordered pairs $(p,q)$, $1\leq p\leq l$, $1\leq q\leq m$ and denote by $(p,q)'$ its successor in this order.
\begin{defi}\label{partialmult}
    A $(p,q)$-partial multiplication map for $(p,q)\leq (l,m)$ is an assignment of a vector of integers of length $n$, $v^{ij}$, for each pair $(i,j)\leq (p,q)$ such that
    \begin{enumerate}[(a)]
        \item $d(\xi_i)d(\eta_j)=\sum_{k=1}^n v^{ij}_kd(\mu_k)$,
        \item $(\bar{\xi_{i_1}} \xi_{i_2},\eta_{j_1}\bar{ \eta_{j_2}})=\sum_{k=1}^n v^{i_1 j_1}_k v^{i_2 j_2}_k$,
        \item for all basis elements $\lambda\in A$ and $(i_1,j_1),(i_2,j_2)\leq (p,q)$ we have
        $$((\xi_{i_1}(\eta_{j_1}\bar{\eta}_{j_2})\bar{\xi}_{i_2}),\lambda)=\sum v^{i_1 j_1}_{k_1}v^{i_2 j_2}_{k_2}(\mu_{k_1}\bar{\mu}_{k_2},\lambda),$$
        \item[(c')] for all basis elements $\kappa\in C$ and $(i_1,j_1),(i_2,j_2)\leq (p,q)$ we have
        $$((\bar{\xi}_{i_1}(\bar{\eta}_{j_1}\eta_{j_2})\xi_{i_2}),\lambda)=\sum v^{i_1 j_1}_{k_1}v^{i_2 j_2}_{k_2}(\bar{\mu}_{k_1}\mu_{k_2},\kappa),$$
        \item for all $(i_1,j_1),(i_2,j_2)\leq (p,q)$ we have
        $$v^{i_1 j_1}\cdot v^{i_2 j_2}=(\bar{\xi}_{i_1}\xi_{i_2},\eta_{j_2}\bar{\eta}_{j_1}).$$
    \end{enumerate}
\end{defi}
We define a $(0,0)$-partial multiplication map to be the empty map, and let $(0,0)'=(1,1)$. Using Algorithm 5.2 in \cite{grossman2016brauer} we inductively build all $(p,q)$-partial multiplication maps for each triple, in particular, the $(l,m)$-partial multiplication map will be the multiplication map for the triple. 
\begin{algo}[\cite{grossman2016brauer}]\label{algorithm} Step 1: Start with the $(0,0)$-partial multiplication map. Then inductively find all extensions of a given $(p,q)$-partial multiplication map, $(p,q)<(l,m)$, to a $(p,q)'$-partial multiplication map as follows:
\begin{enumerate}
    \item[1a)] Let $(p,q)'=(p',q')$. Find candidates for $v^{p' q'}$ by checking conditions (a) and (b) above as follows: for each sum of squares decomposition $(a_i,b_i)$ of $(\bar{\xi}_{p'}\xi_{p'},\eta_{q'}\bar{\eta}_{q'})$ such that $\sum b_i\leq n$, form the vector $v$ of size $n$ given by $b_i$ copies of each $a_i$ with the rest of the entries equal to $0$. Then find all distinct vectors $v'$ which arise as permutations of $v$.
    \item[1b)] For each candidate $v'$ for $v^{p'q'}$ found in (1a), check whether $v'\cdot d_M=d(\xi_i)d(\eta_j)$, where $d_M$ is the dimension vector of the bimodule $M$. Finally, if the dimension condition is satisfied for $v'$, check conditions (c), (c'), and (d) for $(i_1,j_1)=(p',q')$ and all $(i_2,j_2)\leq (p',q')$, using $v^{p'q'}=v'$.

    Step 2: For each $(l,m)$-partial multiplication map found in Step 1, check whether 
    \begin{align*}
        (\xi \rho)\eta&= \xi(\rho \eta),& \lambda(\xi \eta)&=(\lambda\xi)\eta, &  (\xi\eta)\kappa&=\xi(\eta\kappa),
    \end{align*}
    for all basis vectors
    $$\lambda\in A,\; \rho\in B,\; \kappa\in C,\; \xi\in K,\; \eta\in L,$$
    where multiplication between elements of $K$ and $L$ is defined on basis elements by the partial multiplication map and extended biadditively.
\end{enumerate}
 
\end{algo}

Using the combinatorial data for the fusion modules and the algorithm above we compute the fusion graph of each fusion module with respect to $X$:
\begin{table}[H]
\centering
\begin{tblr}{
colspec={rlrl},
hline{2-5} = {1-2}{leftpos = -1, rightpos = -1, endpos},
hline{2-5} = {3-4}{leftpos = -1, rightpos = -1, endpos},
hline{6,7} = {1-4}{leftpos = -1, rightpos = -1, endpos},
vline{3},
row{2-7}={abovesep=5pt,belowsep=5pt},
cell{6,7}{2}={c=2}{l},
}
$9_{12}\cdot 8_{22}=2_{12}$&\begin{tikzpicture}[main/.style={circle,fill=black,inner sep=0pt,minimum size=1.5mm},baseline=0,scale=0.6]
			\node[main] at (0,0) (A0) {};
			\node[main] at (-1,0.5) (A1) {};
			\node[main] at (-1,-0.5) (B1) {};
			\node[main]  at (-2,0.5) (A2) {};
			\node[main]  at (-2,-0.5) (B2) {};
			\node[main]  at (1,0.5) (A3) {};
			\node[main]  at (1,-0.5) (B3) {};
			\draw (A2) -- (A1) -- (A0) -- (A3);
			\draw (B2) -- (B1) -- (A0) -- (B3);
			\end{tikzpicture} & $9_{12}\cdot3_{23}=6_{13}$ &\begin{tikzpicture}[main/.style={circle,fill=black,inner sep=0pt,minimum size=1.5mm},baseline=0,scale=0.6]
			\node[main] at (1,-0.5) (A1) {};
			\node[main] at (2,0) (A2) {};
			\node[main] at (3,0) (A3) {};
			\node[main] at (1,0.5) (B2) {};
			\node[main] at (4,-0.5) (B4) {};
			\node[main] at (4,0.5) (A4) {};
			\node[main] at (5,-0.5) (B5) {};
			\node[main] at (5,0.5) (A5) {};
			\node[main] at (6,-0.5) (B6) {};
			\node[main] at (6,0.5) (A6) {};
			\draw (B2) -- (A2) -- (A3) -- (B4) -- (B5) -- (B6);
			\draw (A1) --  (A2)-- (A3) -- (A4) -- (A5) -- (A6);
			\end{tikzpicture}\\ 
$9_{12}\cdot 11_{22}=5_{12}$&\begin{tikzpicture}[main/.style={circle,fill=black,inner sep=0pt,minimum size=1.5mm},baseline=0,scale=0.6]
			\node[main] at (0,0) (A0) {};
			\node[main] at (1,0) (A1) {};
			\node[main] at (2,0) (A2) {};
			\node[main] at (3,0.5) (A3) {};
			\node[main] at (3,-0.5) (B3) {};
			\node[main] at (4,0.5) (A4) {};
			\node[main] at (4,-0.5) (B4) {};
			\node[main] at (5,0) (A5) {};
			\draw (A0) -- (A1) --  (A2)-- (A3) -- (A4) -- (A5) -- (B4) -- (B3) -- (A2);
	\end{tikzpicture} &$9_{12}\cdot 4_{23}=3_{13}$&\begin{tikzpicture}[main/.style={circle,fill=black,inner sep=0pt,minimum size=1.5mm},baseline=0,scale=0.6]
			\node[main] at (1,-0.5) (A1) {};
			\node[main] at (2,0) (A2) {};
			\node[main] at (3,0) (A3) {};
			\node[main] at (1,0.5) (B2) {};
			\node[main] at (4,0) (A4) {};
			\node[main] at (5,0) (A5) {};
			\node[main] at (6,-0.5) (B6) {};
			\node[main] at (6,0.5) (A6) {};
			\node[main] at (6,0) (C) {};
			\draw (A5) -- (C);
			\draw (B2) -- (A2) -- (A3) -- (A4) -- (A5) -- (B6);
			\draw (A1) --  (A2)-- (A3) -- (A4) -- (A5) -- (A6);
			\end{tikzpicture}\\
$9_{12}\cdot 12_{22}=8_{12}$&\begin{tikzpicture}[main/.style={circle,fill=black,inner sep=0pt,minimum size=1.5mm},baseline=0,scale=0.6]
			\node[main] at (0.5,0) (A0) {};
			\node[main] at (1,0.5) (A1) {};
			\node[main] at (2,0.5) (A2) {};
			\node[main] at (3,0.5) (A3) {};
			\node[main] at (4,0.5) (A4) {};
			\node[main] at (4.5,0) (A5) {};
			\node[main] at (1,-0.5) (A9) {};
			\node[main] at (2,-0.5) (A8) {};
			\node[main] at (3,-0.5) (A7) {};
			\node[main] at (4,-0.5) (A6) {};
			\node[main] at (0,0.5) (B1) {};
			\node[main] at (0,-0.5) (B2) {};
			\draw (B1) -- (A1);
			\draw (B2) -- (A9);
			\draw (A0) -- (A1) --  (A2) -- (A3) -- (A4) -- (A5) -- (A6) -- (A7) -- (A8) -- (A9) -- (A0);
			\end{tikzpicture}&$9_{12}\cdot6_{23}=7_{13}$&\begin{tikzpicture}[main/.style={circle,fill=black,inner sep=0pt,minimum size=1.5mm},baseline=0,scale=0.6]
			\node[main] at (0,0) (A0) {};
			\node[main] at (0.3,0.65) (A1) {};
			\node[main] at (1,1) (A2) {};
			\node[main] at (1.7,0.65) (A4) {};
			\node[main] at (2,0) (A5) {};
			\node[main] at (0.3,-0.65) (A9) {};
			\node[main] at (1,-1) (A7) {};
			\node[main] at (1.7,-0.65) (A6) {};
			\node[main] at (-1,0) (B1) {};
			\node[main] at (-2,0) (B2) {};
			\node[main] at (-3,0) (B3) {};
			\node[main] at (-4,0) (B4) {};
			\draw (B4) -- (B3) -- (B2) -- (B1) -- (A0) -- (A1) --  (A2) -- (A4) -- (A5) -- (A6) -- (A7) -- (A9) -- (A0);
			\end{tikzpicture}\\
$9_{12}\cdot5_{21}=13_{11}$&\begin{tikzpicture}[main/.style={circle,fill=black,inner sep=0pt,minimum size=1.5mm},baseline=0,scale=0.6]
			\node[main] at (0,-0.5) (A0) {};
			\node[main] at (1,-0.5) (A1) {};
			\node[main] at (2,0) (A2) {};
			\node[main] at (3,0) (A3) {};
			\node[main] at (1,0.5) (B2) {};
			\node[main] at (4,-0.5) (B4) {};
			\node[main] at (4,0.5) (A4) {};
			\node[main] at (5,-0.5) (B5) {};
			\node[main] at (5,0.5) (A5) {};
			\draw (B2) -- (A2) -- (A3) -- (B4) -- (B5);
			\draw (A0) -- (A1) --  (A2)-- (A3) -- (A4) -- (A5) ;
			\end{tikzpicture}&$9_{12}\cdot8_{21}=12_{11}$&\begin{tikzpicture}[main/.style={circle,fill=black,inner sep=0pt,minimum size=1.5mm},baseline=0,scale=0.6]
			\node[main] at (0,-0.5) (A0) {};
			\node[main] at (1,-0.5) (A1) {};
			\node[main] at (2,0) (A2) {};
			\node[main] at (3,0) (C) {};
			\node[main] at (4,0) (A3) {};
			\node[main] at (1,0.5) (B2) {};
			\node[main] at (5,0) (B4) {};
			\node[main] at (4,1) (A4) {};
			\node[main] at (6,-0.5) (B5) {};
			\node[main] at (6,0.5) (A5) {};
			\draw (A4) -- (A3);
			\draw (B2) -- (A2) -- (A3) -- (B4) -- (B5);
			\draw (A0) -- (A1) --  (A2)-- (A3) -- (B4) -- (A5) ;
			\end{tikzpicture}\\
$9_{12}\cdot 1_{23}=2_{13}$&\begin{tikzpicture}[main/.style={circle,fill=black,inner sep=0pt,minimum size=1.5mm},baseline=0,scale=0.6]
			\node[main] at (0,0) (A0) {};
			\node[main] at (-1,0.5) (A1) {};
			\node[main] at (-1,-0.5) (B1) {};
			\node[main]  at (-2,0) (A2) {};
			\node[main]  at (1,0) (A3) {};
			\draw (A0) -- (B1) -- (A2) -- (A1) -- (A0) -- (A3);
			\end{tikzpicture}&$9_{12}\cdot2_{21}=10_{11}$&\begin{tikzpicture}[main/.style={circle,fill=black,inner sep=0pt,minimum size=1.5mm},baseline=0,scale=0.6]
			\node[main] at (0,0) (A0) {};
			\node[main] at (-1,0.5) (A1) {};
			\node[main] at (-1,-0.5) (B1) {};
			\node[main]  at (-2,0.5) (A2) {};
			\node[main]  at (-2,-0.5) (B2) {};
			\node[main]  at (1,0.5) (A3) {};
			\node[main]  at (1,-0.5) (B3) {};
			\draw (A2) -- (A1) -- (A0) -- (A3);
			\draw (B2) -- (B1) -- (A0) -- (B3);
			\end{tikzpicture}\\
$9_{12}\cdot 13_{22}=9_{12}$&\begin{tikzpicture}[main/.style={circle,fill=black,inner sep=0pt,minimum size=1.5mm},baseline=0.9,scale=0.6]
\node[main] at (0,0) (A0) {};
\node[main] at (1,0) (A1) {};
\node[main] at (2,0) (A2) {};
\node[main] at (3,0) (A3) {};
\node[main] at (4,0) (A4) {};
\node[main] at (5,0) (A5) {};
\node[main] at (6,0) (A6) {};
\node[main] at (7,0) (A7) {};
\node[main] at (8,0) (A8) {};
\node[main] at (9,0) (A9) {};
\node[main] at (10,0) (A10) {};
\node[main] at (5,0.6) (B1) {};
\node[main] at (5,1.2) (B2) {};
\node[main] at (4.4,1.8) (B3) {};
\node[main] at (5.6,1.8) (B4) {};

\draw (A5) -- (B1) -- (B2) -- (B3);
\draw (B2) -- (B4);
\draw (A0) -- (A1) --  (A2)-- (A3) -- (A4) -- (A5) -- (A6) -- (A7) -- (A8) -- (A9) -- (A10);
\end{tikzpicture} & &Principal graph\\
$9_{12}\cdot9_{21}=14_{11}$&\begin{tikzpicture}[main/.style={circle,fill=black,inner sep=0pt,minimum size=1.5mm},baseline=1.5,scale=0.6]
\node[main] at (0,0) (A0) {};
\node[main] at (1,0) (A1) {};
\node[main] at (2,0) (A2) {};
\node[main] at (3,0) (A3) {};
\node[main] at (4,0) (A4) {};
\node[main] at (5,0) (A5) {};
\node[main] at (6,0) (A6) {};
\node[main] at (7,0) (A7) {};
\node[main] at (8,0) (A8) {};
\node[main] at (9,0) (A9) {};
\node[main] at (5,1) (B1) {};
\node[main] at (6,1) (B2) {};

\draw (A5) -- (B1) ;
\draw (A6) -- (B2);
\draw (A0) -- (A1) --  (A2)-- (A3) -- (A4) -- (A5) -- (A6) -- (A7) -- (A8) -- (A9) ;
\end{tikzpicture} & & Dual Principal graph
\end{tblr}
 \caption{Asaeda-Haagerup module graphs}
    \label{tblr:AHmg}
\end{table}
Unlike the Extended-Haagerup case, there might be other fusion categories $\CC$ in the Morita equivalence class of $\mathcal{AH}_i$. Consequently, we cannot give a full description of the Brauer-Picard groupoid for the Asaeda-Haagerup fusion category. However, in \cite{grossman2016brauer} the authors prove: 

\begin{theo}
    Let $\CC$ be a fusion category which is Morita equivalent to the $\mathcal{AH}_i$, $i=1,2,3$, but not isomorphic as fusion categories to any of them. Then exactly one of the following four cases holds:
    \begin{enumerate}[(a)]
        \item Every $\mathcal{AH}_1$-$\CC$ Morita equivalence realizes $9_1$, every $\mathcal{AH}_2$-$\CC$ Morita equivalence realizes $19_2$, and every $\mathcal{AH}_3$-$\CC$ Morita equivalence realizes $16_3$.
        \item Every $\mathcal{AH}_1$-$\CC$ Morita equivalence realizes $16_1$, every $\mathcal{AH}_2$-$\CC$ Morita equivalence realizes $4_2$, and every $\mathcal{AH}_3$-$\CC$ Morita equivalence realizes $18_3$.
        \item Every $\mathcal{AH}_1$-$\CC$ Morita equivalence realizes $21_1$, every $\mathcal{AH}_2$-$\CC$ Morita equivalence realizes $17_2$, and every $\mathcal{AH}_3$-$\CC$ Morita equivalence realizes $2_3$.
        \item Every $\mathcal{AH}_1$-$\CC$ Morita equivalence realizes $1_1$, every $\mathcal{AH}_2$-$\CC$ Morita equivalence realizes $3_2$, and every $\mathcal{AH}_3$-$\CC$ Morita equivalence realizes $2_3$.
    \end{enumerate}
\end{theo}
Here $a_i$ denotes the $a^{th}$ fusion bimodule on a list of potential $AH_i$ fusion bimodules. Verifying the compatibility of these modules with the bimodule that realizes $9_{12}$ we get
\begin{align*}
    9_{12}\cdot 19_2&=9_1 & 9_{12}\cdot 4_2&=16_1 & 9_{12}\cdot 17_2&=21_1 & 9_{12}\cdot 3_2&=1_1
\end{align*}
Using Algorithm \ref{algorithm} without checking for condition (1c') in step 1b) we recover the multiplication maps for the above modules. We then obtain the following fusion graphs with respect to $X$:
\begin{table}[H]
\centering
\begin{tblr}{
colspec={rlrl},
hline{2} = {1-2}{leftpos = -1, rightpos = -1, endpos},
hline{2} = {3-4}{leftpos = -1, rightpos = -1, endpos},
vline{3},
row{1-2}={abovesep=5pt,belowsep=5pt}
}
$9_{12}\cdot 19_2=9_1$&\begin{tikzpicture}[main/.style={circle,fill=black,inner sep=0pt,minimum size=1.5mm},baseline=0,scale=0.6]
			\node[main] at (0,0) (A0) {};
			\node[main] at (0.3,0.65) (A1) {};
			\node[main] at (1,1) (A2) {};
			\node[main] at (1.7,0.65) (A4) {};
			\node[main] at (2,0) (A5) {};
			\node[main] at (0.3,-0.65) (A9) {};
			\node[main] at (1,-1) (A7) {};
			\node[main] at (1.7,-0.65) (A6) {};
			\node[main] at (-1,0) (B1) {};
			\node[main] at (3,0) (B2) {};
			\draw (B2) -- (A5);
			\draw  (B1) -- (A0) -- (A1) --  (A2) -- (A4) -- (A5) -- (A6) -- (A7) -- (A9) -- (A0);
			\end{tikzpicture}&$9_{12}\cdot 4_2=16_1$&\begin{tikzpicture}[main/.style={circle,fill=black,inner sep=0pt,minimum size=1.5mm},baseline=0,scale=0.6]
			\node[main] at (1,-0.5) (A1) {};
			\node[main] at (2,0) (A2) {};
			\node[main] at (4,0) (A3) {};
			\node[main] at (3,0) (C) {};
			\node[main] at (1,0.5) (B2) {};
			\node[main] at (5,-0.5) (B4) {};
			\node[main] at (5,0.5) (A4) {};
			\node[main] at (3,1) (C2) {};
			\draw (C2) -- (C);
			\draw  (B2) -- (A2) -- (A3) -- (B4);
			\draw  (A1) --  (A2)-- (A3) -- (A4) ;
			\end{tikzpicture}\\
$9_12\cdot 17_2=21_1$&\begin{tikzpicture}[main/.style={circle,fill=black,inner sep=0pt,minimum size=1.5mm},baseline=0,scale=0.6]
			\node[main] at (0,-0.5) (A0) {};
			\node[main] at (1,-0.5) (A1) {};
			\node[main] at (2,0) (A2) {};
			\node[main] at (4,0) (A3) {};
			\node[main] at (3,0) (C) {};
			\node[main] at (0,0.5) (B1) {};
			\node[main] at (1,0.5) (B2) {};
			\node[main] at (5,-0.5) (B4) {};
			\node[main] at (5,0.5) (A4) {};
			\node[main] at (6,-0.5) (B5) {};
			\node[main] at (6,0.5) (A5) {};
			\draw  (B1) -- (B2) -- (A2) -- (A3) -- (B4) -- (B5);
			\draw (A0) -- (A1) --  (A2)-- (A3) -- (A4) -- (A5) ;
			\end{tikzpicture}&$9_{12}\cdot 3_2=1_1$&\begin{tikzpicture}[main/.style={circle,fill=black,inner sep=0pt,minimum size=1.5mm},baseline=0,scale=0.6]
			\node[main] at (0,0) (A0) {};
			\node[main] at (-1,0.5) (A1) {};
			\node[main] at (-1,-0.5) (B1) {};
			\node[main]  at (-2,0) (A2) {};
			\node[main]  at (1,0) (A3) {};
			\draw (A0) -- (B1) -- (A2) -- (A1) -- (A0) -- (A3);
			\end{tikzpicture}
\end{tblr}
\caption{Asaeda-Haagerup potential module graphs}
\label{tblr:AHpmg}
\end{table}
Thus, there are a total of $14$ \emph{different} potential module graphs in Table \ref{tblr:AHmg} and Table \ref{tblr:AHpmg} that are associated to the Asaeda-Haagerup subfactor.

    \section{New commuting squares}\label{newcs}

We construct nondegenerate commuting squares of the form
\begin{equation}\label{consquare}
\begin{array}{ccc}
A_{1,0}&\stackrel{\text{K}}{\sst} &A_{1,1}\\
\rotatebox[origin=c]{90}{$\sst$}_G& &\rotatebox[origin=c]{90}{$\sst$}_G\\
A_{0,0}&\stackrel{\text{H}}{\sst} &A_{0,1}
\end{array}
\end{equation}
where $G$ is one of the following four graphs:
\begin{figure}[h]
	\centering
	\begin{tikzpicture}[main/.style={circle,fill=black,inner sep=0pt,minimum size=1.5mm},
        sub/.style={circle,fill=white,draw=black,inner sep=0pt,minimum size=1.5mm}]
	\node[sub] at (0,0.5) (A0) {};
	\node[sub] at (0,0) (B0) {};
	\node[sub] at (0,-0.5) (C0) {};
	\node[main] at (1,0) (A1) {};
	\node[sub] at (2,0) (A2) {};
	\node[main] at (3,0) (A3) {};
	\node[sub] at (4,0.5) (A10) {};
	\node[sub] at (4,0) (B10) {};
	\node[sub] at (4,-0.5) (C10) {};
	
	\draw (B0) -- (A1) -- (C0);
	\draw (B10) -- (A3) -- (C10);
	\draw (A0) -- (A1) -- (A2) -- (A3) -- (A10);
	\end{tikzpicture}
	\caption{Small double broom}
    \label{smallbroom}
\end{figure}
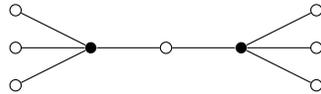

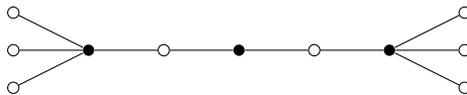
\begin{figure}[h]
	\centering
	\begin{tikzpicture}[main/.style={circle,fill=black,inner sep=0pt,minimum size=1.5mm},
        sub/.style={circle,fill=white,draw=black,inner sep=0pt,minimum size=1.5mm}]
	\node[sub] at (0,0.5) (A0) {};
	\node[sub] at (0,0) (B0) {};
	\node[sub] at (0,-0.5) (C0) {};
	\node[main] at (1,0) (A1) {};
	\node[sub] at (2,0) (A2) {};
	\node[main] at (3,0) (A3) {};
	\node[sub] at (4,0) (A4) {};
	\node[main] at (5,0) (A5) {};
	\node[sub] at (6,0.5) (A10) {};
	\node[sub] at (6,0) (B10) {};
	\node[sub] at (6,-0.5) (C10) {};
	
	\draw (B0) -- (A1) -- (C0);
	\draw (B10) -- (A5) -- (C10);
	\draw (A0) -- (A1) -- (A2) -- (A3) -- (A4) -- (A5) -- (A10);
	\end{tikzpicture}
	\caption{Medium double broom}
    \label{mediumbroom}
\end{figure}

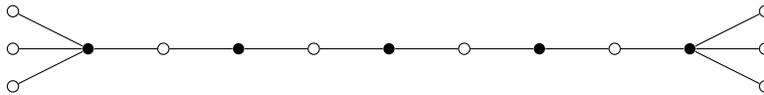
\begin{figure}[h]
	\centering
	\begin{tikzpicture}[main/.style={circle,fill=black,inner sep=0pt,minimum size=1.5mm},
        sub/.style={circle,fill=white,draw=black,inner sep=0pt,minimum size=1.5mm}]
	\node[sub] at (0,0.5) (A0) {};
	\node[sub] at (0,0) (B0) {};
	\node[sub] at (0,-0.5) (C0) {};
	\node[main] at (1,0) (A1) {};
	\node[sub] at (2,0) (A2) {};
	\node[main] at (3,0) (A3) {};
	\node[sub] at (4,0) (A4) {};
	\node[main] at (5,0) (A5) {};
	\node[sub] at (6,0) (A6) {};
	\node[main] at (7,0) (A7) {};
	\node[sub] at (8,0) (A8) {};
	\node[main] at (9,0) (A9) {};
	\node[sub] at (10,0.5) (A10) {};
	\node[sub] at (10,0) (B10) {};
	\node[sub] at (10,-0.5) (C10) {};
	
	\draw (B0) -- (A1) -- (C0);
	\draw (B10) -- (A9) -- (C10);
	\draw (A0) -- (A1) -- (A2) -- (A3) -- (A4) -- (A5) -- (A6) -- (A7) -- (A8) -- (A9) -- (A10);
	\end{tikzpicture}
	\caption{Large double broom}
    \label{largebroom}
\end{figure}

\begin{figure}[H]
	\centering
	\begin{tikzpicture}[main/.style={circle,fill=black,inner sep=0pt,minimum size=1.5mm},
        sub/.style={circle,fill=white,draw=black,inner sep=0pt,minimum size=1.5mm}]
	\node[main] at (-6,0) (A3) {};
    \node[sub] at (-5,0) (a4) {};
    \node[main] at (-4,0) (A2) {};
    \node[sub] at (-3,0) (a2) {};
    \node[main] at (-2,0) (A1) {};
    \node[sub] at (-1,0) (a1) {};
    \node[main] at (0,0) (A0) {};
    \node[sub] at (1,0) (b1) {};
    \node[main] at (2,0) (B1) {};
    \node[sub] at (3,0) (b2) {};
    \node[main] at (4,0) (B2) {};
    \node[sub] at (5,0) (b4) {};
    \node[main] at (6,0) (B3) {};
    \node[sub] at (-4,-0.7) (a3) {};
    \node[sub] at (0,-0.7) (a0) {};
    \node[sub] at (4,-0.7) (b3) {};
	
    \draw (A3) -- (a4) -- (A2) -- (a2) -- (A1) -- (a1) -- (A0) -- (b1) -- (B1) -- (b2) -- (B2) -- (b4) -- (B3);
    \draw (A2) -- (a3);
    \draw (A0) -- (a0);
    \draw (B2) -- (b3);
	\end{tikzpicture}
	\caption{Quipu}
    \label{quipu}
\end{figure}

The graphs have norms $\sqrt{5}$, $\sqrt{3+\sqrt{3}}$, $\sqrt{\frac{5+\sqrt{17}}{2}}$ and $\sqrt{4.37720\dots}$ respectively. Moreover, all of them satisfy Wenzl's criterion of irreducibility (Theorem \ref{wenzcrit}). We will follow the set up in Chapter 8 of \cite{schou2013commuting} to present the bi-unitary connections that yield the commuting squares.

Let $\ZZ_{0,0}$ (respective $\ZZ_{0,1}$, $\ZZ_{1,0}$ and $\ZZ_{1,1}$) denote the set of minimal central projections in  $A_{0,0}$ (respectively $A_{0,1}$, $A_{1,0}$ and $A_{1,1}$). The elements of $\ZZ_{0,0}$ and $\ZZ_{0,1}$ are labelled by $\V_+$ (denoted by black vertices) and the elements of $\ZZ_{1,0}$ and $\ZZ_{1,1}$ are labelled by $\V_-$ (denoted by white vertices), where $\V=\V_+\cup\V_-$ are the vertices of the bipartite graph $G$.

It is shown in \cite[Chapter 1]{schou2013commuting} that the existence of a commuting square of finite dimensional C*-algebras of the form \ref{consquare} is equivalent to constructing a unitary matrix $u$, satisfying the so-called bi-unitary condition. $u$ is a block-diagonal matrix of the form
\[ u=\bigoplus_{(p,s)} u^{(p,s)} \]
where the labels of the blocks $(p,s)$ runs over all $p\in \ZZ_{0,0}$ and all $s\in \ZZ_{1,1}$ that are connected by a path on $H$ and $G$ through $\ZZ_{0,1}$ (respectively, a path on $G$ and $K$ through $\ZZ_{1,0}$). Each direct summand $u^{(p,s)}$ is a $n(p,s)\times n(p,s)$-matrix, where $n(p,s)$ is the number of paths on $H$ and $G$ from $p$ to $s$ through $\ZZ_{0,1}$ (respectively, the number of paths on $G$ and $K$ through $\ZZ_{1,0}$), so each block $u^{(p,s)}$ is a matrix indexed as follows
\[ u^{(p,s)} = \left(u^{(p,s)}_{q,r}\right)_{q,r}, \]
where $(q,r)$ runs over all possible $q\in\ZZ_{1,0}$ and $r\in \ZZ_{0,1}$, that a path from $p$ to $s$ can go through. Since the vertical inclusions, $A_{0,0}\sst A_{1,0}$ and $A_{0,1}\sst A_{1,1}$, do not have multiple edges, each $u^{(p,s)}_{q,r}$ is a $m\times n$-matrix, where $m$ is the multiplicity of the edge $qs$ in the inclusion $A_{1,0}\sst A_{1,1}$ and $n$ is the multiplicity of the edge $pr$ in the inclusion $A_{0,0}\sst A_{0,1}$. We may assume $A_{0,0}$ is Abelian, i.e. of the form $A_{0,0}=\C^{|\ZZ_{0,0}|}$.

Let $\lambda(\cdot)$ (respectively $\eta(\cdot)$) denote a fixed Perron-Frobenius vector for the inclusion graph of $A_{0,0}\sst A_{1,0}$ (respectively $A_{0,1}\sst A_{1,1}$). Set
\begin{equation}\label{biunitfactor}
w(p,q,r,s)=\sqrt{\frac{\lambda(p)\eta(s)}{\lambda(q)\eta(r)}}.    
\end{equation}
Define a matrix $v$ by
\[ v=\bigoplus_{(q,r)}v^{(q,r)}. \]
$v^{(q,r)}$ is a square matrix, which can be written as a block matrix $v^{(q,r)}=\left(v^{(q,r)}_{p,s}\right)_{p,s}$ with each block given by
\begin{equation}\label{biunit}
v^{(q,r)}_{p,s}=w(p,q,r,s)(u^{(p,s)}_{q,r})^*,
\end{equation}
where $(p,q,r,s)$ runs through all quadruples in $\ZZ_{0,0}\times \ZZ_{1,0}\times \ZZ_{0,1}\times \ZZ_{1,1}$, which can be completed to a cycle $p-r-s-q-p$. Here $X^*$ denotes the conjugate transpose of the matrix $X$.

The existence of a commuting square as in \ref{consquare} is then equivalent to the existence of $u$ and $v$ as above such that both are unitaries (this is referred to as the \emph{bi-unitary condition}). In this case, we say $u$ and $v$ are a \emph{bi-unitary connection} for the inclusions \ref{consquare}.

\begin{rem}
    In our case, both vertical inclusions are described by the same graph $G$ hence we can use $\lambda(\cdot)$ to denote a Perron-Frobenius vector of the inclusion graph for $A_{0,0}\sst A_{1,0}$ (or $A_{0,1}\sst A_{1,1}$).
\end{rem}

Surprisingly, we only needed to consider real unitaries $u$ and $v$ to determine the bi-unitary connections for the desired commuting squares.

\subsection{Horizontal inclusions}

In \cite{schou2013commuting} a large family of commuting squares is constructed where $H=GG^t-I$ and $K=G^tG-I$. Ocneanu constructed a commuting square of the form (\ref{consquare}) where $G=E_{10}$ and the horizontal inclusions are of the form $H=p(GG^t)$ and $K=p(G^tG)$ where $p$ is a polynomial with integer coefficients. These commuting squares are particularly nice as they are nondegenerate, hence, they can be used to construct a subfactor with index $\|G\|^2$. Ocneanu's example produces an irreducible hyperfinite subfactor with the smallest known index greater than $4$. It has index $\|E_{10}\|^2\approx 4.026418$, see \cite[Chapter 8]{schou2013commuting}.

In the case where $G$ is one of the double brooms above, we could show that there is no nondegenerate commuting square of the form (\ref{consquare}) if $K$ is a polynomial of $G^t G$. In this case, the eigenspaces of $GG^t$ are all 1-dimensional. The nondegeneracy of the commuting square implies that $HGG^t=GG^tH$. By simultaneously diagonalizing both $H$ and $GG^t$ one can show that $H=p(GG^t)$ where $p$ is a polynomial with potentially non-integer coefficients. 

To obtain a nondegenerate commuting square based on one of these double brooms we only need to consider $H=p(GG^t)$ and $K\neq p(G^tG)$ such that both $H$ and $K$ are incidence matrices for connected bipartite graphs. To simplify this search (and future computations) we assume $K$ to be symmetric about both its diagonals. Under these assumptions, given a fixed polynomial $p$ one obtains relations for the entries of $K$. To further simplify the computations we choose $p$ to be monic as this guarantees that $HG$ will have small entries. We then proceed to perturb $p(G^tG)$ according to these relations and use the bi-unitary condition to test if our choice is compatible with the blocks of $u$ and $v$ being unitary. Most choices of $p$ fail this test very quickly as they generally imply that one row or column of some block has norm greater than $1$.

For example, let $\Gamma$ be the small broom in Figure \ref{smallbroom} with the following labelling of the vertices:
\begin{figure}[H]
	\centering
	$\Gamma=$\begin{tikzpicture}[main/.style={circle,fill=black,inner sep=0pt,minimum size=1.5mm},
        sub/.style={circle,fill=white,draw=black,inner sep=0pt,minimum size=1.5mm},baseline=-2]
	\node[sub,label=left:{$a_1$}] at (0,0.5) (A0) {};
	\node[sub,label=left:{$a_2$}] at (0,0) (B0) {};
	\node[sub,label=left:{$a_3$}] at (0,-0.5) (C0) {};
	\node[main,label=above:{$A$}] at (1,0) (A1) {};
	\node[sub,label=above:{$a_0$}] at (2,0) (A2) {};
	\node[main,label=above:{$B$}] at (3,0) (A3) {};
	\node[sub,label=right:{$b_1$}] at (4,0.5) (A10) {};
	\node[sub,label=right:{$b_2$}] at (4,0) (B10) {};
	\node[sub,label=right:{$b_3$}] at (4,-0.5) (C10) {};
	
	\draw (B0) -- (A1) -- (C0);
	\draw (B10) -- (A3) -- (C10);
	\draw (A0) -- (A1) -- (A2) -- (A3) -- (A10);
	\end{tikzpicture}
\end{figure}
where $\V_+=(A,B)$ and $\V_-=(a_3,a_2,a_1,a_0,b_1,b_2,b_3)$. Then $G$ is the following $\V_+\times \V_-$ matrix
\[ G=\left(\begin{array}{ccccccc}
1 & 1 & 1 & 1 & 0 & 0 & 0 \\
0 & 0 & 0 & 1 & 1 & 1 & 1
\end{array}
\right) \]
We computed and fixed a Perron-Frobenius eigenvector for the graph $\Gamma$ as follows:
\begin{figure}[H]
	\centering
	$\lambda(\cdot)=$\begin{tikzpicture}[main/.style={circle,fill=black,inner sep=0pt,minimum size=1.5mm},
        sub/.style={circle,fill=white,draw=black,inner sep=0pt,minimum size=1.5mm},baseline=-2]
	\node[sub,label=left:{$1$}] at (0,0.5) (A0) {};
	\node[sub,label=left:{$1$}] at (0,0) (B0) {};
	\node[sub,label=left:{$1$}] at (0,-0.5) (C0) {};
	\node[main,label=above:{$5$}] at (1,0) (A1) {};
	\node[sub,label=above:{$2$}] at (2,0) (A2) {};
	\node[main,label=above:{$5$}] at (3,0) (A3) {};
	\node[sub,label=right:{$1$}] at (4,0.5) (A10) {};
	\node[sub,label=right:{$1$}] at (4,0) (B10) {};
	\node[sub,label=right:{$1$}] at (4,-0.5) (C10) {};
	
	\draw (B0) -- (A1) -- (C0);
	\draw (B10) -- (A3) -- (C10);
	\draw (A0) -- (A1) -- (A2) -- (A3) -- (A10);
	\end{tikzpicture}
\end{figure}

Since we are assuming the horizontal inclusions are described by matrices that are symmetric about both their diagonals, they are of the form
\[ H=\left(\begin{array}{cc}
h_{11} & h_{12} \\
h_{12} & h_{11}
\end{array}\right), \quad K=\left(\begin{array}{ccccccc}
k_{11} & k_{12} & k_{13} & k_{14} & k_{15} & k_{16} & k_{17} \\
k_{12} & k_{22} & k_{23} & k_{24} & k_{25} & k_{26} & k_{16} \\
k_{13} & k_{23} & k_{33} & k_{34} & k_{35} & k_{25} & k_{15} \\
k_{14} & k_{24} & k_{34} & k_{44} & k_{34} & k_{24} & k_{14} \\
k_{15} & k_{25} & k_{35} & k_{34} & k_{33} & k_{23} & k_{13} \\
k_{16} & k_{26} & k_{25} & k_{24} & k_{23} & k_{22} & k_{12} \\
k_{17} & k_{16} & k_{15} & k_{14} & k_{13} & k_{12} & k_{11} 
\end{array}
\right)\]
where the vertices are in the order $(A,B)$ and $(a_3,a_2,a_1,a_0,b_1,b_2,b_3)$ respectively, i.e. $H$ and $K$ are $\V_+\times \V_+$ and $\V_-\times \V_-$ matrices respectively.

To understand the block structure of the summand $u^{(A,a_0)}$ we look at the possible $A-r-a_0-q-A$ cycles, where $r\in \V_+$ and $q\in \V_-$:
\begin{center}
\begin{tikzpicture}
	\node at (0,0) (L) {$A$};
	\node at (0.66,3.33) (a2) {$a_2$};
    \node at (1.33,2.33) (a1) {$a_1$};
    \node at (-0.33,4) (a3) {$a_3$};
    \node at (2,1.5) (a0) {$a_0$};
	\node at (4,4) (11) {$a_0$};
	\node at (3.66,1) (A) {$A$};
    \node at (4.33,0) (B) {$B$};
	\draw (L) -- (a3);
    \draw (L) -- (a2);
    \draw (L) -- (a1);
    \draw (L) -- (a0);
	\draw (a3) to node[midway,circle,inner sep=0pt,fill=white] {$k_{14}$} (11);
    \draw (a2) to node[midway,circle,inner sep=0pt,fill=white] {$k_{24}$} (11);
    \draw (a1) to node[midway,circle,inner sep=0pt,fill=white] {$k_{34}$} (11);
    \draw (a0) to node[midway,circle,inner sep=0pt,fill=white] {$k_{44}$} (11);
	\draw (L) to node[midway,circle,inner sep=0pt,fill=white] {$h_{11}$} (A);
    \draw (L) to node[midway,circle,inner sep=0pt,fill=white] {$h_{12}$} (B);
	\draw (11) -- (A);
    \draw (11) -- (B);
\end{tikzpicture}
\end{center}
This implies that $u^{(A,a_0)}$ has the following block structure
$$
u^{(A,a_0)}=\begin{pNiceArray}{W{c}{1cm}|W{c}{1cm}}[first-row,first-col,last-row,last-col]
& A & B\\
a_3& \ast & \ast & k_{14}\\ \hline
a_2& \ast & \ast & k_{24}\\ \hline
a_1& \ast & \ast & k_{34}\\ \hline
a_0& \ast & \ast & k_{44}\\
 & h_{11} & h_{12} &
\end{pNiceArray}
$$

where the labels on the bottom and the right denote the number of columns and rows for the blocks respectively. A similar analysis implies the following block structure for the following summands of $v$,
\begin{align*}
v^{(a_3,B)}&=\begin{pNiceArray}{c|c|c|c}[first-row,first-col,last-row,last-col]
& a_0 & b_1 & b_2 & b_3 &\\
A& \ast & \ast & \ast & \ast & h_{12}\\
 & k_{14} & k_{15} & k_{16} & k_{17} &
\end{pNiceArray}\;,\\
v^{(a_2,B)}&=\begin{pNiceArray}{c|c|c|c}[first-row,first-col,last-row,last-col]
& a_0 & b_1 & b_2 & b_3 &\\
A& \ast & \ast & \ast & \ast & h_{12}\\
 & k_{24} & k_{25} & k_{26} & k_{27} &
\end{pNiceArray}\;,\\
v^{(a_1,B)}&=\begin{pNiceArray}{c|c|c|c}[first-row,first-col,last-row,last-col]
& a_0 & b_1 & b_2 & b_3 &\\
A& \ast & \ast & \ast & \ast & h_{12}\\
 & k_{34} & k_{35} & k_{36} & k_{37} &
\end{pNiceArray}\;.
\end{align*}
In particular, since each block has to be unitary, the sum of the squares of the entries of each column must add up to 1. Therefore $\|v^{(a_3,B)}_{A,a_0}\|_2^2=k_{14}$, $\|v^{(a_2,B)}_{A,a_0}\|_2^2=k_{24}$ and $\|v^{(a_1,B)}_{A,a_0}\|_2^2=k_{34}$, where $\|X\|_2^2$ denotes the sum of the squares of the absolute values of the entries of the block $X$. The bi-unitary condition (\ref{biunit}) then implies that 
$$\|u^{(A,a_0)}_{a_3,B}\|_2^2=\frac{1}{2}k_{14},\;\|u^{(A,a_0)}_{a_2,B}\|_2^2=\frac{1}{2}k_{24},\; \|u^{(A,a_0)}_{a_1,B}\|_2^2=\frac{1}{2}k_{34}.$$
Since $h_{12}=\sum_{i=0}^3\|u^{(A,a_0)}_{a_i,B}\|_2^2$, we have $2h_{12}\geq k_{14}+k_{24}+k_{34}$. Similarly, we obtain other inequalities that relate the entries of $K$ with the entries of $H$. We then use these inequalities to perturb $p(G^tG)$ to obtain a suitable candidate for $K$.

A similar analysis is used to obtain the horizontal inclusions for commuting squares based on the ``Quipu''.

\subsection{Small double broom}

For the ``small double broom'' we will use the following labeling of the vertices:
\begin{figure}[H]
	\centering
	$\Gamma=$\begin{tikzpicture}[main/.style={circle,fill=black,inner sep=0pt,minimum size=1.5mm},
        sub/.style={circle,fill=white,draw=black,inner sep=0pt,minimum size=1.5mm},baseline=-2]
	\node[sub,label=left:{$a_1$}] at (0,0.5) (A0) {};
	\node[sub,label=left:{$a_2$}] at (0,0) (B0) {};
	\node[sub,label=left:{$a_3$}] at (0,-0.5) (C0) {};
	\node[main,label=above:{$A$}] at (1,0) (A1) {};
	\node[sub,label=above:{$a_0$}] at (2,0) (A2) {};
	\node[main,label=above:{$B$}] at (3,0) (A3) {};
	\node[sub,label=right:{$b_1$}] at (4,0.5) (A10) {};
	\node[sub,label=right:{$b_2$}] at (4,0) (B10) {};
	\node[sub,label=right:{$b_3$}] at (4,-0.5) (C10) {};
	
	\draw (B0) -- (A1) -- (C0);
	\draw (B10) -- (A3) -- (C10);
	\draw (A0) -- (A1) -- (A2) -- (A3) -- (A10);
	\end{tikzpicture}
\end{figure}
where $\V_+=(A,B)$ and $\V_-=(a_3,a_2,a_1,a_0,b_1,b_2,b_3)$.

Writing the vertices in the order $(A,B,a_3,a_2,a_1,a_0,b_1,b_2,b_3)$ the adjacency matrix of $\Gamma$ is 
\[ \Delta_\Gamma=\left(\begin{array}{cc}
0 & G \\
G^t & 0
\end{array}
\right) \]
where 
\[ G=\left(\begin{array}{ccccccc}
1 & 1 & 1 & 1 & 0 & 0 & 0 \\
0 & 0 & 0 & 1 & 1 & 1 & 1
\end{array}
\right) \]
describes the vertical inclusions of (\ref{consquare}). In this case, we will try to build a commuting square with horizontal inclusion matrices given by
\[ H=\left(\begin{array}{cc}
3 & 1 \\
1 & 3
\end{array}\right), \quad K=\left(\begin{array}{ccccccc}
0 & 1 & 1 & 1 & 0 & 0 & 0 \\
1 & 0 & 1 & 1 & 0 & 0 & 0 \\
1 & 1 & 1 & 0 & 1 & 0 & 0 \\
1 & 1 & 0 & 2 & 0 & 1 & 1 \\
0 & 0 & 1 & 0 & 1 & 1 & 1 \\
0 & 0 & 0 & 1 & 1 & 0 & 1 \\
0 & 0 & 0 & 1 & 1 & 1 & 0
\end{array}
\right) \]
where the vertices are in the order $(A,B)$ for $H$, respectively $(a_3,a_2,a_1,a_0,b_1,b_2,b_3)$ for $K$. Observe that $H$ and $K$ have some multiple edges.

\begin{figure}[H]
	\centering
	\includegraphics[page=1]{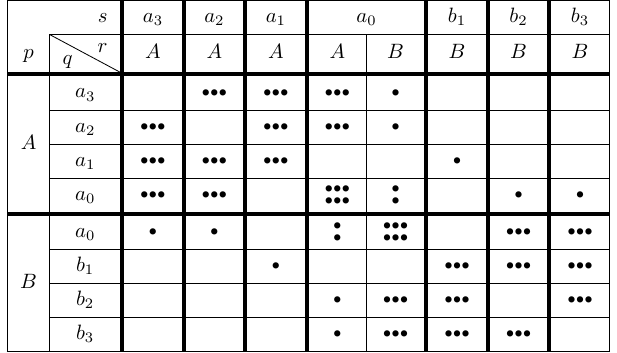}
	\caption{Small double broom - structure of $u$}
    \label{Diagu}
\end{figure}

Note that the $(p,s)$-entry of $GK=HG$ counts the number of paths on $G$ and $K$ going from $p$ to $s$ through $\ZZ_{1,0}$. Since
\[ GK=\left(\begin{array}{ccccccc}
3 & 3 & 3 & 4 & 1 & 1 & 1 \\
1 & 1 & 1 & 4 & 3 & 3 & 3
\end{array}
\right) \]
we deduce that $u$ will have two $4\times 4$ blocks, six $3\times 3$ blocks and six $1\times 1$ blocks. Figure \ref{Diagu} describes the block structure of each direct summand of $u$.

In this table, the thick lines delimit each summand of $u$ and the dots describe the block structure of each summand. For example the summand $u^{A,a_3}$ is a $3\times 3$ matrix made up of three blocks. Each block is a $1\times 3$ arrangement correspoding to each of the following $p-r-s-q-p$ cycles:

\begin{figure}[H]
	\centering
	\begin{tikzpicture}
	\node at (0,0) (L) {$A$};
	\node at (0,1) (10) {$a_2$};
	\node at (1,1) (11) {$a_3$};
	\node at (1,0) (R) {$A$};
	\draw (L) -- (10);
	\draw (R) -- (11);
	\draw (10) -- (11);
	\draw (L) -- (R);
	\path (L.15) edge (R.165);
	\path (L.-15) edge (R.-165);
	\end{tikzpicture}
	\quad
	\begin{tikzpicture}
	\node at (0,0) (L) {$A$};
	\node at (0,1) (10) {$a_1$};
	\node at (1,1) (11) {$a_3$};
	\node at (1,0) (R) {$A$};
	\draw (L) -- (10);
	\draw (R) -- (11);
	\draw (10) -- (11);
	\draw (L) -- (R);
	\path (L.15) edge (R.165);
	\path (L.-15) edge (R.-165);
	\end{tikzpicture}
	\quad
	\begin{tikzpicture}
	\node at (0,0) (L) {$A$};
	\node at (0,1) (10) {$a_0$};
	\node at (1,1) (11) {$a_3$};
	\node at (1,0) (R) {$A$};
	\draw (L) -- (10);
	\draw (R) -- (11);
	\draw (10) -- (11);
	\draw (L) -- (R);
	\path (L.15) edge (R.165);
	\path (L.-15) edge (R.-165);
	\end{tikzpicture}
\end{figure}
Similarly, we can describe the block structure of each direct summand of $v$:

\begin{figure}[H]
	\centering
	\includegraphics[page=2]{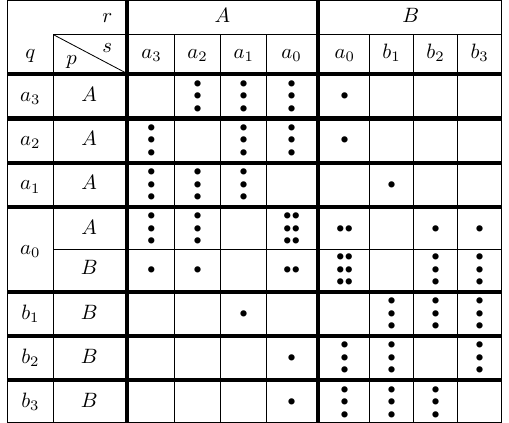}
	\caption{Small double broom - structure of $v$}
\end{figure}

Recall the values of the fixed Perron-Frobenius eigenvector for the graph $\Gamma$ that we will use:
\begin{figure}[H]
	\centering
	$\lambda(\cdot)=$\begin{tikzpicture}[main/.style={circle,fill=black,inner sep=0pt,minimum size=1.5mm},
        sub/.style={circle,fill=white,draw=black,inner sep=0pt,minimum size=1.5mm},baseline=-2]
	\node[sub,label=left:{$1$}] at (0,0.5) (A0) {};
	\node[sub,label=left:{$1$}] at (0,0) (B0) {};
	\node[sub,label=left:{$1$}] at (0,-0.5) (C0) {};
	\node[main,label=above:{$5$}] at (1,0) (A1) {};
	\node[sub,label=above:{$2$}] at (2,0) (A2) {};
	\node[main,label=above:{$5$}] at (3,0) (A3) {};
	\node[sub,label=right:{$1$}] at (4,0.5) (A10) {};
	\node[sub,label=right:{$1$}] at (4,0) (B10) {};
	\node[sub,label=right:{$1$}] at (4,-0.5) (C10) {};
	
	\draw (B0) -- (A1) -- (C0);
	\draw (B10) -- (A3) -- (C10);
	\draw (A0) -- (A1) -- (A2) -- (A3) -- (A10);
	\end{tikzpicture}
\end{figure}

For both $u$ and $v$ we choose the $1\times 1$ blocks to be $1$ and using the bi-unitary condition (\ref{biunit}) we obtain the following values for both $u$ and $v$:
\begin{figure}[H]
	\centering
	\includegraphics[page=3]{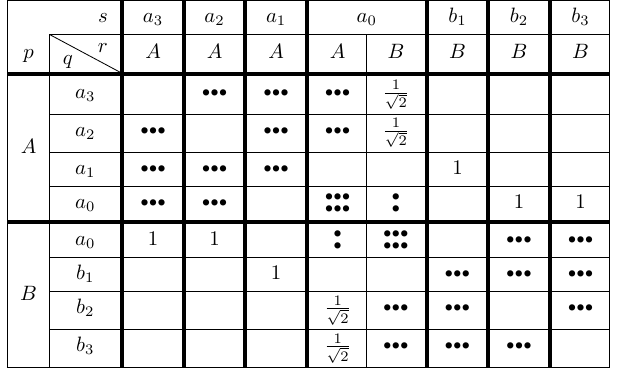}
\end{figure}

\begin{figure}[H]
	\centering
	\includegraphics[page=4]{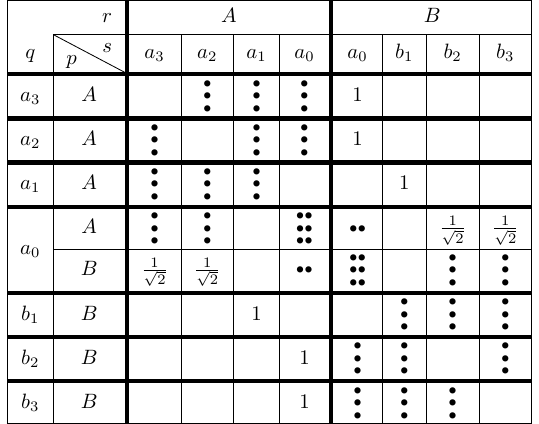}
\end{figure}

Since the columns and rows of a unitary have norm $1$ we get that the blocks $u^{(A,a_0)}_{a_0,B}$, $u^{(B,a_0)}_{a_0,A}$, $v^{(a_0,A)}_{B,a_0}$ and $v^{(a_0,B)}_{A,a_0}$ are all $0$. At this point, we make some choices for $u$ and $v$ and set
\begin{align*}
\begin{pNiceArray}{c}
u^{(A,a_0)}_{a_3,A} \\ \hline
u^{(A,a_0)}_{a_2,A}\\ \hline
u^{(A,a_0)}_{a_0,A}
\end{pNiceArray} & = \begin{pNiceArray}{ccc}[small]
0 & 0 & \sfrac{1}{\sqrt{2}}\\ \hline
0 & 0 & \sfrac{-1}{\sqrt{2}}\\ \hline
1 & 0 & 0 \\
0 & 1 & 0
\end{pNiceArray} & \begin{pNiceArray}{c}
u^{(B,a_0)}_{a_0,B} \\ \hline
u^{(B,a_0)}_{b_2,B}\\ \hline
u^{(B,a_0)}_{b_3,B}
\end{pNiceArray} & = \begin{pNiceArray}{ccc}[small]
0 & 1 & 0\\
0 & 0 & 1\\ \hline
\sfrac{-1}{\sqrt{2}} & 0 & 0 \\ \hline
\sfrac{1}{\sqrt{2}} & 0 & 0
\end{pNiceArray}\\
\begin{pNiceArray}{c|c}
v^{(a_0,A)}_{A,a_3} & v^{(a_0,A)}_{A,a_2}
\end{pNiceArray} & =\begin{pNiceArray}{c|c}[small]
0 & 0 \\ 
0 & 0 \\
\sfrac{1}{\sqrt{2}} & \sfrac{-1}{\sqrt{2}}
\end{pNiceArray} & \begin{pNiceArray}{c|c}
v^{(a_0,B)}_{B,b_2} & v^{(a_0,B)}_{B,b_3}
\end{pNiceArray} & =\begin{pNiceArray}{c|c}[small]
\sfrac{-1}{\sqrt{2}} & \sfrac{1}{\sqrt{2}} \\ 
0 & 0 \\
0 & 0
\end{pNiceArray}
\end{align*}

Using the bi-unitary condition we obtain some of the entries of $u$ and $v$ listed in the next tables:
\begin{figure}[H]
	\centering
	\includegraphics[page=5]{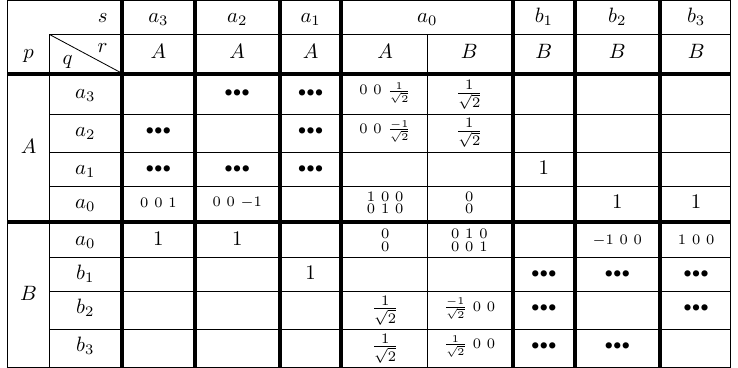}
	\caption{Small double broom - partial entries of $u$}
\end{figure}

\begin{figure}[H]
	\centering
	\includegraphics[page=6]{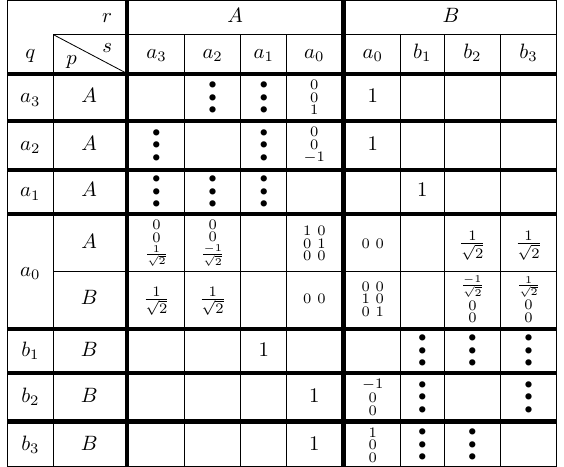}
	\caption{Small double broom - partial entries of $v$}
\end{figure}
We now have multiple ways to select the remaining entries. We present one way to do this:
\begin{figure}[H]
	\centering
	\includegraphics[page=7]{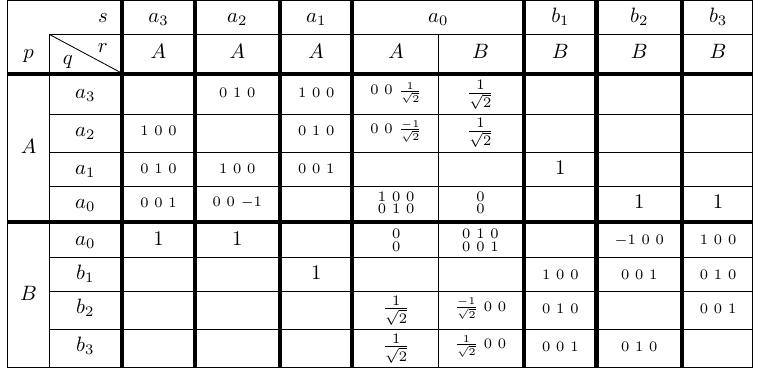}
	\caption{Small double broom - entries of $u$}
\end{figure}

\begin{figure}[H]
	\centering
	\includegraphics[page=8]{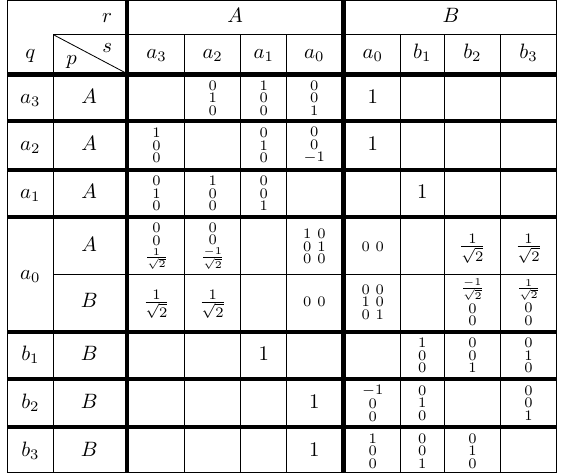}
	\caption{Small double broom - entries of $v$}
\end{figure}

Since both $H$ and $K$ correspond to connected bipartite graphs and since $G$ satisfies Wenzl's criterion for irreducibility (Theorem \ref{wenzcrit}), we have constructed an irreducible subfactor of the hyperfinite II$_1$ factor with index $5$.

\begin{rem}
	We made several choices to determine the entries of $u$ and $v$. One can show that any other choice for $u$ (and consequently for $v$) will be \emph{gauge equivalent} to the one above, i.e. it will be of the form $w_1 u w_2$ where $w_i$ is a unitary matrix (see \cite[Remark 5.3.3]{jones1997introduction}). This means that any subfactor constructed using any of these choices will be isomorphic to the one constructed above.
\end{rem}

\subsection{Medium double broom}

For the ``medium double broom'' we will use the following labeling of the vertices:

\begin{center}
	$\Gamma=$\begin{tikzpicture}[main/.style={circle,fill=black,inner sep=0pt,minimum size=1.5mm},
        sub/.style={circle,fill=white,draw=black,inner sep=0pt,minimum size=1.5mm},baseline=-2]
	\node[sub,label=left:{$a_1$}] at (0,0.5) (A0) {};
	\node[sub,label=left:{$a_2$}] at (0,0) (B0) {};
	\node[sub,label=left:{$a_3$}] at (0,-0.5) (C0) {};
	\node[main,label=above:{$A_1$}] at (1,0) (A1) {};
	\node[sub,label=above:{$a_0$}] at (2,0) (A2) {};
	\node[main,label=above:{$A_0$}] at (3,0) (A3) {};
	\node[sub,label=above:{$b_0$}] at (4,0) (A4) {};
	\node[main,label=above:{$B_1$}] at (5,0) (A5) {};
	\node[sub,label=right:{$b_1$}] at (6,0.5) (A10) {};
	\node[sub,label=right:{$b_2$}] at (6,0) (B10) {};
	\node[sub,label=right:{$b_3$}] at (6,-0.5) (C10) {};
	
	\draw (B0) -- (A1) -- (C0);
	\draw (B10) -- (A5) -- (C10);
	\draw (A0) -- (A1) -- (A2) -- (A3) -- (A4) -- (A5) -- (A10);
	\end{tikzpicture}
\end{center}

and the following Perron-Frobenius eigenvector for the graph $\Gamma$:

\begin{center}
	$\lambda(\cdot)=$\begin{tikzpicture}[main/.style={circle,fill=black,inner sep=0pt,minimum size=1.5mm},
        sub/.style={circle,fill=white,draw=black,inner sep=0pt,minimum size=1.5mm},baseline=-2]
	\node[sub,label=left:{$1$}] at (0,0.5) (A0) {};
	\node[sub,label=left:{$1$}] at (0,0) (B0) {};
	\node[sub,label=left:{$1$}] at (0,-0.5) (C0) {};
	\node[main,label={[label distance=0.1cm]above:$3+\sqrt{3}$}] at (1,0) (A1) {};
	\node[sub,label=below:{$\sqrt{3}$}] at (2,0) (A2) {};
	\node[main,label=above:{$2\sqrt{3}$}] at (3,0) (A3) {};
	\node[sub,label=below:{$\sqrt{3}$}] at (4,0) (A4) {};
	\node[main,label={[label distance=0.1cm]:$3+\sqrt{3}$}] at (5,0) (A5) {};
	\node[sub,label=right:{$1$}] at (6,0.5) (A10) {};
	\node[sub,label=right:{$1$}] at (6,0) (B10) {};
	\node[sub,label=right:{$1$}] at (6,-0.5) (C10) {};
	
	\draw (B0) -- (A1) -- (C0);
	\draw (B10) -- (A5) -- (C10);
	\draw (A0) -- (A1) -- (A2) -- (A3) -- (A4) -- (A5) -- (A10);
	\end{tikzpicture}	
\end{center}

Using the same notation as in the previous section, we thus have
\begin{align*}
G &= \begin{pmatrix}
1 & 1 & 1 & 1 & 0 & 0 & 0 & 0\\
0 & 0 & 0 & 1 & 1 & 0 & 0 & 0\\
0 & 0 & 0 & 0 & 1 & 1 & 1 & 1 
\end{pmatrix}
\end{align*}
for the vertical inclusions. We computed
\begin{align*}
H &= \begin{pmatrix}
2 & 1 & 1 \\
1 & 1 & 1 \\
1 & 1 & 2
\end{pmatrix}, & K &= \begin{pmatrix}
1 & 0 & 0 & 1 & 0 & 0 & 0 & 1 \\
0 & 1 & 0 & 1 & 0 & 0 & 1 & 0 \\
0 & 0 & 2 & 0 & 1 & 0 & 0 & 0 \\
1 & 1 & 0 & 1 & 1 & 1 & 0 & 0 \\
0 & 0 & 1 & 1 & 1 & 0 & 1 & 1 \\
0 & 0 & 0 & 1 & 0 & 2 & 0 & 0 \\
0 & 1 & 0 & 0 & 1 & 0 & 1 & 0 \\
1 & 0 & 0 & 0 & 1 & 0 & 0 & 1
\end{pmatrix}
\end{align*}
for the horizontal ones.

Recall that the $(p,s)$-entry of $GK$ counts the number of paths from $p$ to $s$. Since
\[ GK= \begin{pmatrix}
2 & 2 & 2 & 3 & 2 & 1 & 1 & 1 \\
1 & 1 & 1 & 2 & 2 & 1 & 1 & 1 \\
1 & 1 & 1 & 2 & 3 & 2 & 2 & 2
\end{pmatrix}
 \]
we deduce that $u$ has two $3\times 3$ blocks, ten $2\times 2$ blocks and twelve $1\times 1$ blocks. The following tables describe the block structure of each direct summand of $u$ and $v$:
\begin{figure}[h]
	\centering
	\includegraphics[page=9]{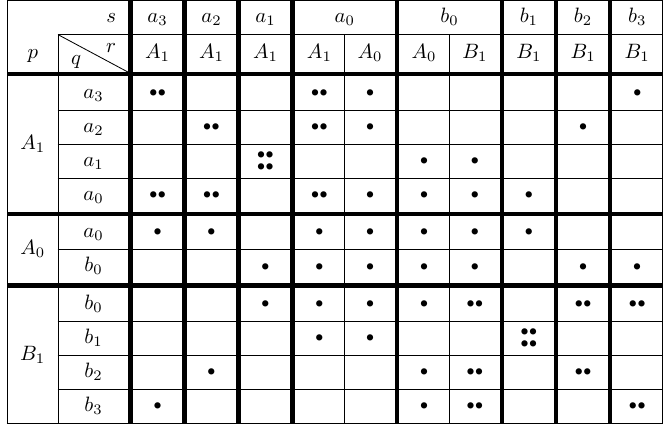}
	\caption{Medium double broom - structure of $u$}
	\label{MedDiagU}
\end{figure}

\begin{figure}[H]
	\centering
	\includegraphics[page=10]{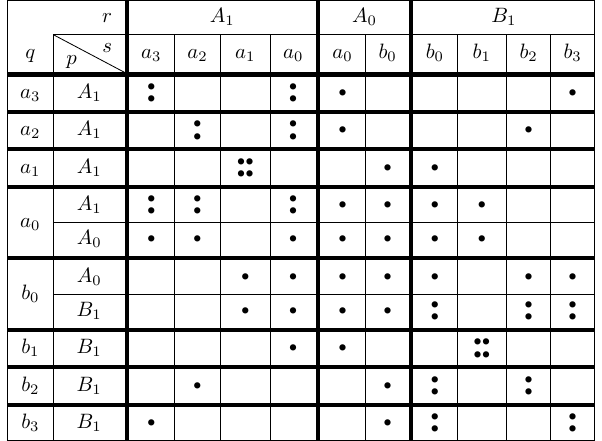}
	\caption{Medium double broom - structure of $v$}
\end{figure}
\NiceMatrixOptions{cell-space-limits = 2pt}
We will outline how to determine the $3 \times 3$ blocks as the remaining entries can be determined using the bi-unitary condition. Consider the $3\times 3$ blocks:
$$u^{(A_1,a_0)}=
\begin{pNiceArray}{W{c}{1cm}c|c}[first-row,first-col]
       & \Block{1-2}{A_1}& & A_0\\
    a_3& \Block{1-2}{\alpha_3}& &x_3 \\ \hline
    a_2& \Block{1-2}{\alpha_2}& &x_2 \\ \hline
    a_0& \Block{1-2}{\alpha_0}& &x_0
\end{pNiceArray},\quad
v^{(a_0,A_1)}=
\begin{pNiceArray}{c|c|c}[first-col,first-row]
&a_3&a_2&a_0\\
\Block{2-1}{A_1}&\Block{2-1}{\beta_3^t}&\Block{2-1}{\beta_2^t}&\Block{2-1}{\beta_0^t}\\
& & &\\ \hline
A_0&y_3 &y_2 &y_0\\
\end{pNiceArray}
$$

where $\alpha_i,\beta_j\in M_{1\times 2}(\R)$.
and observe that $v^{(a_i,A_0)}$ and $u^{(A_0,a_i)}$ are $1\times 1$ blocks for $i=2,3$. Hence we can assume they are both equal to $1$. Using the bi-unitary condition we obtain $y_i=x_i=\sqrt{\frac{\lambda(a_i)\lambda(A_0)}{\lambda(A_1)\lambda(a_0)}}$ for $i=2,3$ and therefore
{$$u^{(A_1,a_0)}=
\begin{pNiceArray}{W{c}{1cm}c|c}[first-row,first-col]
       & \Block{1-2}{A_1}& & A_0\\
    a_3& \Block{1-2}{\alpha_3}& &\sqrt{\frac{3-\sqrt{3}}{3}} \\ \hline
    a_2& \Block{1-2}{\alpha_2}& &\sqrt{\frac{3-\sqrt{3}}{3}} \\ \hline
    a_0& \Block{1-2}{\alpha_0}& &x_0
\end{pNiceArray},\quad 
v^{(a_0,A_1)}=
\begin{pNiceArray}{c|c|c}[first-col,first-row]
&a_3&a_2&a_0\\
\Block{2-1}{A_1}&\Block{2-1}{\beta_3^t}&\Block{2-1}{\beta_2^t}&\Block{2-1}{\beta_0^t}\\
& & &\\ \hline
A_0&\sqrt{\frac{3-\sqrt{3}}{3}} &\sqrt{\frac{3-\sqrt{3}}{3}} &y_0\\
\end{pNiceArray}.$$}

Since the last row (respective column) of $u^{(A_1,a_0)}$ (respective $v^{(a_0,A_1)}$) must have norm 1, we can assume $y_0=x_0=\sqrt{\frac{-3+2\sqrt{3}}{3}}$. Because the rows of $u^{(A_1,a_0)}$ form an orthonormal basis of $\C^3$ we have
\begin{align*}
    \pint{\alpha_3}{\alpha_2}&=-\frac{3-\sqrt{3}}{3} & \|\alpha_3\|^2=\|\alpha_2\|^2&=\frac{\sqrt{3}}{3}\\
    \pint{\alpha_3}{\alpha_0}=\pint{\alpha_3}{\alpha_0}&=-\frac{-5+3\sqrt{3}}{3} & \|\alpha_0\|^2&=\frac{8-3\sqrt{3}}{3}\\
 \end{align*}
 We now look at the $2\times 2$ blocks $v^{(a_i,A_1)}$ and $u^{(A_1,a_i)}$ for $i=2,3$. Using the bi-unitary condition we get 
 
 {$$v^{(a_3,A_1)}=
 \begin{pNiceArray}{c|c}[first-row,first-col]
                & a_3   & a_0\\
\Block{2-1}{A_1}& \ast  & \Block{2-1}{\sqrt{\lambda(a_0)}\alpha_3^t}\\
                & \ast  &
 \end{pNiceArray},\quad
 v^{(a_2,A_1)}=
 \begin{pNiceArray}{c|c}[first-row,first-col]
                & a_2   & a_0\\
\Block{2-1}{A_1}& \ast  & \Block{2-1}{\sqrt{\lambda(a_0)}\alpha_2^t}\\
                & \ast  &
 \end{pNiceArray}
 $$}

 Since both of these have to be unitary, they can be chosen to be of the form:

{$$v^{(a_3,A_1)}=\sqrt{\lambda(a_0)}
 \begin{pNiceArray}{c|c}
 \alpha_{32} & \alpha_{31}\\
 -\alpha_{31} & \alpha_{32}
 \end{pNiceArray},\quad
 v^{(a_2,A_1)}=\sqrt{\lambda(a_0)}
 \begin{pNiceArray}{c|c}
 \alpha_{22} & \alpha_{21}\\
 -\alpha_{21} & \alpha_{22}
 \end{pNiceArray}
$$}

Once again, using the bi-unitary condition and the fact that $u^{(A_1,a_i)}$ for $i=2,3$ are unitaries, we have
{$$u^{(A_1,a_3)}=\sqrt{\lambda(a_0)}
 \begin{pNiceArray}{cc}
 \alpha_{32} & -\alpha_{31}\\ \hline
  \alpha_{31} & \alpha_{32}
 \end{pNiceArray},\quad
 u^{(A_1,a_2)}=\sqrt{\lambda(a_0)}
 \begin{pNiceArray}{cc}
 \alpha_{22} & -\alpha_{21}\\ \hline
 \alpha_{21} & \alpha_{22}
 \end{pNiceArray}
$$}

The bi-unitary condition applied to the blocks $u^{(A_1,a_i)}_{a_0,A_1}$ for $i=0,2,3$ imply that $\beta_i=\alpha_i$. We can now conclude that any choice of $\alpha_i$'s which make $u^{(A_1,a_0)}$ unitary will do the same for $v^{(a_0,A_1)}$. By completing the last row of $u^{(A_1,a_0)}$ to an orthonormal basis we obtain such a choice. A similar reasoning works for the blocks $u^{(B_1,b_0)}$ and $v^{(b_0,B_1)}$. 

The remaining $2\times 2$ blocks are easily obtained by using the bi-unitary condition. Hence we have a bi-unitary connection on the medium double broom.

\subsection{Large double broom}

For the ``large double broom'' we will use the following labeling of the vertices:

\begin{center}
	$\Gamma=$ \begin{tikzpicture}[main/.style={circle,fill=black,inner sep=0pt,minimum size=1.5mm},
        sub/.style={circle,fill=white,draw=black,inner sep=0pt,minimum size=1.5mm},baseline=-2]
	\node[sub,label=left:{$a_2$}] at (0,0.5) (A0) {};
	\node[sub,label=left:{$a_3$}] at (0,0) (B0) {};
	\node[sub,label=left:{$a_4$}] at (0,-0.5) (C0) {};
	\node[main,label=above:{$A_2$}] at (1,0) (A1) {};
	\node[sub,label=above:{$a_1$}] at (2,0) (A2) {};
	\node[main,label=above:{$A_1$}] at (3,0) (A3) {};
	\node[sub,label=above:{$a_0$}] at (4,0) (A4) {};
	\node[main,label=above:{$A_0$}] at (5,0) (A5) {};
	\node[sub,label=above:{$b_0$}] at (6,0) (A6) {};
	\node[main,label=above:{$B_1$}] at (7,0) (A7) {};
	\node[sub,label=above:{$b_1$}] at (8,0) (A8) {};
	\node[main,label=above:{$B_2$}] at (9,0) (A9) {};
	\node[sub,label=right:{$b_2$}] at (10,0.5) (A10) {};
	\node[sub,label=right:{$b_3$}] at (10,0) (B10) {};
	\node[sub,label=right:{$b_4$}] at (10,-0.5) (C10) {};
	
	\draw (B0) -- (A1) -- (C0);
	\draw (B10) -- (A9) -- (C10);
	\draw (A0) -- (A1) -- (A2) -- (A3) -- (A4) -- (A5) -- (A6) -- (A7) -- (A8) -- (A9) -- (A10);
	\end{tikzpicture}
\end{center}

and the following Perron-Frobenius eigenvector for the graph $\Gamma$:

\begin{center}
	$\lambda(\cdot)=$\begin{tikzpicture}[main/.style={circle,fill=black,inner sep=0pt,minimum size=1.5mm},
        sub/.style={circle,fill=white,draw=black,inner sep=0pt,minimum size=1.5mm},baseline=-2]
	\node[sub,label=left:{$1$}] at (0,0.5) (A0) {};
	\node[sub,label=left:{$1$}] at (0,0) (B0) {};
	\node[sub,label=left:{$1$}] at (0,-0.5) (C0) {};
	\node[main,label=above:{$\frac{5+\sqrt{17}}{2}$}] at (1,0) (A1) {};
	\node[sub,label=below:{$\frac{-1+\sqrt{17}}{2}$}] at (2,0) (A2) {};
	\node[main,label=above:{$\frac{1+\sqrt{17}}{2}$}] at (3,0) (A3) {};
	\node[sub,label=below:{$1$}] at (4,0) (A4) {};
	\node[main,label=above:{$2$}] at (5,0) (A5) {};
	\node[sub,label=below:{$1$}] at (6,0) (A6) {};
	\node[main,label=above:{$\frac{1+\sqrt{17}}{2}$}] at (7,0) (A7) {};
	\node[sub,label=below:{$\frac{-1+\sqrt{17}}{2}$}] at (8,0) (A8) {};
	\node[main,label=above:{$\frac{5+\sqrt{17}}{2}$}] at (9,0) (A9) {};
	\node[sub,label=right:{$1$}] at (10,0.5) (A10) {};
	\node[sub,label=right:{$1$}] at (10,0) (B10) {};
	\node[sub,label=right:{$1$}] at (10,-0.5) (C10) {};
	
	\draw (B0) -- (A1) -- (C0);
	\draw (B10) -- (A9) -- (C10);
	\draw (A0) -- (A1) -- (A2) -- (A3) -- (A4) -- (A5) -- (A6) -- (A7) -- (A8) -- (A9) -- (A10);
	\end{tikzpicture}
\end{center}

\begin{figure}[H]
\centering
\includegraphics[page=11,scale=0.75]{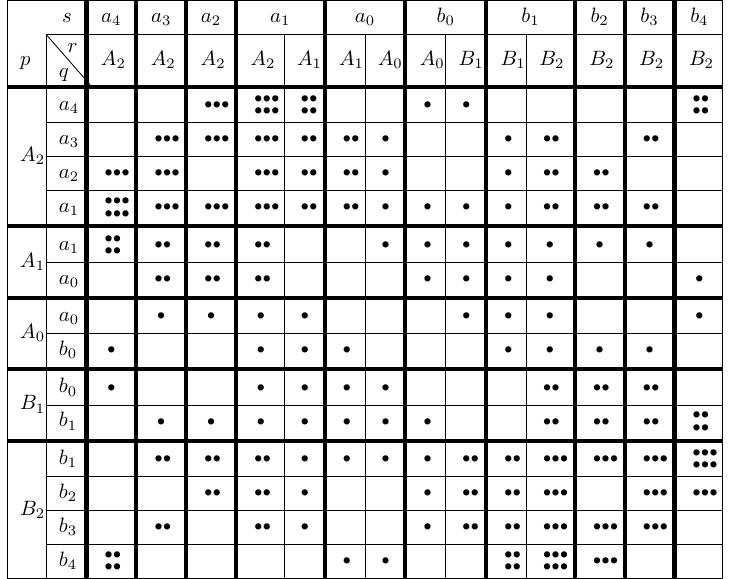}
\caption{Large double broom - structure of $u$}
\label{ublocks}
\end{figure}

\begin{figure}[H]
\centering
\includegraphics[page=12,scale=0.75]{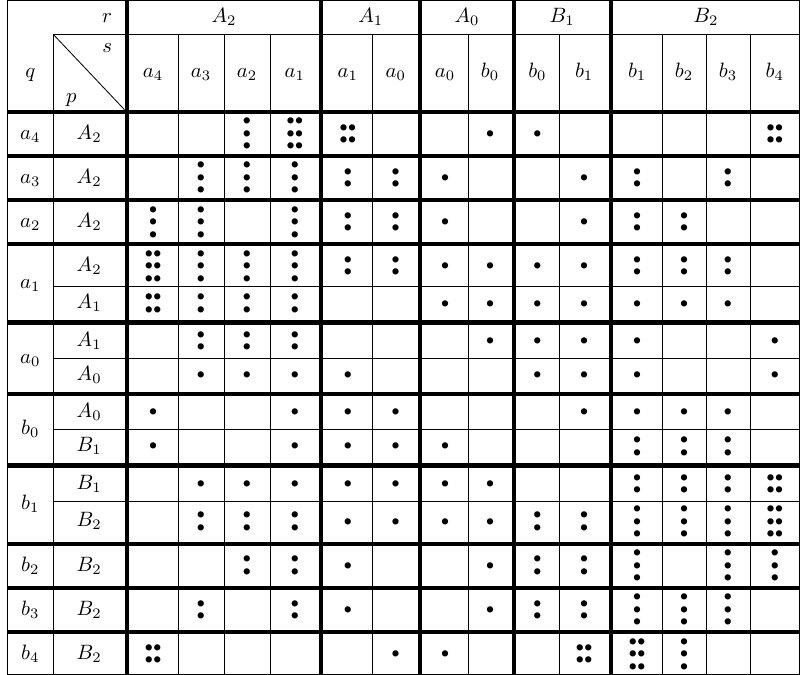}
\caption{Large double broom - structure of $v$}
\label{vblocks}
\end{figure}

Using the same notation as above, we have
\begin{align*}
G &= \begin{pmatrix}
1 & 1 & 1 & 1 & 0 & 0 & 0 & 0 & 0 & 0\\
0 & 0 & 0 & 1 & 1 & 0 & 0 & 0 & 0 & 0\\
0 & 0 & 0 & 0 & 1 & 1 & 0 & 0 & 0 & 0\\
0 & 0 & 0 & 0 & 0 & 1 & 1 & 0 & 0 & 0\\
0 & 0 & 0 & 0 & 0 & 0 & 1 & 1 & 1 & 1 
\end{pmatrix}
\end{align*}
for the vertical inclusions, and 
\begin{align*}
H &= \begin{pmatrix}
3 & 2 & 1 & 1 & 2\\
2 & 0 & 1 & 1 & 2\\
1 & 1 & 0 & 1 & 1\\
1 & 1 & 1 & 0 & 2\\
2 & 1 & 0 & 2 & 3\\
\end{pmatrix}, & K &= \begin{pmatrix}
0 & 0 & 1 & 2 & 0 & 1 & 0 & 0 & 0 & 2\\
0 & 1 & 1 & 1 & 1 & 0 & 1 & 0 & 1 & 0\\
1 & 1 & 0 & 1 & 1 & 0 & 1 & 1 & 0 & 0\\
2 & 1 & 1 & 1 & 1 & 1 & 1 & 1 & 1 & 0\\
0 & 1 & 1 & 1 & 0 & 1 & 1 & 0 & 0 & 1\\
1 & 0 & 0 & 1 & 1 & 0 & 1 & 1 & 1 & 0\\
0 & 1 & 1 & 1 & 1 & 1 & 1 & 1 & 1 & 2\\
0 & 0 & 1 & 1 & 0 & 1 & 1 & 0 & 1 & 1\\
0 & 1 & 0 & 1 & 0 & 1 & 1 & 1 & 1 & 0\\
2 & 0 & 0 & 0 & 1 & 0 & 2 & 1 & 0 & 0
\end{pmatrix}
\end{align*}
for the horizontal ones. Since
\begin{align*}
    GK= \begin{pmatrix}
3 & 3 & 3 & 5 & 3 & 2 & 3 & 2 & 2 & 2\\
2 & 2 & 2 & 2 & 1 & 2 & 2 & 1 & 1 & 1\\
1 & 1 & 1 & 2 & 1 & 1 & 2 & 1 & 1 & 1\\
1 & 1 & 1 & 2 & 2 & 1 & 2 & 2 & 2 & 2\\
2 & 2 & 2 & 3 & 2 & 3 & 5 & 3 & 3 & 3
\end{pmatrix}
\end{align*}
we deduce that $u$ has two $5\times 5$ blocks, ten $3\times 3$ blocks, twenty-two $2\times2$ blocks and sixteen $1\times 1$ blocks, same for $v$. Figure \ref{ublocks} and figure \ref{vblocks} (found at the end of this section) describe the block structure for both $u$ and $v$. The main difficulty in computing the bi-unitary connection will be to determine the entries of the four $5\times 5$ blocks: $u^{(A_2,a_1))}$, $u^{(B_2,b_1)}$, $v^{(a_1,A_2)}$ and $v^{(b_1,B_2)}$. \label{reftotables}

Note that in Figure \ref{ublocks}, $u^{(A_2,a_1))}$ and $u^{(B_2,b_1)}$ are in a symmetric position, hence we will try to assume we can obtain the entries of $u^{(B_2,b_1)}$ by ```reflecting'' the entries of $u^{(A_2,a_1)}$. 

Consider
\begin{equation*}
u^{(A_2,a_1)}=\begin{pNiceArray}{ccc|cc}[first-row,first-col]
    & \Block{1-3}{A_2}       &    & &\Block{1-2}{A_1} &  \\
\Block{2-1}{a_4} & \gamma_{11} & \gamma_{12} & \gamma_{13}  & w_{11} & w_{12}      \\
    & \gamma_{21} & \gamma_{22} & \gamma_{23} & w_{21} & w_{22}      \\
\hline
a_3 & \xi_{11} & \xi_{12} & \xi_{13} & \alpha_{11} & \alpha_{12} \\
\hline
a_2 & \xi_{21} & \xi_{22} & \xi_{23} & \alpha_{21} & \alpha_{22} \\
\hline
a_1 & \xi_{31} & \xi_{32} & \xi_{33} & \alpha_{31} & \alpha_{32}
\end{pNiceArray}
\end{equation*}
and set $\gamma_i=(\gamma_{i1},\gamma_{i2},\gamma_{i3})$, $\xi_i=(\xi_{i1},\xi_{i2},\xi_{i3})$, $w_i=(w_{i1},w_{i2})$ and $\alpha_i=(\alpha_{i1},\alpha_{i2})$. From the bi-unitary condition we have
$$
v^{(a_4,A_1)}=\sqrt{\frac{8}{7-\sqrt{17}}}
\begin{pNiceMatrix}
\Block{2-1}{w_1^t} & \Block{2-1}{w_2^t}\\
 &
\end{pNiceMatrix}
$$
which implies that $\|w_i\|=\sqrt{\frac{7-\sqrt{17}}{8}}$ and $\pint{w_1}{w_2}=0$ as $v^{(a_4,A_1)}$ has to be unitary. Since the first two rows of $u^{(A_2,a_1)}$ are orthonormal, we also get $\|\gamma_i\|=\sqrt{\frac{1+\sqrt{17}}{8}}$ and $\pint{\gamma_1}{\gamma_2}=0$. Without loss of generality we can suppose $v^{(a_4,A_1)}=\ID_2$ and $\gamma_i=\sqrt{\frac{1+\sqrt{17}}{8}} e_i$, where $e_i$ is the canonical $i$-th basis vector in $\C^3$. Thus, we have the following entries:
\begin{equation*}
u^{(A_2,a_1)}=\begin{pNiceArray}{ccc|cc}[first-row,first-col]
    & \Block{1-3}{A_2}       &    & &\Block{1-2}{A_1} &  \\
\Block{2-1}{a_4} & \sqrt{\frac{1+\sqrt{17}}{8}} & 0 & 0  & \sqrt{\frac{7-\sqrt{17}}{8}} & 0      \\
    & 0 & \sqrt{\frac{1+\sqrt{17}}{8}} & 0 & 0 & \sqrt{\frac{7-\sqrt{17}}{8}}      \\
\hline
a_3 & \xi_{31} & \xi_{32} & \xi_{33} & \alpha_{31} & \alpha_{32} \\
\hline
a_2 & \xi_{21} & \xi_{22} & \xi_{23} & \alpha_{21} & \alpha_{22} \\
\hline
a_1 & \xi_{11} & \xi_{12} & \xi_{13} & \alpha_{11} & \alpha_{12}
\end{pNiceArray}
\end{equation*}
Now, from the orthogonality between the rows of $u^{(A_2,a_1)}$ we obtain
\begin{align}
    \xi_{ij}&=-\sqrt{\frac{7-\sqrt{17}}{1+\sqrt{17}}}\alpha_{ij} \label{xivalpha}
\end{align}
for any $i=1,2,3$ and $j=1,2$. Since every row must have norm $1$ we have the identity
\begin{equation}\label{xi3eq}
    1=\xi_{i3}^2+\frac{8}{1+\sqrt{17}}\|\alpha_i\|^2.
\end{equation}
From the bi-unitary condition we obtain
\begin{align*}
    v^{(a_i,A_1)}=\sqrt{\frac{\lambda(A_2)\lambda(a_1)}{\lambda(a_i)\lambda(A_1)}}
    \begin{pNiceArray}{c|c}
    \alpha_{i1} & \ast\\
    \alpha_{i2} & \ast
    \end{pNiceArray}
\end{align*}
for $i=1,2,3$. Since $v^{a_i,A_1}$ has to be unitary, we can assume it is of the form
\begin{align*}
    v^{(a_i,A_1)}=\sqrt{\frac{\lambda(A_2)\lambda(a_1)}{\lambda(a_i)\lambda(A_1)}}
    \begin{pNiceArray}{c|c}
    \alpha_{i1} & -\alpha_{i2}\\
    \alpha_{i2} & \alpha_{i1}
    \end{pNiceArray}
\end{align*}
and $\|\alpha_i\|^2=\frac{\lambda(a_i)\lambda(A_1)}{\lambda(A_2)\lambda(a_1)}$. In particular $\|\alpha_2\|^2=\|\alpha_3\|^2=\frac{7-\sqrt{17}}{8}$ and $\|\alpha_1\|^2=\frac{-3+\sqrt{17}}{2}$. Using this in equation \ref{xi3eq} we obtain that $\xi_{23}^2=\xi_{33}^2=\frac{5-\sqrt{17}}{2}$ and $\xi_{13}^2=-4+\sqrt{17}$. 

We now will determine the entries $\alpha_{ij}$. Consider the block $u^{(A_2,a_0)}$. Using the bi-unitary condition we get

\begin{equation*}
    u^{(A_2,a_0)}=
    \begin{pNiceArray}{W{c}{1.5cm}c|c}[first-row,first-col]
        & \Block{1-2}{A_1}  &               & A_0\\
    a_3 & \Block[c]{1-2}{\sqrt{\frac{\lambda(a_1)}{\lambda(a_0)}}\Tilde{\alpha_3}}&   & \ast\\
    \hline
    a_2 & \Block{1-2}{\sqrt{\frac{\lambda(a_1)}{\lambda(a_0)}}\Tilde{\alpha_2}}&    & \ast\\
    \hline
    a_1 & \Block{1-2}{\sqrt{\frac{\lambda(a_1)}{\lambda(a_0)}}\Tilde{\alpha_1}}&    & \ast
    \end{pNiceArray}
\end{equation*}
where $\Tilde{\alpha_i}=(-\alpha_{i,2},\alpha_{i,1})$. Since $v^{(a_i,A_0)}$ are $1\times 1$ blocks for $i=1,2,3$, we can suppose $v^{(a_2,A_0)}=v^{(a_3,A_0)}=1$ and 
$v^{(a_1,A_0)}=-1$. The bi-unitary condition then implies
\begin{equation}\label{uA2a0}
    u^{(A_2,a_0)}=
    \begin{pNiceArray}{W{c}{1.5cm}c|c}[first-row,first-col]
        & \Block{1-2}{A_1}  &               & A_0\\
    a_3 & \Block[c]{1-2}{\sqrt{\frac{\lambda(a_1)}{\lambda(a_0)}}\Tilde{\alpha_3}}&   & \sqrt{\frac{5-\sqrt{17}}{2}}\\
    \hline
    a_2 & \Block{1-2}{\sqrt{\frac{\lambda(a_1)}{\lambda(a_0)}}\Tilde{\alpha_2}}&    & \sqrt{\frac{5-\sqrt{17}}{2}}\\
    \hline
    a_1 & \Block{1-2}{\sqrt{\frac{\lambda(a_1)}{\lambda(a_0)}}\Tilde{\alpha_1}}&    & -\sqrt{-4+\sqrt{17}}
    \end{pNiceArray}
\end{equation}
Since $u^{(A_2,a_0)}$ is unitary, its columns must form an orthonormal basis of $\C^3$. Given the last column, we can complete it to an orthonormal basis:
$$
x_1=\begin{pNiceMatrix}
    \frac{\sqrt{5-\sqrt{17}}}{2}\\
    0\\
    \frac{\sqrt{-1+\sqrt{17}}}{2}
\end{pNiceMatrix},\;
x_2=\begin{pNiceMatrix}
    -\sqrt{\frac{-11+3\sqrt{17}}{4}}\\
    \sqrt{\frac{-3+\sqrt{17}}{2}}\\
    \frac{\sqrt{21-5\sqrt{17}}}{2}
\end{pNiceMatrix},\;
x_3=\begin{pNiceMatrix}
    \sqrt{\frac{5-\sqrt{17}}{2}}\\
    \sqrt{\frac{5-\sqrt{17}}{2}}\\
    -\sqrt{-4+\sqrt{17}}
\end{pNiceMatrix}
$$
We can now assume that
\begin{equation}\label{uA2a0param}
    u^{(A_2,a_0)}=
    \begin{pNiceArray}{c:c|c}[first-row,first-col]
        & \Block{1-2}{A_1}  &               & A_0\\
    a_3 & \Block{3-1}{q_1(s)}& \Block{3-1}{q_2(s)}  & \sqrt{\frac{5-\sqrt{17}}{2}}\\
    \cline{3-3}
    a_2 & &    & \sqrt{\frac{5-\sqrt{17}}{2}}\\
    \cline{3-3}
    a_1 & &    & -\sqrt{-4+\sqrt{17}}
    \end{pNiceArray}
\end{equation}
where $q_1(s)=\cos(s)x_1+\sin(s)x_2$ and $q_2(s)=-\sin(s)x_1+\cos(s)x_2$. Combining \ref{uA2a0} and \ref{uA2a0param} we obtain $\alpha_{ij}$ as a function of $s$. From \ref{xivalpha} we also obtain $\xi_{ij}$ as a function of $s$ for $i=1,2,3$ and $j=1,2$. 

The fact that the rows of $u^{(A_2,a_1)}$ are orthogonal and \ref{xivalpha} imply that
\begin{equation}\label{xii3eq}
0=\pint{\xi_i}{\xi_j}+\pint{\alpha_i}{\alpha_j}=\xi_{i3}\xi_{j3}+\frac{8}{1+\sqrt{17}}\pint{\alpha_i}{\alpha_j}    
\end{equation}
for $i\neq j$. From the unitarity of \ref{uA2a0} we get
\begin{align}
    \pint{\alpha_2}{\alpha_3}=\pint{\tilde{\alpha_2}}{\tilde{\alpha_3}}&=-\frac{\lambda(a_0)}{\lambda(a_1)}\frac{5-\sqrt{17}}{2} \label{pintalpha23}\\
    \pint{\alpha_i}{\alpha_1}=\pint{\tilde{\alpha_i}}{\tilde{\alpha_1}}&=\frac{\lambda(a_0)}{\lambda(a_1)}\sqrt{\frac{5-\sqrt{17}}{2}}\sqrt{-4+\sqrt{17}},\quad i=2,3\label{pintalpha1i}
\end{align}
From \ref{xii3eq}, \ref{pintalpha23}, \ref{pintalpha1i}, and the fact that $\xi_{23}^2=\xi_{33}^2=\frac{5-\sqrt{17}}{2}$ and $\xi_{13}^2=-4+\sqrt{17}$ we compute
$$\xi_{23}=\xi_{33}=\pm\sqrt{\frac{5-\sqrt{17}}{2}},\quad \xi_{13}=\mp\sqrt{-4+\sqrt{17}}.$$
We will work with $\xi_{13}=-\sqrt{-4+\sqrt{17}}$. 

Note that all entries of $u^{(A_2,a_0)}$ depend on a single parameter $s$. Moreover, since we chose the entries so that the rows form an orthonormal basis of $\C^5$, we conclude that it has to be unitary. To determine the possible values $s$ can take we need to investigate the blocks $v^{(a_i,A_2)}$ for $i=1,\dots,4$ and make sure they can be completed to unitaries using the entries from $u^{(A_2,a_0)}$. It is important to note that $v^{(a_1,A_2)}$ is also a $5\times 5$ block and hence will be the main source of obstructions for $s$.

Using the bi-unitary condition for the blocks $v^{(a_i,A_2)}$ with $i=2,3,4$, we have the following:
\begin{align*}
    v^{(a_4,A_2)}&=
    \begin{pNiceArray}{c|cc}[first-row,first-col]
                    & a_2  &\Block{1-2}{a_1}&\\
    \Block{3-1}{A_2}& \ast & 1 & 0\\
                    & \ast & 0 & 1\\
                    & \ast & 0 & 0
    \end{pNiceArray},&
    v^{(a_3,A_2)}&=
    \begin{pNiceArray}{c|c|W{c}{1.5cm}}[first-row,first-col]
                    & a_3  &a_2&a_1\\
    \Block{3-1}{A_2}& \Block{3-1}{g_3^t}& \Block{3-1}{g_2^t}&\Block{3-1}{\sqrt{\frac{\lambda(a_1)}{\lambda(a_3)}}\xi_3^t} \\
                    &  &  & \\
                    &  &  & 
    \end{pNiceArray},&
    v^{(a_2,A_2)}&=
    \begin{pNiceArray}{c|c|W{c}{1.5cm}}[first-row,first-col]
                    & a_4  &a_3&a_1\\
    \Block{3-1}{A_2}& \Block{3-1}{f_4^t} & \Block{3-1}{f_3^t} &\Block{3-1}{\sqrt{\frac{\lambda(a_1)}{\lambda(a_2)}}\xi_2^t} \\
                    &  &  & \\
                    &  &  & 
    \end{pNiceArray}
\end{align*}
where the $f_i$'s and $g_i's$ are unit vectors in $\C^3$. Since each of the above matrices has to be unitary we can choose $v^{(a_4,A_2)}=\begin{psmallmatrix}
    0 & 1 & 0\\
    0 & 0 & 1\\
    1 & 0 & 0
\end{psmallmatrix}$. To determine the remaining entries we need to look at the blocks $u^{(A_2,a_i)}$ for $i=2,3,4$. Once again, using the bi-unitary condition we have
\begin{align*}
    u^{(A_2,a_4)}&=
    \begin{pNiceArray}{ccc}[first-row,first-col]
                        &\Block{1-3}{A_2}   & &\\
        a_2             &f_{41}&f_{42} &f_{43}\\
        \hline
        \Block{2-1}{a_1}&\eta_{11}&\eta_{12}&\eta_{13}\\
                        &\eta_{21}&\eta_{22}&\eta_{23}
    \end{pNiceArray},&
    u^{(A_2,a_3)}&=
    \begin{pNiceArray}{ccc}[first-row,first-col]
            &\Block{1-3}{A_2}   & &\\
        a_3 &\Block{1-3}{g_3}& &\\
        \hline
        a_2 &\Block{1-3}{f_3}& &\\
        \hline
        a_1 &\Block{1-3}{h_1}& &
    \end{pNiceArray},&
    u^{(A_2,a_2)}&=
    \begin{pNiceArray}{ccc}[first-row,first-col]
            &\Block{1-3}{A_2}   &       &\\
        a_4 & 0                 & 0     & 1\\
        \hline
        a_3 & \Block{1-3}{g_2}  &       &\\
        \hline
        a_1 & \Block{1-3}{g_1}  &   & 
    \end{pNiceArray}
\end{align*}
Since $u^{(A_2,a_2)}$ is unitary, we can assume it is of the form
\begin{equation*}
u^{(A_2,a_2)}=
    \begin{pNiceArray}{ccc}[first-row,first-col]
            &\Block{1-3}{A_2}   &       &\\
        a_4 & 0                 & 0     & 1\\
        \hline
        a_3 & \cos(t)           &\sin(t)&0\\
        \hline
        a_1 & -\sin(t)          &\cos(t)&0
\end{pNiceArray}    
\end{equation*}
and consequently $g_2$ is a function of $t$. Using the fact that $v^{(a_3,A_2)}$ is unitary, we have that $g_3=\sqrt{\frac{\lambda(a_1)}{\lambda(a_3)}}g_2\times \xi_3$ provided $g_2 \perp \xi_3$ (here $\times$ denotes the cross product). Similarly we have that $f_3\perp g_3$ and $f_3\perp \xi_2$, therefore 
\begin{equation*}
    f_3=\frac{g_3\times \xi_2}{\|g_3\times \xi_2\|}
\end{equation*}
and consequently $h_1=g_3\times f_3$. Note that $g_2$ depends on $t$ and $\xi_3$ depends on $s$, in fact $$\pint{g_2
}{\xi_3}=r_1 \cos(t+s)+r_2 \sin(t+s)$$ where $r_1=\sqrt{-4+\sqrt{17}}$ and $r_2=\frac{\sqrt{13-3\sqrt{17}}}{2}$. We can plot this function:
\begin{figure}[H]
    \centering
    \includegraphics[scale=0.5]{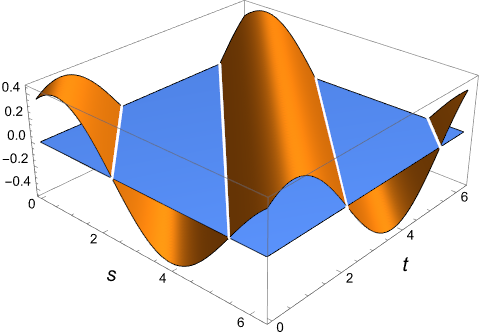}
    \label{g2xi3}
    \caption{$\pint{g_2}{\xi_3}(t,s)$}
\end{figure}
Observer that there are straight lines on which $\pint{g_2}{\xi_3}=0$.

The last set of restrictions will come from looking at the $5\times 5$ block $v^{(a_1,A_2)}$. Using the bi-unitary condition we have
\begin{equation*}
    v^{(a_1,A_2)} = 
    \begin{pNiceArray}{cc|W{c}{1.5cm}|W{c}{1.5cm}|c}[first-row,first-col]
                    &\Block{1-2}{a_4}       & & a_3 & a_2 & a_1\\
    \Block{3-1}{A_2}&\Block{3-1}{\zeta_1^t}  &\Block{3-1}{\zeta_2^t}&\Block{3-1}{\sqrt{\frac{1}{\lambda(a_1)}}h_1^t}&\Block{3-1}{\sqrt{\frac{1}{\lambda(a_1)}}g_1^t}&\Block{3-1}{\xi_1^t}\\
    & & & & &\\
    & & & & &\\
    \hline
    \Block{2-1}{A_1}&\Block{2-1}{w_3^t}&\Block{2-1}{w_4^t}&\Block{2-1}{\beta_3^t}&\Block{2-1}{\beta_2^t}&\Block{2-1}{\beta_1^t}\\
    & & & & &
    \end{pNiceArray}
\end{equation*}
where $\zeta_i=\sqrt{\frac{1}{\lambda(a_1)}}\eta_i$. This block is similar to the block $u^{(A_2,a_1)}$, in particular the orthogonality of the last three columns implies
\begin{align}
    \frac{1}{\lambda(a_1)}\pint{h_1}{g_1}+\pint{\beta_3}{\beta_2}&=0, \label{eq1}\\
    \sqrt{\frac{1}{\lambda(a_1)}}\pint{h_1}{\xi_1}+\pint{\beta_3}{\beta_1}&=0, \label{eq2}\\
    \sqrt{\frac{1}{\lambda(a_1)}}\pint{g_1}{\xi_1}+\pint{\beta_2}{\beta_1}&=0. \label{eq3}
\end{align}
If we now suppose $u^{(A_0,a_3)}=-1$, $u^{(A_0,a_2)}=1$ and $u^{(A_0,a_1)}_{a_0,A2}>0$ then we obtain
\begin{equation*}
    v^{(a_0,A_2)}=
    \begin{pNiceArray}{c|c|c}[first-row,first-col]
        & a_3 & a_2 & a_1\\
        \Block{2-1}{A_1}&\ast &\ast &\ast\\
        &\ast &\ast &\ast\\
        \hline
        A_0&\sqrt{\frac{5-\sqrt{17}}{2}} & -\sqrt{\frac{5-\sqrt{17}}{2}} & \sqrt{-4+\sqrt{17}}
    \end{pNiceArray}
\end{equation*}
and a reasoning similar to the one used in \ref{pintalpha23} and \ref{pintalpha1i} implies that 
\begin{align*}
    \pint{\beta_3}{\beta_2}&=\frac{-3+\sqrt{17}}{4},&
    \pint{\beta_3}{\beta_1}&=\sqrt{\frac{-45+11\sqrt{17}}{16}},&
    \pint{\beta_2}{\beta_1}&=-\sqrt{\frac{-45+11\sqrt{17}}{16}}
\end{align*}
Plugging this into (\ref{eq1}), (\ref{eq2}) and (\ref{eq3}) we obtain the following: 
\begin{align}
    F_{1}(s,t)&=\frac{1}{\lambda(a_1)}\pint{h_1}{g_1}+\frac{-3+\sqrt{17}}{4}=0,\\
    F_{2}(s,t)&=\sqrt{\frac{1}{\lambda(a_1)}}\pint{h_1}{\xi_1}+\sqrt{\frac{-45+11\sqrt{17}}{16}} =0,\\
    F_{3}(s,t)&=\sqrt{\frac{1}{\lambda(a_1)}}\pint{g_1}{\xi_1}-\sqrt{\frac{-45+11\sqrt{17}}{16}} =0.
\end{align}
On the other hand, orthonormality of the first three rows implies
\begin{align*}
    \zeta_{1i}\zeta_{1j}+\zeta_{2i}\zeta_{2j}+\frac{h_{1i} h_{1j}+g_{1i}g_{1j}}{\lambda(a_1)}+\xi_{1i}\xi_{1j}&=0,\quad i\neq j\in \{1,2,3\}\\
    \zeta_{1i}^2+\zeta_{2i}^2+\frac{h_{1i} ^2+g_{1i}^2}{\lambda(a_1)}+\xi_{1i}\xi_{1j}&=1,\quad i\in \{1,2,3\}
\end{align*}
Since $u^{(A_2,a_4)}$ is also orthonormal, we have $f_{4i}f_{4j}+\eta_{1i}\eta_{2j}+\eta_{2i}\eta_{1j}=0$ for $i\neq j$ and $f_{4i}^2+\eta_{1i}^2+\eta_{1i}^2=1$, consequently we obtain
\begin{align*}
    G_{i,j}(s,t)&=\frac{-f_{4i}f_{4j}+h_{1i} h_{1j}+g_{1i}g_{1j}}{\lambda(a_1)}+\xi_{1i}\xi_{1j}=0,\quad i\neq j\in \{1,2,3\}\\
    G_{i,i}(s,t)&=\frac{1-f_{4i}^2+h_{1i} ^2+g_{1i}^2}{\lambda(a_1)}+\xi_{1i}^2-1=0,\quad i\in \{1,2,3\}
\end{align*}

Finding a common root for all $F_i$'s and $G_{i,j}$'s algebraically is very complicated. Below we plot $F_1$, $F_2$ and $F_3$:
\begin{figure}[H]
    \centering
    \includegraphics[scale=0.5]{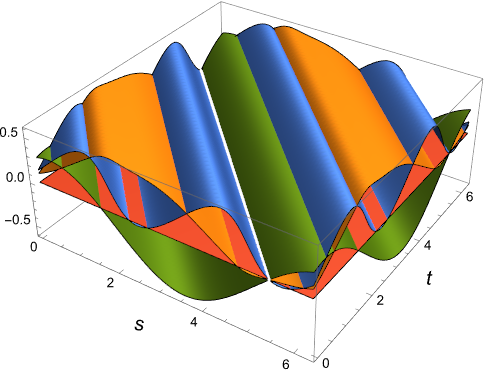}
    \caption{$F_i(s,t)$}
    \label{Fis}
\end{figure}
Observe that the $F_i(s,t)$'s appear to be zero in one of the lines for which $\pint{g_2}{\xi_3}=0$. Moreover, the plot suggests that there is a solution for $s=0$. Solving for $\pint{g_2}{\xi_3}(0,t)=0$ we obtain $$t_0=2\pi -\arctan\left(\frac{\sqrt{-1+\sqrt{17}}}{2}\right)\approx 5.3873.$$
We now consider the plot of $F_i(0,t)$ and $G_{i,j}(0,t)$:
\begin{figure}[H]
    \centering
    \begin{tikzpicture}
    \node at (0,0) {\includegraphics[scale=0.5]{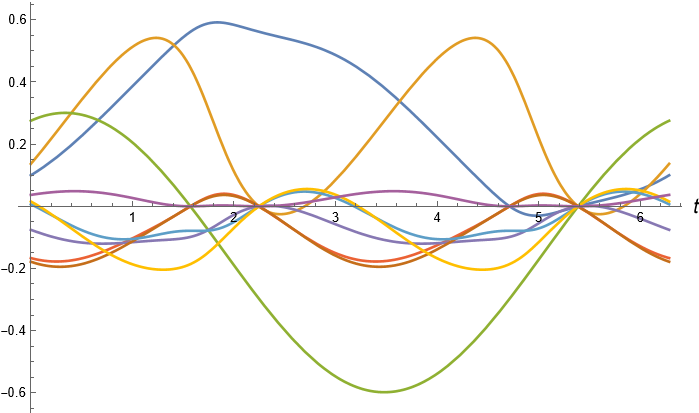}};
    \node[circle,fill=black,inner sep=1pt,label=below:{$t_0$}] at (3.03,0) {};
    \end{tikzpicture}
    \caption{$F_i(0,t)$ and $G_{i,j}(0,t)$}
    \label{FisGijs}
\end{figure}
Observe that they seem to share a single common root at $t_0$. Algebraically we can verify that this is indeed a common root and therefore 
\begin{equation*}
    v^{(a_1,A_2)}=
    \begin{pNiceArray}{cc|c|c|c}[first-row,first-col]
        &\Block{1-2}{a_4} & &a_3 & a_2 & a_1 \\
    \Block{3-1}{A_2}& \Block{3-1}{\zeta_1^t}& \Block{3-1}{\zeta_2^t} & \sqrt{\frac{-169+41\sqrt{17}}{32}} & \sqrt{\frac{-1+\sqrt{17}}{8}} & -\frac{\sqrt{29-7\sqrt{17}}}{2}\\
    & & & \frac{-9+\sqrt{17}}{8} & \frac{1}{2} & \frac{\sqrt{-3+\sqrt{17}}}{2}\\
    & &  & -\sqrt{\frac{31-7\sqrt{17}}{8}} & 0 & -\sqrt{-4+\sqrt{17}}\\ \hline
    \Block{2-1}{A_1}&\Block{2-1}{w_3^t} & \Block{2-1}{w_4^t} & \Block{2-1}{\beta_3^t} & \Block{2-1}{\beta_2^t} & \Block{2-1}{\beta_1^t}\\
    & &  & & & 
    \end{pNiceArray}.
\end{equation*}
To determine the $\zeta_i$'s we observe that any choice of $\eta_i$'s that makes $u^{(A_2,a_4)}$ into a unitary will suffice, for example
$$u^{(A_2,a_4)}=
\begin{pNiceArray}{ccc}[first-row,first-col]
&\Block{1-3}{A_2} & &\\
    a_2 &\sqrt{\frac{-75+19\sqrt{17}}{32}} & \sqrt{\frac{23+\sqrt{17}}{32}} & \frac{5-\sqrt{17}}{4}\\ \hline
    \Block{2-1}{a_1}& \frac{3-\sqrt{17}}{2}& 0 & \sqrt{\frac{-11+3\sqrt{17}}{2}}\\ 
    & -\sqrt{\frac{-101+29\sqrt{17}}{32}}& \frac{-1+\sqrt{17}}{8} & -\sqrt{\frac{31-7\sqrt{17}}{8}}
\end{pNiceArray}.
$$ 

To determine $w_i$'s and the $\beta_i$'s in the block $v^{(a_1,A_2)}$ we complete the first three rows of $v^{(a_1,A_2)}$ to an orthonormal basis of $\C^5$. Once again, the choice of these vectors does not matter as the inner products $\pint{w_i}{w_j}$ and $\pint{\beta_i}{\beta_j}$ are already determined by the rest of the matrix. Finally, the two $5\times 5$ blocks we considered are $u^{(A_2,a_1)}=$
\begin{equation*}
  \begin{pNiceArray}{ccc|cc}[first-row,first-col]
                        & \Block{1-3}{A_2} &  &  & \Block{1-2}{A_1} &  \\
     \Block{2-1}{a_4}&\sqrt{\frac{1+\sqrt{17}}{8}}& 0 & 0 &\sqrt{\frac{7-\sqrt{17}}{8}} &0 \\
 &0 & \sqrt{\frac{1+\sqrt{17}}{8}} & 0 &0&\sqrt{\frac{7-\sqrt{17}}{8}} \\ \hline
 a_3&\sqrt{-4+\sqrt{17}} & \frac{-3+\sqrt{17}}{4} & \sqrt{\frac{5-\sqrt{17}}{2}} & -\frac{\sqrt{5-\sqrt{17}}}{2}&-\sqrt{\frac{-3+\sqrt{17}}{8}} \\ \hline
a_2&-\sqrt{\frac{-19+5\sqrt{17}}{8}} & 0 & \sqrt{\frac{5-\sqrt{17}}{2}}& \sqrt{\frac{7-\sqrt{17}}{8}}& 0 \\ \hline
a_1&-\frac{\sqrt{29-7\sqrt{17}}}{2} & \frac{\sqrt{-3+\sqrt{17}}}{2} & -\sqrt{-4+\sqrt{17}}& \sqrt{\frac{-4+\sqrt{17}}{2}}& -\frac{1}{\sqrt{2}}
  \end{pNiceArray}
\end{equation*}
and $v^{(a_1,A_2)}=$
\begin{equation*}
    \begin{pNiceArray}{cc|c|c|c}[first-row,first-col]
        &\Block{1-2}{a_4} & &a_3 & a_2 & a_1 \\
    \Block{3-1}{A_2}&-\sqrt{\frac{-19+5\sqrt{17}}{8}} & \frac{-9+\sqrt{17}}{8} & \sqrt{\frac{-169+41\sqrt{17}}{32}} & \sqrt{\frac{-1+\sqrt{17}}{8}} & -\frac{\sqrt{29-7\sqrt{17}}}{2}\\
    &0 & \sqrt{\frac{-1+\sqrt{17}}{32}} & \frac{-9+\sqrt{17}}{8} & \frac{1}{2} & \frac{\sqrt{-3+\sqrt{17}}}{2}\\
    &\sqrt{\frac{5-\sqrt{17}}{2}} & -\sqrt{\frac{-11+3\sqrt{17}}{8}} & -\sqrt{\frac{31-7\sqrt{17}}{8}} & 0 & -\sqrt{-4+\sqrt{17}}\\ \hline
    \Block{2-1}{A_1}&\sqrt{\frac{-3+\sqrt{17}}{8}} & -\frac{\sqrt{5-\sqrt{17}}}{2} & \sqrt{\frac{-3+\sqrt{17}}{8}} & 0 & \frac{1}{\sqrt{2}}\\
    &\frac{\sqrt{5-\sqrt{17}}}{2} & \sqrt{\frac{-3+\sqrt{17}}{8}} & \frac{\sqrt{5-\sqrt{17}}}{2} & \sqrt{\frac{7-\sqrt{17}}{8}} & -\sqrt{\frac{-4+\sqrt{17}}{2}}
    \end{pNiceArray}.
\end{equation*}
In a similar manner we obtain the $5\times 5$ blocks $u^{(B_2,b_1)}$ and $v^{(b_1,B_2)}$. The remaining $3\times 3$ blocks: $u^{(A_2,b_1)}$, $u^{(B_2,a_1)}$, $v^{(a_1,B_2)}$ and $v^{(b_1,A_2)}$ work exactly as the $3\times 3$ blocks in the medium broom. The rest of $u$ and $v$ can be easily filled using the bi-unitary condition. In Appendix \ref{appconnections} we present all of the entries of $u$ and $v$ for the bi-unitary connection.

\begin{rem}
From figure (\ref{Fis}) we have that $F_i(s,t)=0$ on the line $s+t=t_0$. Moreover, we can verify that $G_{i,j}(s,t)=0$ on this line and therefore we have a 1-parameter family of solutions for $u$ and $v$. 
\end{rem}

\begin{theo}
    There is an irreducible, hyperfinite $A_\infty$-subfactor with index $\frac{5+\sqrt{17}}{2}$.
\end{theo}
\begin{proof}
    Note that the large double broom is not one of the module graphs in Table \ref{tblr:AHmg} or Table \ref{tblr:AHpmg}. By the embedding theorem \ref{embeddingT} we deduce that the subfactor obtained from this commuting square is not the Asaeda-Haagerup subfactor and consequently has infinite depth. Since the $1$-norm of the first column of $G$ is $1$, Wenzl's criterion implies this subfactor is irreducible. Therefore we have constructed an irreducible hyperfinite subfactor with index $\frac{5+\sqrt{17}}{2}$ and infinite depth. By classification it has to have trivial standard invariant. 
\end{proof}
Next, we list tables that record the structure of the bi-unitary connection for the large double broom (see page \pageref{reftotables}).

\subsection{Quipu}
For the ``quipu'' we will use the following labeling of the vertices:

\begin{center}
	$\Gamma=$\begin{tikzpicture}[main/.style={circle,fill=black,inner sep=0pt,minimum size=1.5mm},
        sub/.style={circle,fill=white,draw=black,inner sep=0pt,minimum size=1.5mm},baseline=-2]
	\node[main,label=above:{$A_3$}] at (-6,0) (A3) {};
    \node[sub,label=below:{$a_4$}] at (-5,0) (a4) {};
    \node[main,label=above:{$A_2$}] at (-4,0) (A2) {};
    \node[sub,label=below:{$a_2$}] at (-3,0) (a2) {};
    \node[main,label=above:{$A_1$}] at (-2,0) (A1) {};
    \node[sub,label=below:{$a_1$}] at (-1,0) (a1) {};
    \node[main,label=above:{$A_0$}] at (0,0) (A0) {};
    \node[sub,label=below:{$a_{-1}$}] at (1,0) (b1) {};
    \node[main,label=above:{$A_{-1}$}] at (2,0) (B1) {};
    \node[sub,label=below:{$a_{-2}$}] at (3,0) (b2) {};
    \node[main,label=above:{$A_{-2}$}] at (4,0) (B2) {};
    \node[sub,label=below:{$a_{-4}$}] at (5,0) (b4) {};
    \node[main,label=above:{$A_{-3}$}] at (6,0) (B3) {};
    \node[sub,label=below:{$a_{3}$}] at (-4,-0.7) (a3) {};
    \node[sub,label=below:{$a_0$}] at (0,-0.7) (a0) {};
    \node[sub,label=below:{$a_{-3}$}] at (4,-0.7) (b3) {};
	
    \draw (A3) -- (a4) -- (A2) -- (a2) -- (A1) -- (a1) -- (A0) -- (b1) -- (B1) -- (b2) -- (B2) -- (b4) -- (B3);
    \draw (A2) -- (a3);
    \draw (A0) -- (a0);
    \draw (B2) -- (b3); 
	\end{tikzpicture}
\end{center}
Since the entries of the Perron-Frobenius eigenvector for the graph $\Gamma$ do not have simple algebraic expressions, we will just use $\lambda(v)$ to refer to them, where $v$ is a vertex of $\Gamma$. In particular, the symmetry of $\Gamma$ implies $\lambda(A_i)=\lambda(A_{-i})$ for $i=1,2,3$ and $\lambda(a_j)=\lambda(a_{-j})$ for $j=1,2,3,4$. Since $(\lambda(A_i))_i$ is a Perron-Frobenius eigenvector of $G^tG$ and $(\lambda(a_j))_j=G(\lambda(A_i))_i$, we also have the following relations between the $\lambda(A_i)'s$ and the $\lambda(a_j)'s$:
\begin{align}
    \lambda(a_0)&=\lambda(A_0),&t \lambda(A_0)&= 2 \lambda(a_1)+\lambda(a_0), \nonumber\\ 
    \lambda(a_1)&=\lambda(A_1)+\lambda(A_0),&t \lambda(A_1)&= \lambda(a_2)+\lambda(a_1), \nonumber\\
    \lambda(a_2)&=\lambda(A_2)+\lambda(A_1),&t \lambda(A_2)&= \lambda(a_4)+\lambda(a_3)+\lambda(a_2), \nonumber\\
    \lambda(a_3)&=\lambda(A_2),&t \lambda(A_3)&= \lambda(a_4),\nonumber\\
    \lambda(a_4)&=\lambda(A_3)+\lambda(A_2), \label{quipueigrel}
\end{align}
where $t$ is the Perron-Frobenius eigenvalue of $G^tG$. Setting $\lambda(A_0)=1$, we can use all but the last of the above relations to obtain
\begin{align}
    \lambda(a_0)&=1, & \lambda(A_0)&=1, \nonumber\\
    \lambda(a_1)&=\frac{t-1}{2}, & \lambda(A_1)&=\frac{t-3}{2}, \nonumber\\
    \lambda(a_2)&=\frac{t^2-4t+1}{2}, & \lambda(A_2)&=\frac{t^2-5t+4}{2}, \nonumber\\
    \lambda(a_3)&=\frac{t^2-5t+4}{2}, & \lambda(A_3)&=\frac{t^3-7t^2+13t-5}{2t}, \nonumber\\
    \lambda(a_4)&=\frac{t^3-7t^2+13t-5}{2}. \label{quipueigvec}
\end{align}
It can be shown that $G^tG$ has minimal polynomial $p(x)=x^3-8x^2+17x-5$, and hence
\begin{equation}\label{quiput^3}
    t^3=8t^2-17t+5.
\end{equation}

Using the same notation as in the previous section we have
\begin{align*}
G &= \begin{pmatrix}
1 & 0 & 0 & 0 & 0 & 0 & 0 & 0 & 0\\
1 & 1 & 1 & 0 & 0 & 0 & 0 & 0 & 0\\
0 & 0 & 1 & 1 & 0 & 0 & 0 & 0 & 0\\
0 & 0 & 0 & 1 & 1 & 1 & 0 & 0 & 0\\
0 & 0 & 0 & 0 & 0 & 1 & 1 & 0 & 0\\
0 & 0 & 0 & 0 & 0 & 0 & 1 & 1 & 1\\
0 & 0 & 0 & 0 & 0 & 0 & 0 & 0 & 1\\
\end{pmatrix}
\end{align*}
for the vertical inclusions, and 
\begin{align*}
H &= \begin{pmatrix}
0 & 0 & 1 & 0 & 1 & 0 & 0\\
0 & 1 & 1 & 2 & 1 & 1 & 0\\
1 & 1 & 1 & 2 & 1 & 1 & 1\\
0 & 2 & 2 & 2 & 2 & 2 & 0\\
1 & 1 & 1 & 2 & 1 & 1 & 1\\
0 & 1 & 1 & 2 & 1 & 1 & 0\\
0 & 0 & 1 & 0 & 1 & 0 & 0\\
\end{pmatrix}, & K &= \begin{pmatrix}
0 & 0 & 1 & 1 & 0 & 1 & 1 & 0 & 0\\
0 & 1 & 0 & 1 & 0 & 1 & 0 & 1 & 0\\
1 & 0 & 1 & 1 & 2 & 1 & 1 & 0 & 1\\
1 & 1 & 1 & 2 & 0 & 2 & 1 & 1 & 1\\
0 & 0 & 2 & 0 & 2 & 0 & 2 & 0 & 0\\
1 & 1 & 1 & 2 & 0 & 2 & 1 & 1 & 1\\
1 & 0 & 1 & 1 & 2 & 1 & 1 & 0 & 1\\
0 & 1 & 0 & 1 & 0 & 1 & 0 & 1 & 0\\
0 & 0 & 1 & 1 & 0 & 1 & 1 & 0 & 0\\
\end{pmatrix}
\end{align*}
for the horizontal ones. Since
\begin{align*}
    GK &= \begin{pmatrix}
0 & 0 & 1 & 1 & 0 & 1 & 1 & 0 & 0\\
1 & 1 & 2 & 3 & 2 & 3 & 2 & 1 & 1\\
2 & 1 & 2 & 3 & 2 & 3 & 2 & 1 & 2\\
2 & 2 & 4 & 4 & 2 & 4 & 4 & 2 & 2\\
2 & 1 & 2 & 3 & 2 & 3 & 2 & 1 & 2\\
1 & 1 & 2 & 3 & 2 & 3 & 2 & 1 & 1\\
0 & 0 & 1 & 1 & 0 & 1 & 1 & 0 & 0
\end{pmatrix} 
\end{align*}
we deduce that $u$ has four $4\times 4$ blocks, eight $3\times 3$ blocks, twenty-one $2\times2$ blocks and twenty $1\times 1$ blocks, same for $v$. Figure \ref{ublocksq} and figure \ref{vblocksq} describe the block structure for both $u$ and $v$.

\begin{figure}[H]
\centering
$$\begin{tblr}{
colspec={cccccc},
column{1-2} = {colsep=4pt},
hlines,
vlines,
hline{3,18} = {7-8,12-13}{2pt},
hline{4,17} = {6-9,11-14}{2pt},
hline{5,16} = {5,15}{2pt},
hline{6,15} = {4-5,10,15-16}{2pt},
hline{7,14} = {3-4,6-14,16-17}{2pt},
hline{8,13} = {5,10,15}{2pt},
hline{9,12} = {3-9,11-17}{2pt},
hline{10,11} = {10}{2pt},
vline{3,18} = {7-8,12-13}{2pt},
vline{4,17} = {6,9-11,14}{2pt},
vline{5,16} = {5-15}{2pt},
vline{6,15} = {4-16}{2pt},
vline{7,14} = {3,17}{2pt},
vline{8,13} = {3-17}{2pt},
vline{9,12} = {3,17}{2pt},
vline{10,11} = {4-16}{2pt},
cell{1}{1} = {r=2,c=2}{c},
cell{1}{3,6,8,11,13,16} = {c=2}{c},
cell{7,12}{1} = {r=2}{c},
cell{4,9,14}{1} = {r=3}{c},
}
\diagboxthree[width=5em,height=5em]{pq}{u}{\rot{rs}}& & \rot{a_4} & & \rot{a_3} & \rot{a_2} &  & \rot{a_1} & & \rot{a_0} & \rot{a_{-1}} & & \rot{a_{-2}} & & \rot{a_{-3}} & \rot{a_{-4}} & \\
& & \rot{A_3} & \rot{A_2} & \rot{A_2} & \rot{A_2} & \rot{A_1} & \rot{A_1} & \rot{A_0} & \rot{A_0} & \rot{A_0} & \rot{A_{-1}} & \rot{A_{-1}} & \rot{A_{-2}} & \rot{A_{-2}} & \rot{A_{-2}} & \rot{A_{-3}}\\
A_3&a_4 & & & & & \bulmat{1}{1} & \bulmat{1}{1} & & & & \bulmat{1}{1} & \bulmat{1}{1} & & & &\\
A_2&a_4 & & & & \bulmat{1}{1} & \bulmat{1}{1} & \bulmat{1}{1} & \bulmat{1}{2} & & \bulmat{1}{2} & \bulmat{1}{1} & \bulmat{1}{1} & \bulmat{1}{1} & & &\\
&a_3 & & & \bulmat{1}{1} & & & \bulmat{1}{1} & \bulmat{1}{2} & & \bulmat{1}{2} & \bulmat{1}{1} & & & \bulmat{1}{1} & &\\
&a_2 & & \bulmat{1}{1} & & \bulmat{1}{1} & \bulmat{1}{1} & \bulmat{1}{1} & \bulmat{1}{2} & \bulmat{2}{2} & \bulmat{1}{2} & \bulmat{1}{1} & \bulmat{1}{1} & \bulmat{1}{1} & & \bulmat{1}{1} &\\
A_1&a_2 & \bulmat{1}{1} & \bulmat{1}{1} & & \bulmat{1}{1} & \bulmat{1}{1} & \bulmat{1}{1} & \bulmat{1}{2} & \bulmat{2}{2} & \bulmat{1}{2} & \bulmat{1}{1} & \bulmat{1}{1} & \bulmat{1}{1} & & \bulmat{1}{1} & \bulmat{1}{1}\\
&a_1& \bulmat{1}{1} & \bulmat{1}{1} & \bulmat{1}{1} & \bulmat{1}{1} & \bulmat{1}{1} & \bulmat{2}{1} & \bulmat{2}{2} & & \bulmat{2}{2} & \bulmat{2}{1} & \bulmat{1}{1} & \bulmat{1}{1} & \bulmat{1}{1} & \bulmat{1}{1} & \bulmat{1}{1}  \\
A_0&a_1 & & \bulmat{1}{2} & \bulmat{1}{2} & \bulmat{1}{2} & \bulmat{1}{2} & \bulmat{2}{2} & \bulmat{2}{2} & & \bulmat{2}{2} & \bulmat{2}{2} & \bulmat{1}{2} & \bulmat{1}{2} & \bulmat{1}{2} & \bulmat{1}{2} & \\
&a_0 & & & & \bulmat{2}{2} & \bulmat{2}{2} & & & \bulmat{2}{2} & & & \bulmat{2}{2} & \bulmat{2}{2} & & & \\
&a_{-1} & & \bulmat{1}{2} & \bulmat{1}{2} & \bulmat{1}{2} & \bulmat{1}{2} & \bulmat{2}{2} & \bulmat{2}{2} & & \bulmat{2}{2} & \bulmat{2}{2} & \bulmat{1}{2} & \bulmat{1}{2} & \bulmat{1}{2} & \bulmat{1}{2} & \\
A_{-1}&a_{-1} & \bulmat{1}{1} & \bulmat{1}{1} & \bulmat{1}{1} & \bulmat{1}{1} & \bulmat{1}{1} & \bulmat{2}{1} & \bulmat{2}{2} & & \bulmat{2}{2} & \bulmat{2}{1} & \bulmat{1}{1} & \bulmat{1}{1} & \bulmat{1}{1} & \bulmat{1}{1} & \bulmat{1}{1}\\
&a_{-2} & \bulmat{1}{1} & \bulmat{1}{1} & & \bulmat{1}{1} & \bulmat{1}{1} & \bulmat{1}{1} & \bulmat{1}{2} & \bulmat{2}{2} & \bulmat{1}{2} & \bulmat{1}{1} & \bulmat{1}{1} & \bulmat{1}{1} & & \bulmat{1}{1} & \bulmat{1}{1}\\
A_{-2}&a_{-2} & & \bulmat{1}{1} & & \bulmat{1}{1} & \bulmat{1}{1} & \bulmat{1}{1} & \bulmat{1}{2} & \bulmat{2}{2} & \bulmat{1}{2} & \bulmat{1}{1} & \bulmat{1}{1} & \bulmat{1}{1} & & \bulmat{1}{1} &\\
&a_{-3} & & & \bulmat{1}{1} & & & \bulmat{1}{1} & \bulmat{1}{2} & & \bulmat{1}{2} & \bulmat{1}{1} & & & \bulmat{1}{1} & &\\
&a_{-4} & & & & \bulmat{1}{1} & \bulmat{1}{1} & \bulmat{1}{1} & \bulmat{1}{2} & & \bulmat{1}{2} & \bulmat{1}{1} & \bulmat{1}{1} & \bulmat{1}{1} & & &\\
A_{-3}&a_{-4} & & & & & \bulmat{1}{1} & \bulmat{1}{1} & & & & \bulmat{1}{1} & \bulmat{1}{1} & & & &\\
\end{tblr}$$
\caption{Quipu - structure of $u$}
\label{ublocksq}
\end{figure}

\begin{figure}[H]
\centering
$$\begin{tblr}{
colspec={cccccc},
column{1-2} = {colsep=4pt},
hlines,
vlines,
vline{3,18} = {7-8,12-13}{2pt},
vline{4,17} = {6-9,11-14}{2pt},
vline{5,16} = {5,15}{2pt},
vline{6,15} = {4-5,10,15-16}{2pt},
vline{7,14} = {3-4,6-14,16-17}{2pt},
vline{8,13} = {5,10,15}{2pt},
vline{9,12} = {3-9,11-17}{2pt},
vline{10,11} = {10}{2pt},
hline{3,18} = {7-8,12-13}{2pt},
hline{4,17} = {6,9-11,14}{2pt},
hline{5,16} = {5-15}{2pt},
hline{6,15} = {4-16}{2pt},
hline{7,14} = {3,17}{2pt},
hline{8,13} = {3-17}{2pt},
hline{9,12} = {3,17}{2pt},
hline{10,11} = {4-16}{2pt},
cell{1}{1} = {r=2,c=2}{c},
cell{1}{7,12} = {c=2}{c},
cell{1}{4,9,14} = {c=3}{c},
cell{3,6,8,11,13,16}{1} = {r=2}{c},
}
\diagboxthree[width=5em,height=5em]{qp}{v}{\rot{sr}}& & \rot{A_3} & \rot{A_2} & & & \rot{A_1} & & \rot{A_0} & & \rot{A_0} & \rot{A_{-1}} & & \rot{A_{-2}} & & & \rot{A_{-3}}\\
& & \rot{a_4} & \rot{a_4} & \rot{a_3} & \rot{a_2} & \rot{a_2} & \rot{a_1} & \rot{a_1} & \rot{a_0} & \rot{a_{-1}} & \rot{a_{-1}} & \rot{a_{-2}} & \rot{a_{-2}} & \rot{a_{-3}} & \rot{a_{-4}} & \rot{a_{-4}}\\
a_4&A_3& & & & & \bulmat{1}{1} & \bulmat{1}{1} & & & & \bulmat{1}{1} & \bulmat{1}{1} & & & &\\
&A_2& & & & \bulmat{1}{1} & \bulmat{1}{1} & \bulmat{1}{1} & \bulmat{2}{1} & & \bulmat{2}{1} & \bulmat{1}{1} & \bulmat{1}{1} & \bulmat{1}{1} & & &\\
a_3&A_2& & & \bulmat{1}{1} & & & \bulmat{1}{1} & \bulmat{2}{1} & & \bulmat{2}{1} & \bulmat{1}{1} & & & \bulmat{1}{1} & &\\
a_2&A_2& & \bulmat{1}{1} & & \bulmat{1}{1} & \bulmat{1}{1} & \bulmat{1}{1} & \bulmat{2}{1} & \bulmat{2}{2} & \bulmat{2}{1} & \bulmat{1}{1} & \bulmat{1}{1} & \bulmat{1}{1} & & \bulmat{1}{1} &\\
&A_1 & \bulmat{1}{1} & \bulmat{1}{1} & & \bulmat{1}{1} & \bulmat{1}{1} & \bulmat{1}{1} & \bulmat{2}{1} & \bulmat{2}{2} & \bulmat{2}{1} & \bulmat{1}{1} & \bulmat{1}{1} & \bulmat{1}{1} & & \bulmat{1}{1} & \bulmat{1}{1}\\
a_1&A_1& \bulmat{1}{1} & \bulmat{1}{1} & \bulmat{1}{1} & \bulmat{1}{1} & \bulmat{1}{1} & \bulmat{1}{2} & \bulmat{2}{2} & & \bulmat{2}{2} & \bulmat{1}{2} & \bulmat{1}{1} & \bulmat{1}{1} & \bulmat{1}{1} & \bulmat{1}{1} & \bulmat{1}{1}  \\
&A_0 & & \bulmat{2}{1} & \bulmat{2}{1} & \bulmat{2}{1} & \bulmat{2}{1} & \bulmat{2}{2} & \bulmat{2}{2} & & \bulmat{2}{2} & \bulmat{2}{2} & \bulmat{2}{1} & \bulmat{2}{1} & \bulmat{2}{1} & \bulmat{2}{1} & \\
a_0&A_0& & & & \bulmat{2}{2} & \bulmat{2}{2} & & & \bulmat{2}{2} & & & \bulmat{2}{2} & \bulmat{2}{2} & & & \\
a_{-1}&A_0& & \bulmat{2}{1} & \bulmat{2}{1} & \bulmat{2}{1} & \bulmat{2}{1} & \bulmat{2}{2} & \bulmat{2}{2} & & \bulmat{2}{2} & \bulmat{2}{2} & \bulmat{2}{1} & \bulmat{2}{1} & \bulmat{2}{1} & \bulmat{2}{1} & \\
&A_{-1}& \bulmat{1}{1} & \bulmat{1}{1} & \bulmat{1}{1} & \bulmat{1}{1} & \bulmat{1}{1} & \bulmat{1}{2} & \bulmat{2}{2} & & \bulmat{2}{2} & \bulmat{1}{2} & \bulmat{1}{1} & \bulmat{1}{1} & \bulmat{1}{1} & \bulmat{1}{1} & \bulmat{1}{1}\\
a_{-2}&A_{-1}& \bulmat{1}{1} & \bulmat{1}{1} & & \bulmat{1}{1} & \bulmat{1}{1} & \bulmat{1}{1} & \bulmat{2}{1} & \bulmat{2}{2} & \bulmat{2}{1} & \bulmat{1}{1} & \bulmat{1}{1} & \bulmat{1}{1} & & \bulmat{1}{1} & \bulmat{1}{1}\\
&A_{-2} & & \bulmat{1}{1} & & \bulmat{1}{1} & \bulmat{1}{1} & \bulmat{1}{1} & \bulmat{2}{1} & \bulmat{2}{2} & \bulmat{2}{1} & \bulmat{1}{1} & \bulmat{1}{1} & \bulmat{1}{1} & & \bulmat{1}{1} &\\
a_{-3}&A_{-2} & & & \bulmat{1}{1} & & & \bulmat{1}{1} & \bulmat{2}{1} & & \bulmat{2}{1} & \bulmat{1}{1} & & & \bulmat{1}{1} & &\\
a_{-4}&A_{-2} & & & & \bulmat{1}{1} & \bulmat{1}{1} & \bulmat{1}{1} & \bulmat{2}{1} & & \bulmat{2}{1} & \bulmat{1}{1} & \bulmat{1}{1} & \bulmat{1}{1} & & &\\
&A_{-3} & & & & & \bulmat{1}{1} & \bulmat{1}{1} & & & & \bulmat{1}{1} & \bulmat{1}{1} & & & &\\
\end{tblr}$$
\caption{Quipu - structure of $v$ }
\label{vblocksq}
\end{figure}

We make the following choices for the $1\times 1$ blocks in $u$: 
$$u^{(A_{\pm 3},a_{\pm 2})}=u^{(A_{\pm 2},a_{\pm 4})}=u^{(A_{\pm 2},a_{\pm 3})}=u^{(A_{\pm 1},a_{\pm 3})}=1,\quad u^{A_{\pm3},a_{\pm 1}}=-1.$$
Similarly for the $1\times 1$ blocks in $v$, we choose:
$$v^{(a_{\pm 2},A_{\pm 3})}=v^{(a_{\pm 4},A_{\pm 2})}=v^{(a_{\pm 3},A_{\pm 2})}=v^{(a_{\pm 3},A_{\pm 1})}=1,\quad v^{a_{\pm 1},A_{\pm3}}=-1.$$
Using the bi-unitary condition we start filling in the entries for the $2\times 2$ blocks that are connected to some $1\times 1$ block. To simplify the notation we set 
$$x_1=\sqrt{\frac{\lambda(A_3)\lambda(a_2)}{\lambda(a_4)\lambda(A_1)}},\quad y_1=\sqrt{\frac{\lambda(A_3)\lambda(a_1)}{\lambda(a_4)\lambda(A_1)}},\quad x_2=\sqrt{\frac{\lambda(a_4)}{\lambda(a_2)}},\quad y_2=\sqrt{\frac{\lambda(a_4)\lambda(A_1)}{\lambda(A_2)\lambda(a_2)}},\quad y_1=\sqrt{\frac{\lambda(A_3)\lambda(a_1)}{\lambda(A_2)\lambda(a_2)}}.$$ 
We then have the following for the blocks $v^{a_{\pm 4}, A_{\pm 1}}$:
$$\begin{tblr}{
colspec={cccccc},
column{1-2} = {colsep=4pt},
column{5}={1pt,colsep=1pt},
hline{1-4,6-7} ={1-4,6-7}{solid},
vline{1,2,3,4,7,8} = {1-4,5-6}{solid},
hline{5}= {1}{1-2}{dashed},
hline{5}= {2}{1-2}{dashed},
vline{5}= {1-2}{dashed},
vline{6}= {1-2}{dashed},
hline{5} = {1}{3-4,6-7}{2pt},
hline{5} = {2}{3-4,6-7}{2pt},
hline{3,7} = {3-4,6-7}{2pt},
vline{3,5,6,8} = {3-4,5-6}{2pt},
cell{1}{1} = {r=2,c=2}{c},
cell{1}{3,6} = {c=2}{c},
cell{3,5}{1} = {r=2}{c},
}
\diagboxthree[width=5em,height=5em]{qp}{v}{\rot{sr}}& &\rot{A_1} &&& \rot{A_{-1}}&\\
& &\rot{a_2} & \rot{a_1} && \rot{a_{-1}} & \rot{a_{-2}}\\
a_4&A_3&x_1 &-y_1 && -y_1 &x_1\\
&A_2&\bullet&\bullet &&\bullet &\bullet\\
a_{-4}&A_{-2}&\bullet &\bullet &&\bullet &\bullet\\
&A_{-3}&x_1 &-y_1 && -y_1 &x_1
\end{tblr}$$
Note that using \ref{quipueigrel} we get 
$$x_1^2+y_1^2=\frac{\lambda(A_3)(\lambda(a_1)+\lambda(a_2))}{\lambda(a_4)\lambda(A_1)}=1.$$
We get a similar result for the blocks $u^{(A_{\pm 1},a_{\pm 4})}$. Since every block has to be unitary, we determine the remaining entries in the following way:

$$\begin{tblr}{
colspec={cccccc},
column{1-2} = {colsep=4pt},
column{5}={1pt,colsep=1pt},
hline{1-4,6-7} ={1-4,6-7}{solid},
vline{1,2,3,4,7,8} = {1-4,5-6}{solid},
hline{5}= {1}{1-2}{dashed},
hline{5}= {2}{1-2}{dashed},
vline{5}= {1-2}{dashed},
vline{6}= {1-2}{dashed},
hline{5} = {1}{3-4,6-7}{2pt},
hline{5} = {2}{3-4,6-7}{2pt},
hline{3,7} = {3-4,6-7}{2pt},
vline{3,5,6,8} = {3-4,5-6}{2pt},
cell{1}{1} = {r=2,c=2}{c},
cell{1}{3,6} = {c=2}{c},
cell{3,5}{1} = {r=2}{c},
}
\diagboxthree[width=5em,height=5em]{qp}{v}{\rot{sr}}& &\rot{A_1} &&& \rot{A_{-1}}&\\
& &\rot{a_2} & \rot{a_1} && \rot{a_{-1}} & \rot{a_{-2}}\\
a_4&A_3&x_1 &-y_1 && -y_1 &x_1\\
&A_2&y_1 &x_1 &&x_1 &y_1\\
a_{-4}&A_{-2}&y_1 &x_1 &&x_1 &y_1\\
&A_{-3}&x_1 &-y_1 && -y_1 &x_1
\end{tblr}$$

We make similar choices for all the blocks $v^{(a_{\pm 4},A_{\pm 1})}$ and $u^{(A_{\pm 1}, a_{\pm 4})}$. For the blocks $u^{(A_{\pm 2},a_{\pm 2})}$, the bi-unitary condition implies 
$$\begin{tblr}{
colspec={cccccc},
column{1-2} = {colsep=4pt},
column{5}={1pt,colsep=1pt},
hline{1-4,6-7} ={1-4,6-7}{solid},
vline{1,2,3,4,7,8} = {1-4,5-6}{solid},
hline{5}= {1}{1-2}{dashed},
hline{5}= {2}{1-2}{dashed},
vline{5}= {1-2}{dashed},
vline{6}= {1-2}{dashed},
hline{5} = {1}{3-4,6-7}{2pt},
hline{5} = {2}{3-4,6-7}{2pt},
hline{3,7} = {3-4,6-7}{2pt},
vline{3,5,6,8} = {3-4,5-6}{2pt},
cell{1}{1} = {r=2,c=2}{c},
cell{1}{3,6} = {c=2}{c},
cell{3,5}{1} = {r=2}{c},
}
\diagboxthree[width=5em,height=5em]{pq}{u}{\rot{rs}}& &\rot{a_2} &&& \rot{a_{-2}}&\\
& &\rot{A_2} & \rot{A_1} && \rot{A_{-1}} & \rot{A_{-2}}\\
A_2&a_4&x_2 &y_2 && y_2 &x_2\\
&a_2&\bullet &\bullet &&\bullet &\bullet\\
A_{-2}&a_{-2}&\bullet &\bullet &&\bullet &\bullet\\
&a_{-4}&x_2 &y_2 && y_2 &x_2
\end{tblr}$$

Once again we choose the remaining entries so that every block is unitary. Thus
$$u^{(A_2,a_2)}=
\begin{pNiceArray}{c|c}[first-row,first-col]
 &A_2 & A_1\\
 a_4& x_2& y_2\\
 \hline
 a_2& -y_2& x_2\\
\end{pNiceArray}.$$

The blocks $u^{(A_{\pm 2},a_{\pm 2})}$ and $v^{(a_{\pm 2},A_{\pm 2})}$ are filled in a similar manner. Using \ref{quipueigvec} and \ref{quiput^3} we have 
$$x_2^2+y_2^2=\frac{t^3-7t^2+13t-5}{t^2-4t}=1,$$ 
hence the blocks $u^{(A_{\pm 2},a_{\pm 2})}$ and $v^{(a_{\pm 2},A_{\pm 2})}$ are unitary.

We proceed to handle the blocks $u^{(A_{\pm 1}, a_{\pm 2})}$ and $v^{(a_{\pm 2},A_{\pm 1})}$, which are connected to the previous $2\times 2$ blocks. Set $x_4=\sqrt{\frac{\lambda(A_2)}{\lambda(A_1)}}x_2=\sqrt{\frac{\lambda(A_2)\lambda(a_4)}{\lambda(A_1)\lambda(a_2)}}$, since $u^{(A_{\pm 2},a_{\pm 2})}_{a_{\pm 2},A_{\pm 1}}=x_2$ and $v^{(a_{\pm 2},A_{\pm 2})}_{A_{\pm 1},a_{\pm 2}}=x_2$ then 
$$\begin{tblr}{
colspec={cccccc},
column{1-2} = {colsep=4pt},
column{5}={1pt,colsep=1pt},
hline{1-4,6-7} ={1-4,6-7}{solid},
vline{1,2,3,4,7,8} = {1-4,5-6}{solid},
hline{5}= {1}{1-2}{dashed},
hline{5}= {2}{1-2}{dashed},
vline{5}= {1-2}{dashed},
vline{6}= {1-2}{dashed},
hline{5} = {1}{3-4,6-7}{2pt},
hline{5} = {2}{3-4,6-7}{2pt},
hline{3,7} = {3-4,6-7}{2pt},
vline{3,5,6,8} = {3-4,5-6}{2pt},
cell{1}{1} = {r=2,c=2}{c},
cell{1}{3,6} = {c=2}{c},
cell{3,5}{1} = {r=2}{c},
}
\diagboxthree[width=5em,height=5em]{pq}{u}{\rot{rs}}& &\rot{a_2} &&& \rot{a_{-2}}&\\
& &\rot{A_2} & \rot{A_1} && \rot{A_{-1}} & \rot{A_{-2}}\\
A_1&a_2&x_4 &\bullet&&\bullet&x_4\\
&a_1&\bullet &\bullet &&\bullet &\bullet\\
A_{-1}&a_{-1}&\bullet &\bullet &&\bullet &\bullet\\
&a_{-2}&x_4&\bullet&&\bullet&x_4
\end{tblr}\quad\quad
\begin{tblr}{
colspec={cccccc},
column{1-2} = {colsep=4pt},
column{5}={1pt,colsep=1pt},
hline{1-4,6-7} ={1-4,6-7}{solid},
vline{1,2,3,4,7,8} = {1-4,5-6}{solid},
hline{5}= {1}{1-2}{dashed},
hline{5}= {2}{1-2}{dashed},
vline{5}= {1-2}{dashed},
vline{6}= {1-2}{dashed},
hline{5} = {1}{3-4,6-7}{2pt},
hline{5} = {2}{3-4,6-7}{2pt},
hline{3,7} = {3-4,6-7}{2pt},
vline{3,5,6,8} = {3-4,5-6}{2pt},
cell{1}{1} = {r=2,c=2}{c},
cell{1}{3,6} = {c=2}{c},
cell{3,5}{1} = {r=2}{c},
}
\diagboxthree[width=5em,height=5em]{qp}{v}{\rot{sr}}& &\rot{A_1} &&& \rot{A_{-1}}&\\
& &\rot{a_2} & \rot{a_1} && \rot{a_{-1}} & \rot{a_{-2}}\\
a_2&A_2&x_4 &\bullet&&\bullet&x_4\\
&A_1&\bullet &\bullet &&\bullet &\bullet\\
a_{-2}&A_{-1}&\bullet &\bullet &&\bullet &\bullet\\
&A_{-2}&x_4 &\bullet&&\bullet&x_4
\end{tblr}$$

Let $z_4=\sqrt{1-x_4^2}$ and choose the remaining entries in $v^{(a_{\pm 2},A_{\pm 1})}$ in the following way
$$\begin{tblr}{
colspec={cccccc},
column{1-2} = {colsep=4pt},
column{5}={1pt,colsep=1pt},
hline{1-4,6-7} ={1-4,6-7}{solid},
vline{1,2,3,4,7,8} = {1-4,5-6}{solid},
hline{5}= {1}{1-2}{dashed},
hline{5}= {2}{1-2}{dashed},
vline{5}= {1-2}{dashed},
vline{6}= {1-2}{dashed},
hline{5} = {1}{3-4,6-7}{2pt},
hline{5} = {2}{3-4,6-7}{2pt},
hline{3,7} = {3-4,6-7}{2pt},
vline{3,5,6,8} = {3-4,5-6}{2pt},
cell{1}{1} = {r=2,c=2}{c},
cell{1}{3,6} = {c=2}{c},
cell{3,5}{1} = {r=2}{c},
}
\diagboxthree[width=5em,height=5em]{qp}{v}{\rot{sr}}& &\rot{A_1} &&& \rot{A_{-1}}&\\
& &\rot{a_2} & \rot{a_1} && \rot{a_{-1}} & \rot{a_{-2}}\\
a_2&A_2&x_4 &z_4&&z_4&x_4\\
&A_1&-z_4 &x_4 &&x_4 &-z_4\\
a_{-2}&A_{-1}&-z_4 &x_4 &&x_4 &-z_4\\
&A_{-2}&x_4 &z_4&&z_4&x_4
\end{tblr}$$

Since $v^{(a_{\pm 2},A_{\pm 1})}_{A_{\pm 1},a_{\pm 2}}=-z_4$, the bi-unitary condition implies that $u^{(A_{\pm 1},a_{\pm 2})}_{a_{\pm 2},A_{\pm 1}}=-z_4$ and consequently the blocks $u^{(A_{\pm 1},a_{\pm 2})}$ are obtained in the following way:

$$\begin{tblr}{
colspec={cccccc},
column{1-2} = {colsep=4pt},
column{5}={1pt,colsep=1pt},
hline{1-4,6-7} ={1-4,6-7}{solid},
vline{1,2,3,4,7,8} = {1-4,5-6}{solid},
hline{5}= {1}{1-2}{dashed},
hline{5}= {2}{1-2}{dashed},
vline{5}= {1-2}{dashed},
vline{6}= {1-2}{dashed},
hline{5} = {1}{3-4,6-7}{2pt},
hline{5} = {2}{3-4,6-7}{2pt},
hline{3,7} = {3-4,6-7}{2pt},
vline{3,5,6,8} = {3-4,5-6}{2pt},
cell{1}{1} = {r=2,c=2}{c},
cell{1}{3,6} = {c=2}{c},
cell{3,5}{1} = {r=2}{c},
}
\diagboxthree[width=5em,height=5em]{pq}{u}{\rot{rs}}& &\rot{a_2} &&& \rot{a_{-2}}&\\
& &\rot{A_2} & \rot{A_1} && \rot{A_{-1}} & \rot{A_{-2}}\\
A_1&a_2&x_4 &-z_4&&-z_4&x_4\\
&a_1&z_4 &x_4 &&x_4 &z_4\\
A_{-1}&a_{-1}&z_4 &x_4 &&x_4 &z_4\\
&a_{-2}&x_4&-z_4&&-z_4&x_4
\end{tblr}$$

We proceed to determine the entries of the $3\times 3$ blocks $u^{(A_{\pm2},a_{\pm 1})}$ and $v{(a_{\pm 1},A_{\pm2})}$. Let
\begin{align*}
x_3&=\sqrt{\frac{\lambda(a_4)\lambda(A_1)}{\lambda(A_2)\lambda(a_1)}},& x_1&=\sqrt{\frac{\lambda(a_2)\lambda(A_3)}{\lambda(A_2)\lambda(a_1)}},& y_3&=\sqrt{\frac{\lambda(a_3)\lambda(A_1)}{\lambda(A_2)\lambda(a_1)}},&
z_3&=\sqrt{\frac{\lambda(a_2)\lambda(A_1)}{\lambda(A_2)\lambda(a_1)}}z_4.
\end{align*}
Since 
\begin{align*}
v^{(a_{\pm 4},A_{\pm 1})}_{A_{\pm 2},a_{\pm 1}}&=x_1 & v^{(a_{\pm 3},A_{\pm 1})}_{A_{\pm 2},a_{\pm 1}}&=1 & v^{(a_{\pm 
2},A_{\pm 1})}_{A_{\pm 2},a_{\pm 1}}&=z_4,\\
u^{(A_{\pm 1},a_{\pm 4})}_{a_{\pm 1},A_{\pm 2}}&=x_1 & u^{(A_{\pm 1},a_{\pm 3})}_{a_{\pm 1},A_{\pm 2}}&=1 & u^{(A_{\pm 1},a_{\pm 2})}_{a_{\pm 1},A_{\pm 2}}&=z_4,
\end{align*}
the bi-unitary condition implies that the blocks are given by

$$\begin{tblr}{
colspec={cccccc},
column{1-2} = {colsep=4pt},
column{5}={1pt,colsep=1pt},
hline{1-5,7-9} ={1-4,6-7}{solid},
vline{1,2,3,4,7,8} = {1-5,6-8}{solid},
hline{6}= {1}{1-2}{dashed},
hline{6}= {2}{1-2}{dashed},
vline{5}= {1-2}{dashed},
vline{6}= {1-2}{dashed},
hline{6} = {1}{3-4,6-7}{2pt},
hline{6} = {2}{3-4,6-7}{2pt},
hline{3,9} = {3-4,6-7}{2pt},
vline{3,5,6,8} = {3-5,6-8}{2pt},
cell{1}{1} = {r=2,c=2}{c},
cell{1}{3,6} = {c=2}{c},
cell{3,6}{1} = {r=3}{c},
}
\diagboxthree[width=5em,height=5em]{pq}{u}{\rot{rs}}& &\rot{a_1} &&& \rot{a_{-1}}&\\
& &\rot{A_1} & \rot{A_0} && \rot{A_{0}} & \rot{A_{-1}}\\
A_2&a_4&x_3 &\bulmat{1}{2}&&\bulmat{1}{2}&x_3\\
&a_3&y_3 &\bulmat{1}{2}&&\bulmat{1}{2} &y_3\\
&a_2&z_3 &\bulmat{1}{2}&&\bulmat{1}{2} &z_3\\
A_{-2}&a_{-2}&z_3 &\bulmat{1}{2}&&\bulmat{1}{2} &z_3\\
&a_{-3}&y_3&\bulmat{1}{2}&&\bulmat{1}{2}&y_3\\
&a_{-4}&x_3&\bulmat{1}{2}&&\bulmat{1}{2}&x_3
\end{tblr}\quad\quad
\begin{tblr}{
colspec={cccccc},
column{1-2} = {colsep=4pt},
column{6}={1pt,colsep=1pt},
hline{1-4,6-7} ={1-5,7-9}{solid},
vline{1,2,3,4,5,8,9,10} = {1-4,5-6}{solid},
hline{5}= {1}{1-2}{dashed},
hline{5}= {2}{1-2}{dashed},
vline{6}= {1-2}{dashed},
vline{7}= {1-2}{dashed},
hline{5} = {1}{3-5,7-9}{2pt},
hline{5} = {2}{3-5,7-9}{2pt},
hline{3,7} = {3-5,7-9}{2pt},
vline{3,6,7,10} = {3-4,5-6}{2pt},
cell{1}{1} = {r=2,c=2}{c},
cell{1}{3,7} = {c=3}{c},
cell{3,5}{1} = {r=2}{c},
}
\diagboxthree[width=5em,height=5em]{qp}{v}{\rot{sr}}& &\rot{A_2} &&&& \rot{A_{-2}}&\\
& &\rot{a_4} & \rot{a_3} & \rot{a_2} && \rot{a_{-2}} & \rot{a_{-3}} & \rot{a_{-4}}\\
a_1&A_1&x_3 &y_3&z_3&&z_3&y_3&x_3\\
&A_0&\bulmat{2}{1}&\bulmat{2}{1}&\bulmat{2}{1}&&\bulmat{2}{1}&\bulmat{2}{1}&\bulmat{2}{1}\\
a_{-1}&A_{0}&\bulmat{2}{1}&\bulmat{2}{1}&\bulmat{2}{1}&&\bulmat{2}{1}&\bulmat{2}{1}&\bulmat{2}{1}\\
&A_{-1}&x_3 &y_3&z_3&&z_3&y_3&x_3
\end{tblr}
$$

Using \ref{quipueigvec} and \ref{quiput^3} we get
$$x_3^2+y_3^2+z_3^2=\frac{-t^5+12t^4-49t^3+76t^2-33t+5}{t^3-5t^2+4t}=1.$$
Thus, the known column (or row) of every block has norm $1$.

Let $\alpha_3,\beta_3,\gamma_3\in \C^2$ be such that 
$u^{(A_2,a_1)}=\begin{pNiceArray}{c|cc}[first-row,first-col]
    &A_1&\Block{1-2}{A_0} &\\
a_4&x_3&\alpha_{3,1}&\alpha_{3,2}\\
\hline
a_3&y_3&\beta_{3,1}&\beta_{3,2}\\
\hline
a_2&z_3&\gamma_{3,1}&\gamma_{3,2}
\end{pNiceArray}$ is unitary. We fill the remaining blocks in the following way

$$\begin{tblr}{
colspec={cccccc},
column{1-2} = {colsep=4pt},
column{5}={1pt,colsep=1pt},
hline{1-5,7-9} ={1-4,6-7}{solid},
vline{1,2,3,4,7,8} = {1-5,6-8}{solid},
hline{6}= {1}{1-2}{dashed},
hline{6}= {2}{1-2}{dashed},
vline{5}= {1-2}{dashed},
vline{6}= {1-2}{dashed},
hline{6} = {1}{3-4,6-7}{2pt},
hline{6} = {2}{3-4,6-7}{2pt},
hline{3,9} = {3-4,6-7}{2pt},
vline{3,5,6,8} = {3-5,6-8}{2pt},
cell{1}{1} = {r=2,c=2}{c},
cell{1}{3,6} = {c=2}{c},
cell{3,6}{1} = {r=3}{c},
}
\diagboxthree[width=5em,height=5em]{pq}{u}{\rot{rs}}& &\rot{a_1} &&& \rot{a_{-1}}&\\
& &\rot{A_1} & \rot{A_0} && \rot{A_{0}} & \rot{A_{-1}}\\
A_2&a_4&x_3 &\alpha_3&&\tilde{\alpha_3}&x_3\\
&a_3&y_3 &\beta_3&&\tilde{\beta_3}&y_3\\
&a_2&z_3 &\gamma_3&&\tilde{\gamma_3}&z_3\\
A_{-2}&a_{-2}&z_3 &\gamma_3&&\tilde{\gamma_3}&z_3\\
&a_{-3}&y_3&\beta_3&&\tilde{\beta_3}&x_3\\
&a_{-4}&x_3&\alpha_3&&\tilde{\alpha_3}&x_3
\end{tblr}\quad\quad
\begin{tblr}{
colspec={cccccc},
column{1-2} = {colsep=4pt},
column{6}={1pt,colsep=1pt},
hline{1-4,6-7} ={1-5,7-9}{solid},
vline{1,2,3,4,5,8,9,10} = {1-4,5-6}{solid},
hline{5}= {1}{1-2}{dashed},
hline{5}= {2}{1-2}{dashed},
vline{6}= {1-2}{dashed},
vline{7}= {1-2}{dashed},
hline{5} = {1}{3-5,7-9}{2pt},
hline{5} = {2}{3-5,7-9}{2pt},
hline{3,7} = {3-5,7-9}{2pt},
vline{3,6,7,10} = {3-4,5-6}{2pt},
cell{1}{1} = {r=2,c=2}{c},
cell{1}{3,7} = {c=3}{c},
cell{3,5}{1} = {r=2}{c},
}
\diagboxthree[width=5em,height=5em]{qp}{v}{\rot{sr}}& &\rot{A_2} &&&& \rot{A_{-2}}&\\
& &\rot{a_4} & \rot{a_3} & \rot{a_2} && \rot{a_{-2}} & \rot{a_{-3}} & \rot{a_{-4}}\\
a_1&A_1&x_3 &y_3&z_3&&z_3&y_3&x_3\\
&A_0&\alpha^t_3&\beta^t_3&\gamma^t_3&&\gamma^t_3&\beta^t_3&\alpha^t_3\\
a_{-1}&A_{0}&\tilde{\alpha_3}^t&\tilde{\beta_3}^t&\tilde{\gamma_3}^t&&\tilde{\gamma_3}^t&\tilde{\beta_3}^t&\tilde{\alpha_3}^t\\
&A_{-1}&x_3 &y_3&z_3&&z_3&y_3&x_3
\end{tblr}
$$

where $\tilde{\alpha_i}=(-\alpha_{i,2},\alpha_{i,1})$. Observe that since each row and each column has norm 1 we have
\begin{equation}\label{gamma3norm}
\|\gamma_3\|^2=x_3^2+y_3^2=\frac{\lambda(a_2)\lambda(A_3)+\lambda(a_3)\lambda(A_1)}{\lambda(A_2)\lambda(a_1)}.
\end{equation} 
This allows us to determine the entries of the $2\times 2$ blocks $u^{(A_0,a_{\pm 4})}$, $u^{(A_0,a_{\pm 3})}$, $v^{(a_{\pm 3},A_0)}$ and $v^{(a_{\pm 4},A_0)}$. 
Set 
$$\alpha_4=\sqrt{\frac{\lambda(A_2)\lambda(a_1)}{\lambda(a_4)\lambda(A_0)}}\alpha_3,\quad \beta_4=\sqrt{\frac{\lambda(A_2)\lambda(a_1)}{\lambda(a_3)\lambda(A_0)}}\alpha_3.$$
Then the bi-unitary condition implies

$$
\begin{tblr}{
colspec={cccccc},
column{1-2} = {colsep=4pt},
column{5}={1pt,colsep=1pt},
hline{1-7} ={1-4,6-7}{solid},
vline{1,2,3,4,7,8} = {1-4,5-6}{solid},
vline{5}= {1-2}{dashed},
vline{6}= {1-2}{dashed},
hline{5} = {3-4,6-7}{2pt},
hline{3,7} = {3-4,6-7}{2pt},
vline{3,4,5,6,7,8} = {3-4,5-6}{2pt},
cell{1}{1} = {r=2,c=2}{c},
cell{3}{1} = {r=2}{c},
}
\diagboxthree[width=5em,height=5em]{pq}{u}{\rot{rs}}& &\rot{a_4}&\rot{a_3}& & \rot{a_{-3}}& \rot{a_{-4}}\\
& &\rot{A_2} & \rot{A_2} && \rot{A_{-2}} & \rot{A_{-2}}\\
A_0&a_1&\alpha_4 &\beta_4&&\beta_4&\alpha_4\\
&a_{-1}&\tilde{\alpha_4}&\tilde{\beta_4}&&\tilde{\beta_4}&\tilde{\alpha_4}\\
\end{tblr}\quad\quad
\begin{tblr}{
colspec={cccccc},
column{1-2} = {colsep=4pt},
column{5}={1pt,colsep=1pt},
hline{1-4,6-7} ={1-4,6-7}{solid},
vline{1,2,3,4,7,8} = {1-4,5-6}{solid},
hline{5}= {1}{1-2}{dashed},
hline{5}= {2}{1-2}{dashed},
vline{5}= {1-2}{dashed},
vline{6}= {1-2}{dashed},
hline{5} = {1}{3-4,6-7}{2pt},
hline{5} = {2}{3-4,6-7}{2pt},
hline{3,4,6,7} = {3-4,6-7}{2pt},
vline{3,5,6,8} = {3-4,5-6}{2pt},
cell{1}{1} = {r=2,c=2}{c},
cell{1}{3} = {c=2}{c},
}
\diagboxthree[width=5em,height=5em]{qp}{v}{\rot{sr}}& &\rot{A_0} &\\
& &\rot{a_1} & \rot{a_{-1}} \\
a_4&A_2&\alpha_4^t &\tilde{\alpha_4}^t\\
a_3&A_2&\beta_4^t &\tilde{\beta_4}^t\\
a_{-3}&A_{-2}&\beta_4^t&\tilde{\beta_4}^t\\
a_{-4}&A_{-2}&\alpha_4^t&\tilde{\alpha_4}^t
\end{tblr}
$$

Note that $\|\alpha_3\|^2=1-\frac{\lambda(a_2)\lambda(A_3)}{\lambda(A_2)\lambda(a_1)}$ and consequently $\|\alpha_4\|^2=\frac{\lambda(A_2)\lambda(a_1)-\lambda(A_3)\lambda(a_2)}{\lambda(A_0)\lambda(a_4)}$. Using \ref{quipueigvec} and \ref{quiput^3} we have 
$$\|\alpha_4\|^2=\frac{-t^5+12 t^4-48t^3+73t^2-37 t+5}{2 t^4-14t^3+26 t^2-10 t}=1,$$
similarly we obtain $\|\beta_4\|^2=1$. Hence, all the blocks above are unitary. 

For the remaining four $3\times 3$ blocks of $u$ and $v$, we set $x_5=\sqrt{\frac{\lambda(a_2)}{\lambda(a_1)}}x_4=\sqrt{\frac{\lambda(A_2)\lambda(a_4)}{\lambda(A_1)\lambda(a_1)}}$. Then, by the bi-unitary condition, we have
$$\begin{tblr}{
colspec={cccccc},
column{1-2} = {colsep=4pt},
column{5}={1pt,colsep=1pt},
hline{1-4,6-7} ={1-4,6-7}{solid},
vline{1,2,3,4,7,8} = {1-4,5-6}{solid},
hline{5}= {1}{1-2}{dashed},
hline{5}= {2}{1-2}{dashed},
vline{5}= {1-2}{dashed},
vline{6}= {1-2}{dashed},
hline{5} = {1}{3-4,6-7}{2pt},
hline{5} = {2}{3-4,6-7}{2pt},
hline{3,7} = {3-4,6-7}{2pt},
vline{3,5,6,8} = {3-4,5-6}{2pt},
cell{1}{1} = {r=2,c=2}{c},
cell{1}{3,6} = {c=2}{c},
cell{3,5}{1} = {r=2}{c},
}
\diagboxthree[width=5em,height=5em]{pq}{u}{\rot{rs}}& &\rot{a_1} &&& \rot{a_{-1}}&\\
& &\rot{A_1} & \rot{A_0} && \rot{A_{0}} & \rot{A_{-1}}\\
A_1&a_2&x_5 &\bulmat{1}{2}&&\bulmat{1}{2}&x_5\\
&a_1&\bulmat{2}{1}&\bulmat{2}{2}&&\bulmat{2}{2} &\bulmat{2}{1}\\
A_{-1}&a_{-1}&\bulmat{2}{1}&\bulmat{2}{2}&&\bulmat{2}{2}&\bulmat{2}{1}\\
&a_{-2}&x_5&\bulmat{1}{2}&&\bulmat{1}{2}&x_5
\end{tblr}\quad\quad
\begin{tblr}{
colspec={cccccc},
column{1-2} = {colsep=4pt},
column{5}={1pt,colsep=1pt},
hline{1-4,6-7} ={1-4,6-7}{solid},
vline{1,2,3,4,7,8} = {1-4,5-6}{solid},
hline{5}= {1}{1-2}{dashed},
hline{5}= {2}{1-2}{dashed},
vline{5}= {1-2}{dashed},
vline{6}= {1-2}{dashed},
hline{5} = {1}{3-4,6-7}{2pt},
hline{5} = {2}{3-4,6-7}{2pt},
hline{3,7} = {3-4,6-7}{2pt},
vline{3,5,6,8} = {3-4,5-6}{2pt},
cell{1}{1} = {r=2,c=2}{c},
cell{1}{3,6} = {c=2}{c},
cell{3,5}{1} = {r=2}{c},
}
\diagboxthree[width=5em,height=5em]{qp}{v}{\rot{sr}}& &\rot{A_1} &&& \rot{A_{-1}}&\\
& &\rot{a_2} & \rot{a_1} && \rot{a_{-1}} & \rot{a_{-2}}\\
a_1&A_1&x_5&\bulmat{1}{2}&&\bulmat{1}{2}&x_5\\
&A_0&\bulmat{2}{1}&\bulmat{2}{2}&&\bulmat{2}{2} &\bulmat{2}{1}\\
a_{-1}&A_0&\bulmat{2}{1}&\bulmat{2}{2}&&\bulmat{2}{2} &\bulmat{2}{1}\\
&A_{-1}&x_5&\bulmat{1}{2}&&\bulmat{1}{2}&x_5
\end{tblr}$$
and to determine the remaining entries we make the following choices
$$\begin{tblr}{
colspec={cccccccc},
column{1-2} = {colsep=4pt},
column{6}={1pt,colsep=1pt},
hline{1-4,8-9} ={1-5,7-9}{solid},
vline{1,2,3,4,9,10} = {1-5,6-8}{solid},
hline{6}= {1}{1-2}{dashed},
hline{6}= {2}{1-2}{dashed},
vline{6}= {1-2}{dashed},
vline{7}= {1-2}{dashed},
hline{6} = {1}{3-5,7-9}{2pt},
hline{6} = {2}{3-5,7-9}{2pt},
hline{3,9} = {3-5,7-9}{2pt},
vline{3,6,7,10} = {3-5,6-8}{2pt},
cell{1}{1} = {r=2,c=2}{c},
cell{1}{3,7} = {c=3}{c},
cell{2}{4,7} = {c=2}{c},
cell{4,6}{2} = {r=2}{c},
cell{3,6}{1} = {r=3}{c},
}
\diagboxthree[width=5em,height=5em]{pq}{u}{\rot{rs}}& &\rot{a_1} &&&& \rot{a_{-1}}&&\\
& &\rot{A_1} & \rot{A_0} &&& \rot{A_{0}} && \rot{A_{-1}}\\
A_1&a_2&x_5 &y_5&0&&0&y_5&x_5\\
&a_1&y_5&-x_5&0&&0&-x_5&y_5\\
&&0&0&1&&1&0&0\\
A_{-1}&a_{-1}&0&0&1&&1&0&0\\
&&y_5&-x_5&0&&0&-x_5&y_5\\
&a_{-2}&x_5&y_5&0&&0&y_5&x_5
\end{tblr}\quad\quad
\begin{tblr}{
colspec={cccccccc},
column{1-2} = {colsep=4pt},
column{6}={1pt,colsep=1pt},
hline{1-4,8-9} ={1-5,7-9}{solid},
vline{1,2,3,4,9,10} = {1-5,6-8}{solid},
hline{6}= {1}{1-2}{dashed},
hline{6}= {2}{1-2}{dashed},
vline{6}= {1-2}{dashed},
vline{7}= {1-2}{dashed},
hline{6} = {1}{3-5,7-9}{2pt},
hline{6} = {2}{3-5,7-9}{2pt},
hline{3,9} = {3-5,7-9}{2pt},
vline{3,6,7,10} = {3-5,6-8}{2pt},
cell{1}{1} = {r=2,c=2}{c},
cell{1}{3,7} = {c=3}{c},
cell{2}{4,7} = {c=2}{c},
cell{4,6}{2} = {r=2}{c},
cell{3,6}{1} = {r=3}{c},
}
\diagboxthree[width=5em,height=5em]{qp}{v}{\rot{sr}}& &\rot{A_1} &&&& \rot{A_{-1}}&&\\
& &\rot{a_2} & \rot{a_1} &&& \rot{a_{-1}} && \rot{a_{-2}}\\
a_1&A_1&x_5 &y_5&0&&0&y_5&x_5\\
&A_0&y_5&-x_5&0&&0&-x_5&y_5\\
&&0&0&1&&1&0&0\\
a_{-1}&A_0&0&0&1&&1&0&0\\
&&y_5&-x_5&0&&0&-x_5&y_5\\
&A_{-1}&x_5&y_5&0&&0&y_5&x_5
\end{tblr}$$

where $y_5=\sqrt{1-x_5^2}=\sqrt{\frac{\lambda(A_1)\lambda(a_1)-\lambda(A_2)\lambda(a_4)}{\lambda(A_1)\lambda(a_1)}}$. All blocks are unitary with these choices. We can use these to obtain the entries of the $4\times 4$ blocks $u^{(A_0,a_{\pm 2})}$ and $v^{(a_{\pm 2},A_0)}$. Set $$y_4=\sqrt{\frac{\lambda(a_1)\lambda(A_1)}{\lambda(A_0)\lambda(a_2)}}y_5=\sqrt{\frac{\lambda(A_1)\lambda(a_1)-\lambda(A_2)\lambda(a_4)}{\lambda(A_0)\lambda(a_2)}},\quad \gamma_4=\sqrt{\frac{\lambda(A_2)\lambda(a_1)}{\lambda(a_2)\lambda(A_0)}}\gamma_3,$$ 
then the bi-unitary condition implies
$$\begin{tblr}{
colspec={ccccccccc},
column{1-2} = {colsep=4pt},
column{6}={1pt,colsep=1pt},
hline{1-6} ={1-5,7-9}{solid},
vline{1,2,3,4,9,10} = {1-5,6-8}{solid},
vline{6}= {1-2}{dashed},
vline{7}= {1-2}{dashed},
hline{3,6} = {3-5,7-9}{2pt},
vline{3,6,7,10} = {3-5,6-8}{2pt},
cell{1}{1} = {r=2,c=2}{c},
cell{3}{1} = {r=3}{c},
cell{1}{3,7} = {c=3}{c},
cell{2}{4,7} = {c=2}{c},
}
\diagboxthree[width=5em,height=5em]{pq}{u}{\rot{rs}}& &\rot{a_2} &&&& \rot{a_{-2}}&\\
& &\rot{A_2} & \rot{A_1} &&& \rot{A_{-1}} && \rot{A_{-2}}\\
A_0&a_1&\gamma_{4}&y_4 &0&&y_4&0&\gamma_4\\
&a_0&\bulmat{2}{2}&\bulmat{2}{1}&\bulmat{2}{1}&&\bulmat{2}{1}&\bulmat{2}{1}&\bulmat{2}{2}\\
&a_{-1}&\tilde{\gamma_{4}}&0&y_4&&0&y_4&\tilde{\gamma_4}\\
\end{tblr}\quad\quad
\begin{tblr}{
colspec={ccccc},
column{1-2} = {colsep=4pt},
hline{1-4,8-9} ={1-5,7-9}{solid},
vline{1,2,3,4,5,6} = {1-5,6-8}{solid},
hline{6}= {1}{1-2}{dashed},
hline{6}= {2}{1-2}{dashed},
hline{6} = {1}{3-5,7-9}{2pt},
hline{6} = {2}{3-5,7-9}{2pt},
hline{3,9} = {3-5,7-9}{2pt},
vline{3,6,7,10} = {3-5,6-8}{2pt},
cell{1}{1} = {r=2,c=2}{c},
cell{1}{3} = {c=3}{c},
cell{4,6}{2} = {r=2}{c},
cell{3,6}{1} = {r=3}{c},
}
\diagboxthree[width=5em,height=5em]{qp}{v}{\rot{sr}}& &\rot{A_0}& & \\
& &\rot{a_1} & \rot{a_{0}} & \rot{a_{-1}}\\
a_2&A_2&\gamma_4^t&\bulmat{2}{2}&\tilde{\gamma_4}^t\\
&A_1&y_4&\bulmat{1}{2}&0\\
&&0&\bulmat{1}{2}&y_4\\
a_{-2}&A_1&y_4&\bulmat{1}{2}&0\\
& &0&\bulmat{1}{2}&y_4\\
&A_{-2}&\gamma_4^t&\bulmat{2}{2}&\tilde{\gamma_4}^t
\end{tblr}$$

Note that by \ref{gamma3norm} we have 
\begin{equation}\label{gamma4norm}
\|\gamma_4\|^2= \frac{\lambda(a_2)\lambda(A_3)+\lambda(a_3)\lambda(A_1)}{\lambda(a_2)\lambda(A_0)} . 
\end{equation}
and using \ref{quipueigvec} and \ref{quiput^3} we obtain
$$y_4^2+\|\gamma_4\|^2=\frac{-t^4+9t^3-25t^2+24t-5}{2t}=1.$$
Now, let $X,Y\in M_2(\C)$ be two matrices such that 
$$u^{(A_0,a_2)}=
\begin{pNiceArray}{cc|cc}[first-row,first-col]
&\Block{1-2}{A_2} &&\Block{1-2}{A_1} &\\
a_1&\gamma_{4,1}&\gamma_{4,2}& y_4&0\\
\hline
\Block{2-1}{a_0}&X_{1,1}&X_{1,2}&Y_{1,1}&Y_{1,2}\\
&X_{2,1}&X_{2,2}&Y_{2,1}&Y_{2,2}\\
\hline
a_{-1}&-\gamma_{4,2}&\gamma_{4,1}&0&y_4\\
\end{pNiceArray}
$$

is unitary. Note that such matrices exist as the first and last row of $u^{(A_0,a_2)}$ are orthonormal. Completing these two rows to an orthonormal basis produces such $X$ and $Y$. The entries for $u^{(A_0,a_{\pm 2})}$ and $v^{(a_{\pm 2},A_0)}$ will be 
$$\begin{tblr}{
colspec={cccccccccc},
column{1-2} = {colsep=4pt},
column{6}={1pt,colsep=1pt},
hline{1-6} ={1-5,7-9}{solid},
vline{1,2,3,4,9,10} = {1-6}{solid},
vline{6}= {1-2}{dashed},
vline{7}= {1-2}{dashed},
hline{3,6} = {3-5,7-9}{2pt},
vline{3,6,7,10} = {3-5,6-8}{2pt},
cell{1}{1} = {r=2,c=2}{c},
cell{4}{4,7} = {c=2}{c},
cell{3}{1} = {r=3}{c},
cell{1}{3,7} = {c=3}{c},
cell{2}{4,7} = {c=2}{c},
}
\diagboxthree[width=5em,height=5em]{pq}{u}{\rot{rs}}& &\rot{a_2} &&&& \rot{a_{-2}}\\
& &\rot{A_2} & \rot{A_1} &&& \rot{A_{-1}} && \rot{A_{-2}}&\\
A_0&a_1&\gamma_{4}&y_4 &0&&y_4&0&\gamma_4\\
&a_0&X&Y&&&Y&&X\\
&a_{-1}&\tilde{\gamma_4}&0&y_4&&0&y_4&\tilde{\gamma_4}
\end{tblr}\quad\quad
\begin{tblr}{
colspec={ccccc},
column{1-2} = {colsep=4pt},
hline{1-4,8-9} ={1-5,7-9}{solid},
vline{1,2,3,4,5,6} = {1-5,6-8}{solid},
hline{6}= {1}{1-2}{dashed},
hline{6}= {2}{1-2}{dashed},
hline{6} = {1}{3-5,7-9}{2pt},
hline{6} = {2}{3-5,7-9}{2pt},
hline{3,9} = {3-5,7-9}{2pt},
vline{3,6,7,10} = {3-5,6-8}{2pt},
cell{1}{1} = {r=2,c=2}{c},
cell{1}{3} = {c=3}{c},
cell{4,6}{2} = {r=2}{c},
cell{3,6}{1} = {r=3}{c},
cell{4,6}{4} = {r=2}{c},
}
\diagboxthree[width=5em,height=5em]{qp}{v}{\rot{sr}}& &\rot{A_0}& & \\
& &\rot{a_1} & \rot{a_{0}} & \rot{a_{-1}}\\
a_2&A_2&\gamma_4^t&X^t&\tilde{\gamma_4}^t\\
&A_1&y_4&Y^t&0\\
&&0&&y_4\\
a_{-2}&A_{-1}&y_4&Y^t&0\\
& &0&&y_4\\
&A_{-2}&\gamma_4^t&X^t&\tilde{\gamma_4}^t
\end{tblr}$$

Since $u^{(A_0,a_2)}$ is unitary, its columns are orthonormal and therefore
\begin{align*}
    0&=X_{1,1}X_{1,2}+X_{2,1}X_{2,2}+\gamma_{4,1}\gamma_{4,2}-\gamma_{4,1}\gamma_{4,2}=X_{1,1}X_{1,2}+X_{2,1}X_{2,2}\\
    0&=Y_{1,1}Y_{1,2}+Y_{2,1}Y_{2,2}.
\end{align*}
Hence, the columns of $X$ and $Y$ are orthogonal. Moreover, since every column has norm 1, we have
\begin{align*}
    1&=\|\gamma_4\|^2+X_{1,1}^2+X_{2,1}^2=\|\gamma_4\|^2+X_{1,2}^2+X_{2,2}^2
\end{align*}
and hence
\begin{align}
X_{1,1}^2+X_{2,1}^2&=X_{1,2}^2+X_{2,2}^2=1-\|\gamma_4\|^2=\frac{\lambda(a_2)\lambda(A_0)-\lambda(a_2)\lambda(A_3)-\lambda(a_3)\lambda(A_1)}{\lambda(a_2)\lambda(A_0)}. \label{quipuXcolnorm}
\end{align}
Similarly, we get $1-y_4^2=Y_{1,1}^2+Y_{2,1}^2=Y_{1,2}^2+Y_{2,2}^2$. Set $X'=\sqrt{\frac{\lambda(A_0)\lambda(a_2)}{\lambda(a_0)\lambda(A_2)}}X$ and $Y'=\sqrt{\frac{\lambda(A_0)\lambda(a_2)}{\lambda(a_0)\lambda(A_1)}}Y$. The last four $2\times 2$ blocks for $u$ and $v$ will then be

$$\begin{tblr}{
colspec={ccc},
column{1-2} = {colsep=4pt},
hline{1-4,6,7} ={1-3}{solid},
vline{1,2,3,4,5,6} = {1-5,6-8}{solid},
hline{5}= {1}{1-2}{dashed},
hline{5}= {2}{1-2}{dashed},
hline{5} = {1}{3-5,7-9}{2pt},
hline{5} = {2}{3-5,7-9}{2pt},
hline{3,4,6,7} = {3}{2pt},
vline{3,4} = {3-4,5-6}{2pt},
cell{1}{1} = {r=2,c=2}{c}
}
\diagboxthree[width=5em,height=5em]{pq}{u}{\rot{rs}}& &\rot{a_0}\\
& &\rot{A_0}\\
A_2&a_2&X'\\
A_1&a_2&Y'\\
A_{-1}&a_{-2}&Y'\\
A_{-2}&a_{-2}&X'\\
\end{tblr}\quad\quad
\begin{tblr}{
colspec={ccccccccc},
column{1-2} = {colsep=4pt},
column{5}={1pt,colsep=1pt},
hline{1-4} ={1-4,6-7}{solid},
vline{1,2,3,4,7,8} = {1-3}{solid},
vline{5}= {1-2}{dashed},
vline{6}= {1-2}{dashed},
hline{3,4} = {3-4,6-7}{2pt},
vline{3,4,5,6,7,8} = {3}{2pt},
cell{1}{1} = {r=2,c=2}{c}
}
\diagboxthree[width=5em,height=5em]{qp}{v}{\rot{sr}}& &\rot{A_2}&\rot{A_1}&&\rot{A_1}&\rot{A_2}\\
& &\rot{a_2}&\rot{a_2}&&\rot{a_{-2}}&\rot{a_{-2}}\\
a_0&A_0&(X')^t&(Y')^t&&(Y')^t&(X')^t
\end{tblr}$$

Using \ref{quipueigvec}, \ref{quiput^3} and \ref{quipuXcolnorm} we get
\begin{align*}
    (X'_{1,1})^2+(X'_{2,1})^2=(X'_{1,2})^2+(X'_{2,2})^2&=\frac{\lambda(a_2)\lambda(A_0)-\lambda(a_2)\lambda(A_3)-\lambda(a_3)\lambda(A_1)}{\lambda(a_0)\lambda(A_2)}\\
    &=\frac{t^4-9 t^3+23 t^2-14 t+5}{-2t^2+8t}=1\\
    (Y'_{1,1})^2+(Y'_{2,1})^2=(Y'_{1,2})^2+(Y'_{2,2})^2&=\frac{\lambda(A_0)\lambda(a_2)-\lambda(A_1)\lambda(a_1)+\lambda(A_2)\lambda(a_4)}{\lambda(a_0)\lambda(A_1)}\\
    &=\frac{t^4-9t^3+25t^2-22t+7}{2}=1
\end{align*}

and since the columns of $X'$ and $Y'$ are orthogonal to each other we conclude that both $X'$ and $Y'$ are unitary. We now have only two $4\times 4$ blocks remaining for $u$ and $v$: $u^{(A_0,a_{\pm 1})}$ and $v^{(a_{\pm 1},A_0)}$. Using the entries of $u^{(A_{\pm 1},a_{\pm 1})}_{a_{\pm 1},A_0}$ and $v^{(a_{\pm 1},A_{\pm 1})}_{A_0,a_{\pm1 }}$, the bi-unitary condition implies

$$\begin{tblr}{
colspec={cccccccc},column{1-2} = {colsep=4pt},column{6}={1pt,colsep=1pt},hline{1-3,5,7} ={1-5,7-9}{solid},vline{1,2,3,5,7,8,10} = {1-6}{solid},vline{6}= {1-2}{dashed},vline{7}= {1-2}{dashed},hline{3,7} = {3-5,7-9}{2pt},vline{3,6,7,10} = {3-4,5-6}{2pt},cell{1}{1} = {r=2,c=2}{c},cell{1}{3,7} = {c=3}{c},cell{2}{3,8} = {c=2}{c},cell{3,5}{2,5,7} = {r=2}{c},cell{3}{1} = {r=4}{c}}
\diagboxthree[width=5em,height=5em]{pq}{u}{\rot{rs}}& &\rot{a_1} &&&& \rot{a_{-1}}&\\
& &\rot{A_1} && \rot{A_0} && \rot{A_{0}} & \rot{A_{-1}}\\
A_0&a_1&-x_6 &0&\bulmat{2}{2}&&\bulmat{2}{2}&0&-x_6\\
&&0&y_6&&&&y_6&0\\
&a_{-1}&0&-x_6&\bulmat{2}{2}&&\bulmat{2}{2}&-x_6&0\\
&&y_6&0&&&&0&y_6\\
\end{tblr}\quad\quad\begin{tblr}{
colspec={ccccc},column{1-2} = {colsep=4pt},hline{1-3,5,7,9} ={1-6}{solid},vline{1,2,3,5,7} = {1-5,6-8}{solid},hline{6}= {1}{1-2}{dashed},hline{6}= {2}{1-2}{dashed},hline{6} = {1}{3-6}{2pt},hline{6} = {2}{3-6}{2pt},hline{3,9} = {3-6}{2pt},vline{3,7} = {3-5,6-8}{2pt},cell{1}{1} = {r=2,c=2}{c},cell{1}{3} = {c=4}{c},cell{2}{3,5} = {c=2}{c},cell{3,6}{1} = {r=3}{c},cell{3,7}{2} = {r=2}{c},cell{5,6}{3,5} = {c=2}{c}}
\diagboxthree[width=5em,height=5em]{qp}{v}{\rot{sr}}& &\rot{A_0} &&&\\
& &\rot{a_1}&&\rot{a_{-1}}& \\
a_1&A_1&-x_6&0&0&y_6\\
&&0&y_6&-x_6&0\\
&A_0&\bulmat{2}{2}&&\bulmat{2}{2}&\\
a_{-1}&A_0&\bulmat{2}{2}&&\bulmat{2}{2}&\\
&A_{-1}&0&y_6&-x_6&0\\
&&-x_6&0&0&y_6
\end{tblr}$$

where $x_6=\sqrt{\frac{\lambda(A_1)}{\lambda(A_0)}}x_5=\sqrt{\frac{\lambda(A_2)\lambda(a_4)}{\lambda(A_0)\lambda(a_1)}}$ and $y_6=\sqrt{\frac{\lambda(A_1)}{\lambda(A_0)}}$. Using \ref{quipueigvec} and \ref{quiput^3} we get
$$x_6^2+y_6^2=\frac{t^4-11t^3+41t^2-56t+17}{2}=1$$
Thus, the first two columns of $u^{(A_0, a_{1})}$ have norm $1$, and moreover, they are orthonormal. Hence there exists $Z,W\in M_2(\C)$ such that 
$$u^{(A_0,a_1)}=
\begin{pNiceArray}{cc|cc}[first-row,first-col]
&\Block{1-2}{A_1} &&\Block{1-2}{A_0} &\\
\Block{2-1}{a_1}&-x_6&0&\Block{2-2}{Z}&\\
&0&y_6& &\\
\hline
\Block{2-1}{a_{-1}}&0&-x_6&\Block{2-2}{W}&\\
&y_6&0& &\\
\end{pNiceArray}$$

is unitary. Finally, we can finish the remaining entries of $u^{A_0,a_{\pm 1}}$ and $v^{(a_{\pm1,A_0})}$ in the following way:

$$\begin{tblr}{
colspec={cccccccc},column{1-2} = {colsep=4pt},column{6}={1pt,colsep=1pt},hline{1-3,5,7} ={1-5,7-9}{solid},vline{1,2,3,5,7,8,10} = {1-6}{solid},vline{6}= {1-2}{dashed},vline{7}= {1-2}{dashed},hline{3,7} = {3-5,7-9}{2pt},vline{3,6,7,10} = {3-4,5-6}{2pt},cell{1}{1} = {r=2,c=2}{c},cell{1}{3,7} = {c=3}{c},cell{2}{3,8} = {c=2}{c},cell{3,5}{2,5,7} = {r=2}{c},cell{3}{1} = {r=4}{c}}
\diagboxthree[width=5em,height=5em]{pq}{u}{\rot{rs}}& &\rot{a_1} &&&& \rot{a_{-1}}&\\
& &\rot{A_1} && \rot{A_0} && \rot{A_{0}} & \rot{A_{-1}}\\
A_0&a_1&-x_6 &0&Z&&W&0&-x_6\\
&&0&y_6&&&&y_6&0\\
&a_{-1}&0&-x_6&W&&Z&-x_6&0\\
&&y_6&0&&&&0&y_6\\
\end{tblr}\quad\quad\begin{tblr}{
colspec={ccccc},column{1-2} = {colsep=4pt},hline{1-3,5,7,9} ={1-6}{solid},vline{1,2,3,5,7} = {1-5,6-8}{solid},hline{6}= {1}{1-2}{dashed},hline{6}= {2}{1-2}{dashed},hline{6} = {1}{3-6}{2pt},hline{6} = {2}{3-6}{2pt},hline{3,9} = {3-6}{2pt},vline{3,7} = {3-5,6-8}{2pt},cell{1}{1} = {r=2,c=2}{c},cell{1}{3} = {c=4}{c},cell{2}{3,5} = {c=2}{c},cell{3,6}{1} = {r=3}{c},cell{3,7}{2} = {r=2}{c},cell{5,6}{3,5} = {c=2}{c}}
\diagboxthree[width=5em,height=5em]{qp}{v}{\rot{sr}}& &\rot{A_0} &&&\\
& &\rot{a_1}&&\rot{a_{-1}}& \\
a_1&A_1&-x_6&0&0&y_6\\
&&0&y_6&-x_6&0\\
&A_0&Z^t&&W^t&\\
a_{-1}&A_0&W^t&&Z^t&\\
&A_{-1}&0&y_6&-x_6&0\\
&&-x_6&0&0&y_6
\end{tblr}$$

Observe that this makes all blocks unitary and is consistent with the bi-unitary condition. We conclude that there is a commuting square with vertical inclusions given by $G$ and horizontal inclusions given by $H$ and $K$. 

\begin{theo}
   There is an irreducible, hyperfinite $A_\infty$-subfactor with index $\approx 4.3720$. 
\end{theo}
\begin{proof}
    Note that $\Gamma$ is not one of the graphs into whose graph planar algebra, the Extended-Haagerup subfactor planar algebra embeds. The embedding theorem \ref{embeddingT} implies that the subfactor obtained with the commuting square constructed above is not the Extended-Haagerup subfactor and consequently has infinite depth. Since the $1$-norm of the first row of $G$ is $1$, Wenzl's criterion implies that this subfactor is irreducible. Therefore we have constructed an irreducible hyperfinite subfactor with index $\approx 4.37720$ and infinite depth. By classification it must have trivial standard invariant.    
\end{proof}

    \section{A 1-parameter family of commuting squares for 4-stars}\label{1param}

\begin{defi}[\cite{schou2013commuting}]
We say that a graph $G$ is an $m$-star, if $G$ is connected, and has a ``central'' vertex $p$, of valency $m$, and $m$ rays of the form $\begin{tikzpicture}[main/.style={circle,fill=black,inner sep=0pt,minimum size=1.5mm},
        sub/.style={circle,fill=white,draw=black,inner sep=0pt,minimum size=1.5mm},baseline=-2]
	\node[main,label=below:{$p$}] at (0,0) (A0) {};
	\node[main] at (1,0) (A1) {};
	\node[main] at (2,0) (A2) {};
	\node at (3,0) (A3) {$\dots$};
	\node[main] at (4,0) (A4) {};
	\node[main] at (5,0) (A5) {};
	
	\draw (A0) -- (A1) -- (A2) -- (A3) -- (A4) -- (A5) ;
\end{tikzpicture}$ with $k_i$ vertices (not counting $p$), $i=1,\dots,m$. We will denote an $m$-star by $S(k_1,k_2,\dots,k_m)$, $k_1\leq k_2\leq \dots\leq k_m$.
\end{defi}

Let $G$ be a 4-star of the form $G=S(i,i,j,j)$. We show the existence of a 1-parameter family of non-equivalent bi-unitary connections for inclusions of the form 
\begin{equation}
\begin{array}{ccc}
A_{1,0}&\stackrel{G^t}{\sst} &A_{1,1}\\
\rotatebox[origin=c]{90}{$\sst$}_G& &\rotatebox[origin=c]{90}{$\sst$}_{G^t}\\
A_{0,0}&\stackrel{G}{\sst} &A_{0,1}
\end{array}
\end{equation}
We will use the following labeling of the vertices:
\begin{center}
    \begin{tikzpicture}[main/.style={circle,fill=black,inner sep=0pt,minimum size=1.5mm},sub/.style={circle,fill=white,draw=black,inner sep=0pt,minimum size=1.5mm},end/.style={circle,diagonal fill={black}{white},draw=black,inner sep=0pt,minimum size=1.5mm},baseline=0,scale=0.6]
			\node[main,label=above:{$A$}] at (0,0) (A) {};
			\node[sub,label=above:{$a_1$}] at (-1,0.5) (a1) {};
			\node[sub,label=below:{$b_1$}] at (-1,-0.5) (b1) {};
			\node  at (-2,0.5) (ad) {$\cdots$};
			\node  at (-2,-0.5) (bd) {$\cdots$};
            \node[end,label=above:{$a_i$}]  at (-3,0.5) (ai) {};
			\node[end,label=below:{$b_j$}]  at (-3,-0.5) (bi) {};
			\node[sub,label=above:{$c_1$}]  at (1,0.5) (c1) {};
			\node[sub,label=below:{$d_1$}]  at (1,-0.5) (d1) {};
            \node at (2,0.5) (cd) {$\cdots$};
			\node at (2,-0.5) (dd) {$\cdots$};
            \node[end,label=above:{$c_j$}] at (3,0.5) (ci) {};
			\node[end,label=below:{$d_i$}] at (3,-0.5) (di) {};
			\draw (ai) -- (ad) -- (a1) -- (A) -- (b1) -- (bd) -- (bi);
			\draw (ci) -- (cd) -- (c1) -- (A) -- (d1) -- (dd) -- (di);
			\end{tikzpicture}
\end{center}
where black vertices denote minimal central projections in $\ZZ_{0,0}$, $\ZZ_{1,1}$, and white vertices denote minimal central projections in $\ZZ_{1,0}$, $\ZZ_{0,1}$. The size of the blocks of $u$ is given by the entries of $GG^t$, which are equal to the number of paths of length 2 starting and ending on a black vertex. Similarly, the size of the summands of $v$ is given by the number of paths of length 2 starting and ending on a white vertex. Consequently, we have the following summands for $u$ and $v$:
\begin{center}
    \begin{tblr}{
    hline{2}={2pt},
    vline{2}={1pt},
    }
        Size & Summands\\
        $1\times 1$&$u^{(A,x_{2})},u^{(x_{2},A)},u^{(x_{2k},x_{2(k\pm 1)})},v^{(x_{2k-1},x_{2(k\pm 1)-1})},v^{x_1,y_1}$\\
        $2\times 2$&$u^{(x_{2k},x_{2k})},v^{(x_{2k-1},x_{2k-1})}$\\
        $4\times 4$&$u^{(A,A)}$
    \end{tblr}
\end{center}
where $x\neq y\in\{a,b,c,d\}$, all other summands are $0$. In particular, the only $4\times 4$ summand has the following block structure
$$u^{(A,A)}=\begin{pNiceArray}{c|c|c|c}[first-col,first-row]
    &a_1& b_1& c_1 & d_1\\
    a_1&\ast&\ast&\ast&\ast\\ \hline
    b_1&\ast&\ast&\ast&\ast\\ \hline
    c_1&\ast&\ast&\ast&\ast\\ \hline
    d_1&\ast&\ast&\ast&\ast
\end{pNiceArray}.$$

Note that $v^{(x_1,y_1)}\in \mathbb{T}$, for $x\neq y\in \{a,b,c,d\}$. Using the bi-unitary condition, we can find all the entries off the diagonal of $u^{(A,A)}$. In particular, $u^{(A,A)}$ is unitarily equivalent to a $4\times 4$ unitary of the form
{$$\begin{pNiceMatrix}
     \alpha_1 & \alpha_2  & \alpha_2        & \alpha_3 \\
     \alpha_2 & \xi z_1 & \beta z_2       & \alpha_2 z_3\\
     \alpha_2 & \beta x_1 & \xi x_2       & \alpha_2 x_3\\
     \alpha_3 & \alpha_2 y_1 & \alpha_2 y_2 & \alpha_1 y_3\\
\end{pNiceMatrix}
$$}
where $x_i,y_i,z_i\in \mathbb{T}$ and $\alpha_2,\alpha_3,\beta$ are fixed positive real numbers of the form $w(A,x,y,A)^{-1}$ where $x\neq y\in \{a,b,c,d\}$ (see \ref{biunitfactor}). Since every row of $u$  has to have norm $1$ we have
\begin{equation}\label{ucoefeq}
    \begin{cases}
        \alpha_1^2+2\alpha_2^2+\alpha_3^2&=1, \\
    2\alpha_2^2+\beta^2+\xi^2 &=1.
    \end{cases}
\end{equation}
Schou showed in \cite{schou2013commuting} that one can obtain such $4\times 4$ unitaries for any $(i,j)$. We will prove that one can actually obtain a 1-parameter family of non-equivalent unitaries for every fixed $(i,j)$.

We first note that once we have the first three rows of $u$, these will determine three orthonormal vectors in $\C^4$. By picking an element in the orthogonal complement such that its first entry is $\alpha_3$ we will obtain the values of the $y_i$'s. Thus, we only need to determine $z_i$'s and $x_i$'s such that the first three rows of $u$ are orthonormal. 

We choose $x_3=-z_3$. The orthogonality of the second and third row of $u$ imply that $(z_1,z_2)$ and $(x_1,x_2)$ are orthogonal vectors in $\C^2$. Consequently, we have $x_1=-w \overline{z}_2$ and $x_2=w \overline{z}_1$ where $w\in \mathbb{T}$. Hence, we have the following for $u$:
{$$u=\begin{pNiceMatrix}
     \alpha_1 & \alpha_2        & \alpha_2        & \alpha_3 \\
     \alpha_2 & \xi z_1       & \beta z_2       & \alpha_2 z_3\\
     \alpha_2 & -\beta w \overline{z}_2  & \xi w \overline{z}_1   & -\alpha_2 z_3\\
     \alpha_3 & \ast    & \ast   & \ast\\
\end{pNiceMatrix}=
\begin{pNiceArray}{c}
     u_1 \\
     u_2 \\
     u_3 \\
     u_4
\end{pNiceArray}$$}

Using the orthogonality of the rows, we obtain
\begin{align*}
    \alpha_1 \alpha_2+\alpha_2 (\xi z_1+\beta z_2)+ \alpha_2 \alpha_3 z_3 &= \langle u_2, u_1\rangle = 0\\
    \alpha_1 \alpha_2+\alpha_2 w(\xi\overline{z}_1-\beta\overline{z}_2)-\alpha_2 \alpha_3 z_3 &= \langle u_3, u_1\rangle = 0.
\end{align*}
Dividing the equations by $\alpha_2$ and multiplying the second one by $\overline{w}$ we obtain
\begin{align*}
    \alpha_1+ \xi z_1+\beta z_2+  \alpha_3 z_3 &=  0,\\
    \alpha_1 \overline{w}+ \xi \overline{z}_1-\beta\overline{z}_2- \alpha_3 \overline{w} z_3 &=  0.
\end{align*}
By taking the complex conjugate of the second equation we obtain
\begin{align*}
   \xi z_1+\beta z_2 &= - (\alpha_1+\alpha_3 z_3),\\
   \xi z_1-\beta z_2 &=  -(\alpha_1 w-\alpha_3 w \overline{z}_3).
\end{align*}
We can solve for $z_1$ and $z_2$, namely:
\begin{align*}
    z_1&=-\frac{\alpha_1(1+w)+\alpha_3(z_3-w \overline{z}_3)}{2\xi},\\
    z_2&=-\frac{\alpha_1(1-w)+\alpha_3(z_3+w \overline{z}_3)}{2\beta}.
\end{align*}
Note that the equations above are equivalent to $\langle u_2,u_1\rangle=\langle u_3,u_1\rangle=0$. Now, let $z_3 = e^{it}$ and $w=e^{is}$ and determine for which pairs $(t,s)$ the above formulas give $|z_1|=|z_2|=1$. We have
\begin{align*}
    |z_1|^2&=\frac{\alpha_1^2+\alpha_3^2+\alpha_1^2 \cos(s)-\alpha_3^2\cos(s-2t)}{2 \xi^2}\\
    |z_2|^2&=\frac{\alpha_1^2+\alpha_3^2-\alpha_1^2 \cos(s)+\alpha_3^2\cos(s-2t)}{2 \beta^2}.
\end{align*}
From (\ref{ucoefeq}) we observe that $\beta^2+\xi^2=\alpha_1^2+\alpha_3^2$ and therefore
\begin{align*}
    |z_1|^2&=\frac{1}{2}+\frac{\beta^2+\alpha_1^2 \cos(s)-\alpha_3^2\cos(s-2t)}{2 \xi^2}\\
    |z_2|^2&=\frac{1}{2}+\frac{\xi^2-\alpha_1^2 \cos(s)+\alpha_3^2\cos(s-2t)}{2 \beta^2}.
\end{align*}
This means that $|z_1|=|z_2|=1$ precisely when $\beta^2=\xi^2-\alpha_1^2 \cos(s)+\alpha_3^2 \cos(s-2t)$. Therefore $$t=\frac{1}{2}\left(s-\arccos\left(\frac{\beta^2-\xi^2+\alpha^2_1\cos(s)}{\alpha_3^2}\right)\right).$$
Since $\alpha_1^2+\alpha_3^2=\xi^2+\beta^2$, we always have 
$$-\alpha_3^2\leq \beta^2-\xi^2+\alpha_1^2\cos(s)\leq \alpha_3^2.$$
This shows there is a 1-parameter family of solutions for which $u$ is a unitary. 

Since every graph in the commuting square has simple edges and the first row of $u$ is the same for all of these connections, we conclude that they are all non-equivalent connections. Below we have a table showing the indices for the subfactors obtained from a commuting square associated to $S(i,i,j,j)$.

\begin{table}[H]
    \centering
\begin{tblr}[caption = {Indices of $S(i,i,j,j)$-commuting square subfactors}]{|c|c|c|c|c|c|c|c|c|}
    \hline
    \diagbox{$j$}{$i$}&1&2&3&4&$\cdots$&$\infty$\\
    \hline
    1&4& & & & & \\
    \hline
    2&$\frac{5+\sqrt{17}}{2}$&$5$& & & & \\
    \hline
    3&$3+\sqrt{3}$&$5.1249$&$3+\sqrt{5}$& & & \\
    \hline
    4&$\frac{5+\sqrt{21}}{2}$&$5.1642$&$5.2703$&$\frac{7+\sqrt{13}}{2}$& & \\
    \hline
    $\vdots$&$\vdots$&$\vdots$&$\vdots$&$\vdots$&$\ddots$& \\
    \hline
    $\infty$&$2+2\sqrt{2}$&$5.1844$&$5.2870$&$5.3184$& &$\frac{16}{3}$\\
    \hline
\end{tblr}
\caption{Indices of $S(i,i,j,j)$-commuting square subfactors}
\label{tab:indices4stars}
\end{table}
We will use this 1-parameter family of bi-unitary connections to show the 
existence of $A_\infty$-subfactors at the remaining Jones indices between 
$4$ and $5$ at which there exist finite depth subfactors. This will
follow from a cardinality argument that uses in a crucial way a result
of Kawahigashi that we describe next.

Let $N\sst M$ be a hyperfinite subfactor with finite index and finite depth. 
The following proposition and corollary can be found in 
\cite[Corollary 3.6]{kawahigashi2023characterization}:

\begin{prop}\label{kawa}
    Any finite-dimensional commuting squares $\{A_{0,l}\sst A_{1,l}\}_l$ giving a hyperfinite subfactor $A_{0,\infty}\sst A_{1,\infty}$ that is isomorphic to $N\sst M$ is of the form
    $$A_{k,l}=\begin{cases}
    \End(_S S_1\otimes_S \cdots\otimes_S S_1\otimes_S X\otimes_P Q\otimes_P\cdots\otimes_P Q_P),& \\
    \text{$k/2$ copies of $S_1$ and $l/2$ copies of $Q$,}& \text{if $k$ and $l$ are even,}\\
    \End(_{S_1} S_1\otimes_S \cdots\otimes_S S_1\otimes_S X\otimes_P Q\otimes_P\cdots\otimes_P Q_P),& \\
    \text{$(k+1)/2$ copies of $S_1$ and $l/2$ copies of $Q$,}& \text{if $k$ is odd and $l$ is even,}\\
    \End(_S S_1\otimes_S \cdots\otimes_S S_1\otimes_S X\otimes_P Q\otimes_P\cdots\otimes_P Q_Q),&\\
    \text{$k/2$ copies of $S_1$ and $(l+1)/2$ copies of $Q$,}& \text{if $k$ is even and $l$ is odd,}\\
    \End(_{S_1} S_1\otimes_S \cdots\otimes_S S_1\otimes_S X\otimes_P Q\otimes_P\cdots\otimes_P Q_Q),&\\
    \text{$(k+1)/2$ copies of $S_1$ and $(l+1)/2$ copies of $Q$,}& \text{if $k$ and $l$ are odd,}
    \end{cases}$$
    where $P\sst Q$ is a subfactor Morita equivalent to $R=M^{\opp} \sst N^{\opp}=S$, $S_1$ is the basic construction for $R\sst S$ and $_S X_P$ is the $S-P$ bimodule giving the Morita equivalence between the category of $S-S$ bimodules arising from $S\sst S_1$ and $P-P$ bimodules arising from $P\sst Q$.
\end{prop}

\begin{coro}\label{kawacoro}
    There are only countably many non-equivalent commuting squares $\{A_{i,j},\, 0\leq i,j\leq 1\}$ such that $A_{0,\infty}\sst A_{1,\infty}$ is isomorphic to $N\sst M$.
\end{coro}

\begin{proof}
    Given a fusion category $\CC$, there are only finitely many fusion categories Morita equivalent to $\CC$ (see \cite{longo1994duality}, Section 6). Let $R\sst S$ be as in Proposition \ref{kawa}. Fix a fusion category of $P-P$ bimodules Morita equivalent to the fusion category of $S$-$S$ bimodules arising from $S\sst S_1$. Since there are countably many subfactors $P\sst Q$ and countably many bimodules $_S X_P$, by Proposition \ref{kawa}, we have countably many non-equivalent commuting squares that produce $N\sst M$.    
\end{proof}

Since we have constructed above are 1-parameter families of non-equivalent 
commuting squares, we obtain

\begin{theo}\label{1paramAinf}
    There are irreducible, hyperfinite $A_\infty$-subfactors with index 
$\frac{5+\sqrt{17}}{2}$, $3+\sqrt{3}$, $\frac{5+\sqrt{21}}{2}$ and $5$. 
There is an irreducible, hyperfinite subfactor with index $3+\sqrt{5}$ 
that either has principal graph $A_\infty$ or has an $A_3\ast A_4$ 
standard invariant \cite{bisch1997intermediate}.
\end{theo}
\begin{proof}
 By classification of small index subfactors, we have finitely many 
finite depth subfactors at the indices $\frac{5+\sqrt{17}}{2}$, 
$3+\sqrt{3}$, $\frac{5+\sqrt{21}}{2}$, $5$ and $3+\sqrt{5}$. 
Corollary \ref{kawacoro} then implies that our 1-parameter families 
of non-equivalent bi-unitary connections produce at least one irreducible, 
hyperfinite, infinite depth subfactor at each of these indices. Using 
the classification of small index subfactor planar algebras again, we 
have that all of them but the one with index $3+\sqrt{5}$ must have 
trivial standard invariant. In the case of index $3+\sqrt{5}$, since 
the commuting square subfactor is irreducible, by classification it 
has to be either an $A_\infty$-subfactor or an $A_3\ast A_4$ Fuss-Catalan
subfactor. The latter exists as a hyperfinite subfactor by 
\cite{bisch1997intermediate} (construction of the standard invariant)
and \cite{popa1995freeindependent} (realizability as {\it hyperfinite}
subfactor).
\end{proof}

Note that we have constructed $A_\infty$-subfactors with index
$\frac{5+\sqrt{17}}{2}$ in two different ways in sections 4 and 5.
It is natural to ask if they are isomorphic. We would expect this
not to be the case, but do not currently have the tools to decide
this question. 

The 1-parameter families of bi-unitary connections constructed on
$4$-stars above may lead to many non-isomorphic A$_\infty$-subfactors
at the same index. The
cardinality argument we use does not allow us to determine what
happens at different parameters. This is a very intriguing open 
problem.




    \section{Intermediate commuting squares}
Consider a symmetric commuting square of finite-dimensional C*-algebras as in \ref{csfindim}: 
	\begin{equation}\label{csq2}
	\begin{array}{ccc}
	A_{1,0} & \sst & A_{1,1} \\
	\cup &  & \cup \\
	A_{0,0} & \sst & A_{0,1}
	\end{array}
	\end{equation}
 and let $N=A_{0,\infty}\sst A_{1,\infty}=M$ be the subfactor we obtain by iterating the basic construction horizontally as in \ref{infgrid}. Suppose this subfactor has an intermediate subfactor $P$, that is $N\sst P\sst  M$ and let $p_1$ be Jones projection for $P\sst M$. Note that since $M'\sst P'\sst N'$ in $B(L^2(M))$ then $M\sst P_1\sst M_1$, where $P_1$ and $M_1$ are the result of the basic construction for $P\sst M$ and $N\sst M$ respectively. Thus, we have
 $$p_1\in P'\cap P_1\sst N'\cap M_1=A_{0,1}'\cap A_{2,0}\sst A_{2,0}.$$
 Let $B_{1,n}=P\cap A_{1,n}=\{p_1\}'\cap A_{1,n}$. We will show that 
 \begin{equation}\label{csqint}
	\begin{array}{ccc}
	B_{1,0} & \sst & B_{1,1} \\
	\cup &  & \cup \\
	A_{0,0} & \sst & A_{0,1}
	\end{array},\quad \begin{array}{ccc}
	A_{1,0} & \sst & A_{1,1} \\
	\cup &  & \cup \\
	B_{1,0} & \sst & B_{1,1}
	\end{array}
\end{equation}
are nondegenerate commuting squares. 

\begin{lem}\label{lemcs1}
    If $\begin{array}{ccc}
	C & \sst & D \\
	\cup &  & \cup \\
	A & \sst & B
	\end{array}$ is a commuting square of von Neumann algebras with 
respect to a fixed normal, faithful trace $\tr$ on D, and $R\sst S$ are 
von Neumann algebras such that $A\sst R\sst C$ and $B\sst S\sst D$, 
then 

$$\begin{array}{ccc}
	R & \sst & S \\
	\cup &  & \cup \\
	A & \sst & B
	\end{array}$$

is also a commuting square with respect to the trace $\tr$ restricted 
to the subalgebras.
\end{lem}
\begin{proof}
This is obvious as $L^2(A)^\perp \cap L^2(B) \perp L^2(A)^\perp \cap L^2(C)$,
implies $L^2(A)^\perp \cap L^2(B) \perp L^2(A)^\perp \cap L^2(R)$.

\end{proof}

Setting $A=A_{0,0},B=A_{0,1},C=A_{1,0},D=A_{1,1},R=B_{1,0}, S=B_{1,1}$ in lemma \ref{lemcs1}, it follows that the first set of inclusions in 
\ref{csqint} is a commuting square.

\begin{lem}\label{lemcs2} For all $n\geq 0$, the quadrilaterals 
$\begin{array}{ccc}
	A_{1,n} & \sst & M \\
	\cup &  & \cup \\
	B_{1,n} & \sst & P
	\end{array}$ are commuting squares.
\end{lem}
\begin{proof}
    Let $x\in A_{1,n}$ then, since $p_1$ implements the conditional expectation $E^M_P$ we have 
    $$E^{M_1}_M(p_1 x p_1)=E^{M_1}_M(E^M_P(x) p_1)=E^M_P(x)E^{M_1}_M(p_1)=E^M_P(x)E^{P_1}_M(p_1)=\beta^{-1} E^M_P(x),$$
    where $\beta=[M:P]$. 
    
    Now, since $\begin{array}{ccc}
	A_{2,n} & \sst & M_1 \\
	\cup &  & \cup \\
	A_{1,n} & \sst & M
	\end{array}$ is also a commuting square, we have $E^{M_1}_M(A_{2,n})\sst A_{1,n}$. Recall that $p_1\in A_{2,0}$, therefore $p_1xp_1\in A_{2,n}$ and consequently $$E^M_P(x)=\beta E^{M_1}_M(p_1 x p_1)\in A_{1,n}\cap P=B_{1,n}.$$
\end{proof}

Setting $A=B_{1,0},B=A_{1,0},C=P,D=M,R=B_{1,1}, S=A_{1,1}$ in lemma \ref{lemcs1}, it follows that the second set of inclusions in \ref{csqint} is a 
commuting square.

Since \ref{csq2} is nondegenerate we have $A_{1,1}=\spn A_{1,0}A_{0,1}$ and consequently
$$B_{1,1}=E^{A_{1,1}}_{B_{1,1}}(A_{1,1})=\spn E^{A_{1,1}}_{B_{1,1}}(A_{1,0}A_{0,1})=\spn E^{A_{1,1}}_{B_{1,1}}(A_{1,0})A_{0,1}\sst \spn B_{1,0}A_{0,1}$$
which shows that the first commuting square in \ref{csqint} is nondegenerate. Similarly,
$$A_{1,1}=\spn A_{1,0}A_{0,1}\sst \spn A_{1,0}B_{1,1}$$
shows that the second commuting square is nondegenerate. Thus we have proved:

\begin{theo}\label{intercs}
        Let $N\sst P\sst M$ be an intermediate subfactor of a commuting
square subfactor $N \subset M$ constructed from a nondegenerate commuting
square \ref{csq2}. Set $B_{1,n}=A_{1,n}\cap P$. Then 
        $$\begin{array}{cccccc}
    B_{1,0}&\sst &B_{1,1}&\sst &B_{1,2}&\sst \cdots\\
    \rotatebox[origin=c]{90}{$\sst$}_{G_{0,0}}& &\rotatebox[origin=c]{90}{$\sst$}_{L_{0,1}}& &\rotatebox[origin=c]{90}{$\sst$}_{G_{0,2}}&\\
    A_{0,0}&\stackrel{\text{H}}{\sst} &A_{0,1}&\stackrel{\text{H}^t}{\sst} &A_{0,2}&\sst \cdots
    \end{array},\quad \begin{array}{cccccc}
    A_{1,0}&\stackrel{K}{\sst} &A_{1,1}&\stackrel{K^t}{\sst}  &A_{1,2}&\sst \cdots\\
    \rotatebox[origin=c]{90}{$\sst$}_{G_{1,0}}& &\rotatebox[origin=c]{90}{$\sst$}_{L_{1,1}}& &\rotatebox[origin=c]{90}{$\sst$}_{G_{1,2}}&\\
    B_{1,0}&\sst &B_{1,1}&\sst &B_{1,2}&\sst \cdots
    \end{array}$$
    are sequences of nondegenerate commuting squares approximating $N\sst P$ and $P\sst M$, respectively.
    \end{theo}

Note that in general $G_{i,n}$ and $L_{i,n}$, $i=0$, $1$, may not be connected 
inclusion graphs (see e.g. \cite{bisch1997intermediate}), and consequently 
we need a more general version of the index formula to compute the Jones 
indices of the intermediate subfactor as given 
in \cite[Theorem 4.3.3]{goodman2012coxeter}. It turns out that the index formula as in the case of connected inclusion graphs still holds provided the commuting square is nondegenerate.

\begin{lem}\label{indformula}
	Let $\begin{array}{rcl}
		N&\sst& M\\
		\cup& & \cup\\
		A&\stackrel{G}{\sst} &B
	\end{array}$ be a nondegenerate commuting square, where $A\sst B$ are multi-matrix algebras with inclusion graph $G$ (possibly not connected) and $N\sst M$ is a finite index inclusion of II$_1$ factors. Then $[M:N]=\|G\|^2$ and $G^tG \vec{t}=[M:N] \vec{t}$, where $\vec{t}$ is the trace vector for $B$ induced
from the normalized, normal faithful trace on $M$. In particular, 
$G^tG$ has to have a positive eigenvector.
\end{lem}

\begin{proof}
	Since we have a nondegenerate commuting square, by performing the basic construction horizontally, we obtain the following commuting square
	$$\begin{array}{rcl}
		M&\sst& \langle M, e \rangle=M_1\\
		\cup& & \cup\\
		B&\sst &\langle B, e \rangle=B_1\end{array}$$
	where we can identify the Jones projections $e_N$ and $e_A$ for both
horizontal inclusions and simply denote it by $e$. Let $\tr$ be unique 
normalized trace on $M_1$ and $E_B$ the conditional expectation from $B_1$ onto $B$ with respect to $\tr|_{B_1}$. Observe that 
	$$E_B(e_A)=E_B(e_N)=E_{B_1}E_{M}(e_N)=\lambda 1,$$
	where $\lambda=[M:N]^{-1}$. Hence $\tr|_B$ is a $\lambda$-Markov 
trace and consequently 
	$G^t G \vec{t}=\lambda^{-1}\vec{t}$, where $\vec{t}$ is the 
trace vector on $B$. 
	
	Note that since $\vec{t}$ is a positive vector (i.e. all entries are
stricly positive), by \cite[Corollary 8.1.30]{horn2012matrix} we obtain 
that $\vec{t}$ has to be the eigenvector associated to $\|G^tG\|=\|G\|^2$.
\end{proof}

In theorem \ref{1paramAinf} we constructed an irreducible hyperfinite subfactor $N\sst M$ with index $3+\sqrt{5}$ using a commuting square whose first vertical inclusion graph is $G=S(3,3,3,3)$. If this subfactor is isomorphic to an 
$A_3\ast A_4$ Fuss-Catalan subfactor, then it has an intermediate 
subfactor $N\sst P\sst M$ such that $[P:N]=\frac{3+\sqrt{5}}{2}$ (or $2$,
but we may assume $[P:N]=\frac{3+\sqrt{5}}{2}$ by considering the upwards basic construction with vertical inclusions $G^t$). 
In this case, by theorem \ref{intercs} and lemma \ref{indformula}, there exist 
$G_1$, $G_2$ such that $G=G_1 G_2$, $\|G_1\|^2=\frac{3+\sqrt{5}}{2}$, $\|G_2\|^2=2$, and both $G_1$ and $G_2$ have positive eigenvectors.

\begin{prop}
	Let $G\in M_{m\times n}(\Z_{\geq 0})$ be the adjacency matrix of 
a graph, also called $G$. Then there are only finitely many ways to write 
$G=HK$, where $H$ and $K$ are adjacency matrices for bipartite graphs without 
isolated vertices.
\end{prop}
\begin{proof}
	Let $H\in M_{m\times q}(\Z_{\geq 0})$ and $K\in M_{q\times n}(\Z_{\geq 0})$ be adjacency matrices corresponding to bipartite graphs without isolated vertices such that $G=HK$. Since $H$ and $K$ correspond to graphs without isolated vertices we have $\sum_{i=1}^m H_{ik}\geq 1$ and $\sum_{j=1}^n K_{kj}\geq 1$ for any $(i,j)$. Note that 
	$$S=\sum_{i=1}^m\sum_{j=1}^n G_{ij}=\sum_{i=1}^m\sum_{j=1}^n\sum_{k=1}^q H_{ik}K_{kj}=\sum_{k=1}^q\sum_{i=1}^m H_{ik}\sum_{j=1}^n K_{kj}\geq q,$$
	hence $H$ and $K$ cannot be arbitrarily large. 
	
	Since $\sum_{j=1}^n K_{kj}\geq 1$, for any $k$ there exists $j$ such that $K_{kj}\geq 1$. Observe that $H_{i,k}K_{kj} \leq G_{ij}$ and therefore
	$$H_{i,k} \leq \frac{G_{ij}}{K_{kj}}\leq  G_{ij}\leq S,$$
	hence the entries of $H$ are bounded. Similarly, we conclude that the entries of $K$ are bounded.
\end{proof}

The previous proposition implies that there are only finitely many 
$G_1$, $G_2$ such that $G=G_1 G_2$. To determine the possible $G_1$'s and 
$G_2$'s we will decompose $G$ into \emph{building blocks}. A building block 
is a bipartite graph where all vertices have degree 1 except for a single 
vertex with degree 2. To illustrate this idea, here is an example of 
a building block:

\begin{center}
    \begin{tikzpicture}[main/.style={circle,fill=black,inner sep=0pt,minimum size=1.5mm},sub/.style={circle,fill=white,draw=black,inner sep=0pt,minimum size=1.5mm},end/.style={circle,diagonal fill={black}{white},draw=black,inner sep=0pt,minimum size=1.5mm},baseline=0,scale=0.6]
		\node[main] at (0,0) (A3) {};
		\node[main] at (1.5,0) (A4) {};
		\node[main] at (2.5,0) (A5) {};
		\node[main] at (-1.5,0) (A2) {};
		\node[main] at (-2.5,0) (A1) {};
		\node[main] at (0.5,1) (B3) {};
		\node[main] at (-0.5,1) (B4) {};
		\node[main] at (1.5,1) (B5) {};
		\node[main] at (2.5,1) (B6) {};
		\node[main] at (-1.5,1) (B2) {};
		\node[main] at (-2.5,1) (B1) {};
		\draw (B1) -- (A1);
		\draw (B2) -- (A2);
		\draw (B3) -- (A3) -- (B4);
		\draw (B5) -- (A4);
		\draw (B6) -- (A5);
	\end{tikzpicture}
\end{center}

We will then express $G$ as a product of building blocks, for example
$$G=\begin{tikzpicture}[main/.style={circle,fill=black,inner sep=0pt,minimum size=1.5mm},baseline=0,scale=0.6]
	\node[main] at (0,-1) (A1) {};
	\node[main] at (0,0) (A2) {};
	\node[main] at (0,1) (A3) {};
	\node[main] at (1,-0.5) (B1) {};
	\node[main] at (1,0.5) (B2) {};
	\draw (A1) -- (B1) -- (A2) -- (B2) -- (A3);
\end{tikzpicture} = \begin{tikzpicture}[main/.style={circle,fill=black,inner sep=0pt,minimum size=1.5mm},baseline=0,scale=0.6]
	\node[main] at (0,-1.5) (A1) {};
	\node[main] at (0,0) (A2) {};
	\node[main] at (0,1.5) (A3) {};
	\node[main] at (1,-1.5) (B1) {};
	\node[main] at (1,-0.5) (B2) {};
	\node[main] at (1,0.5) (B3) {};
	\node[main] at (1,1.5) (B4) {};
	\node[main] at (1.4,-1.5) (C1) {};
	\node[main] at (1.4,-0.5) (C2) {};
	\node[main] at (1.4,0.5) (C3) {};
	\node[main] at (1.4,1.5) (C4) {};
	\node[main] at (2.4,-1.5) (D1) {};
	\node[main] at (2.4,-0.5) (D2) {};
	\node[main] at (2.4,1) (D3) {};
	\node[main] at (2.8,-1.5) (E1) {};
	\node[main] at (2.8,-0.5) (E2) {};
	\node[main] at (2.8,1) (E3) {};
	\node[main] at (3.8,-1) (F1) {};
	\node[main] at (3.8,1) (F2) {};
	\draw (A1) -- (B1);
	\draw (B2) -- (A2) -- (B3);
	\draw (A3) -- (B4);
	\draw (C1) -- (D1);
	\draw (C2) -- (D2);
	\draw (C3) -- (D3) -- (C4);
	\draw (E1) -- (F1) -- (E2);
	\draw (E3) -- (F2);
\end{tikzpicture}=B_1 B_2 B_3.$$
We will call $\mathcal{B}=(B_1,B_2,B_3)$ a \emph{chain} for $G$. Two 
chains $\mathcal{B}=(B_1,\dots,B_k)$ and $\mathcal{B'}=(B'_1,\dots,B'_k)$ 
will be {\it isomorphic} if there exist permutation matrices 
$(P_i)_{i=0}^k$ such that $$B'_i=P_{i-1}B_i P^{-1}_i.$$

To determine all possible $G_1$'s and $G_2$'s, we reduce every maximal 
chain to every possible chain of length 2. For the example above, there are two possibilities: 
$$(G_1,G_2)=(B_1,B_2B_3)\text{ and }(G_1,G_2)=(B_1B_2,B_3).$$

For $G=S(3,3,3,3)$, there are $12240$ maximal chains of length $11$ (these computations were done in \verb|Python| and can be found in  \cite{Caceres_Graph_Factorization}). These yield precisely $35$ non-isomorphic ways to write $G=G_1G_2$, where $G_i$ are adjacency matrices for bipartite graphs, not necessarily connected, without isolated vertices. There are precisely two factorizations for which we have $\|G_1\|^2=2$ and $\|G_2\|^2=\frac{3+\sqrt{5}}{2}$:
\begin{align*}
	G_1=\left[\begin{matrix}1 & 1 & 0 & 0 & 0 & 0 & 0 & 0\\0 & 0 & 1 & 1 & 0 & 0 & 0 & 0\\0 & 0 & 0 & 0 & 1 & 1 & 0 & 0\\0 & 0 & 0 & 0 & 0 & 0 & 1 & 0\\0 & 0 & 0 & 0 & 0 & 0 & 0 & 1\end{matrix}\right],\;& G_2=\left[\begin{matrix}0 & 0 & 1 & 0 & 0 & 0 & 0 & 0\\1 & 0 & 0 & 0 & 0 & 0 & 0 & 0\\0 & 0 & 0 & 1 & 0 & 0 & 0 & 0\\0 & 1 & 0 & 0 & 0 & 0 & 0 & 0\\0 & 0 & 1 & 1 & 0 & 1 & 0 & 0\\0 & 0 & 0 & 0 & 1 & 0 & 0 & 0\\0 & 0 & 0 & 0 & 1 & 0 & 1 & 0\\0 & 0 & 0 & 0 & 0 & 1 & 0 & 1\end{matrix}\right]\\
	G_1=\left[\begin{matrix}1 & 1 & 0 & 0 & 0 & 0 & 0\\0 & 0 & 1 & 0 & 0 & 0 & 0\\0 & 0 & 0 & 1 & 1 & 0 & 0\\0 & 0 & 0 & 0 & 0 & 1 & 0\\0 & 0 & 0 & 0 & 0 & 0 & 1\end{matrix}\right],\;& G_2=\left[\begin{matrix}0 & 0 & 1 & 0 & 0 & 0 & 0 & 0\\1 & 0 & 0 & 0 & 0 & 0 & 0 & 0\\0 & 1 & 0 & 1 & 0 & 0 & 0 & 0\\0 & 0 & 1 & 0 & 1 & 1 & 0 & 0\\0 & 0 & 0 & 1 & 0 & 0 & 0 & 0\\0 & 0 & 0 & 0 & 1 & 0 & 1 & 0\\0 & 0 & 0 & 0 & 0 & 1 & 0 & 1\end{matrix}\right].
\end{align*}
In both cases $G_1^t G_1$ has no positive eigenvector, hence we can conclude:
\begin{theo}
    There is an irreducible, hyperfinite $A_\infty$-subfactor with 
index $3+\sqrt{5}$.
\end{theo}
    
\medskip
The method presented in this section allows us to determine
the standard invariant of our commuting square subfactor without computing
any higher relative commutants explicitly. This would be hard to do, even
in this fairly simple case. We established that our
simple combinatorial method suffices to conclude that the commuting 
square subfactor with index $3+\sqrt{5}$ that we constructed in 
theorem 5.4 must have TLJ standard.

    \begin{landscape}

\appendix

\section{Bi-unitary connections}\label{appconnections}
\subsection{Large broom: unitary $u$ ($d=\sqrt{17}$) - part 1}\hfill\\
\hfill\\
\includegraphics[width=\linewidth,page=15]{tables}
\newpage 
\subsection{Large broom: unitary $u$ ($d=\sqrt{17}$) - part 2}\hfill\\
\hfill\\
\includegraphics[width=\linewidth,page=16]{tables}
\newpage
\subsection{Large broom: unitary $v$ ($d=\sqrt{17}$) - part 1}\hfill\\
\hfill\\
\includegraphics[width=\linewidth,page=17]{tables}
\newpage
\subsection{Large broom: unitary $v$ ($d=\sqrt{17}$) - part 2}\hfill\\
\hfill\\
\includegraphics[width=\linewidth,page=18]{tables}

\end{landscape}

	\nocite{*}
	\bibliographystyle{alpha}
	\bibliography{references}{}	
\end{document}